\documentclass[10pt,a4paper]{amsart}
\usepackage[toc,page]{appendix}
\usepackage{a4}
\usepackage[active]{srcltx}
\usepackage{amsmath,bbm,esint}
\usepackage{amssymb}
\usepackage{amsfonts}
\usepackage{mathrsfs} 
\usepackage{graphicx}
\usepackage{tikz}
\usepackage{color}
\usepackage[colorlinks=true,urlcolor=blue,
citecolor=red,linkcolor=blue,linktocpage,pdfpagelabels,bookmarksnumbered,bookmarksopen]{hyperref}

\usepackage{mathtools}
\mathtoolsset{showonlyrefs}

\newcounter{hyp}
\newcounter{hyp_prime}
\def\a{\alpha}
\def\b{\beta}
\def\d{\delta}

\def\l{\lambda}

\def\t{\tau}

\def\Om{\Omega}

\newcommand{\cF}{{\mathcal F}}

\newcommand{\cL}{{\mathcal L}}

\newcommand{\Z}{{\mathbb Z}}
\newcommand{\R}{{\mathbb R}}

\newcommand{\T}{\mathbb{T}}

\newcommand{\E}{{\mathbb E}}

\renewcommand{\P}{{\mathbb P}}

\newcommand{\sP}{{\mathscr P}}

\newcommand{\ds}{\displaystyle}

\newcommand{\spt}{{\rm{spt}}}

\newcommand{\dd}{\hspace{0.7pt}{\rm d}}

\newcommand{\sym}{{\rm sym}}

\DeclareMathOperator*{\esssup}{ess\,sup}

\usepackage{color}

\newcommand{\be}{\begin{equation}}
\newcommand{\ee}{\end{equation}}

\newcommand{\ba}{\begin{array}}
\newcommand{\ea}{\end{array}}

\newtheorem{theorem}{Theorem}[section]
\newtheorem{definition}[theorem]{Definition}
\newtheorem{lemma}[theorem]{Lemma}
\newtheorem{corollary}[theorem]{Corollary}
\newtheorem{proposition}[theorem]{Proposition}

\newtheorem{remark}[theorem]{Remark}

\newtheorem{innercustomgeneric}{\customgenericname}
\providecommand{\customgenericname}{}
\newcommand{\newcustomtheorem}[2]{%
  \newenvironment{#1}[1]
  {%
   \renewcommand\customgenericname{#2}%
   \renewcommand\theinnercustomgeneric{##1}%
   \innercustomgeneric
  }
  {\endinnercustomgeneric}
}

\newcustomtheorem{customthm}{Framework}

\numberwithin{equation}{section}

\title[]{Long Time Behavior and Stabilization for Displacement Monotone Mean Field Games}  

\author[M. Cirant]{Marco Cirant} 
\address{Dipartimento di Matematica ``Tullio Levi-Civita'', Universit\`a di Padova, Via Trieste 63, 35121, Italy.}
\email{cirant@math.unipd.it}

\author[A.R. M\'esz\'aros]{Alp\'ar R. M\'esz\'aros}  
\date{\today}
\address{Department of Mathematical Sciences, University of Durham, Durham DH1 3LE, United Kingdom}
\email{alpar.r.meszaros@durham.ac.uk} 

\begin{document}

\begin{abstract} This paper is devoted to the study of the long time behavior of Nash equilibria in Mean Field Games within the framework of displacement monotonicity. We first show that any two equilibria defined on the time horizon $[0,T]$ must be close as $T \to \infty$, in a suitable sense, independently of initial/terminal conditions. The way this stability property is made quantitative involves the $L^2$ distance between solutions of the associated Pontryagin system of FBSDEs that characterizes the equilibria. Therefore, this implies in particular the stability in the 2-Wasserstein distance for the two flows of probability measures describing the agent population density and the $L^2$ distance between the co-states of agents, that are related to the optimal feedback controls. We then prove that the value function of a typical agent converges as $T \to \infty$, and we describe this limit via an infinite horizon MFG system, involving an ergodic constant. All of our convergence results hold true in a unified way for deterministic and idiosyncratic noise driven Mean Field Games, in the case of strongly displacement monotone non-separable Hamiltonians. All these are quantitative at exponential rates.

\end{abstract}

\maketitle

\section{Introduction}

Since the inception of Mean Field Games (MFGs) theory, these models have been explored in contexts involving finite time horizons, infinite time horizons, and stationary (or time-independent) frameworks. Already the foundational papers in the field highlight the distinct nature of each approach (cf. \cite{LasLio:06-I, LasLio:06-II,HuaMalCai}). Each model type is motivated by specific applications, and the interplay between them has raised intriguing questions. One such question is concerned with the long time behavior of solutions to finite horizon MFGs and in particular whether some sort of `stabilization' phenomena could be observed as the time horizon tends to infinity.  

Before the conception of MFGs, similar questions have been studied intensively in the context of classical control problems and the associated Hamilton--Jacobi--Bellman (HJB) equations. Probably the first account on this matter is documented in \cite{LioPapVar}, which later turned out to be influential in many profound directions in fields as homogenization, the weak KAM theory  \`a la Fathi, the Aubry--Mather theory and elsewhere. Without the intention of being exhaustive, we refer to the classical works \cite{Fat:97, BarSou, FatMad} and to the monographs and lecture notes \cite{Fat:08, Fat:12,  BarCap, MitTra} for an excellent exposition of some of these directions.

At the same time, the study of the long time asymptotic properties and convergence to equilibrium of solutions to Fokker--Planck--Kolmogorov (FPK) type or more general parabolic PDEs received great attention in the literature. It would simply be impossible to mention all the various approaches and the vast amount of deep results in this context. However, we would like to point to \cite{BGG}, to the recent work \cite{Por:24} and to the references therein for a description of some of these results. The latter work has a particular connection to our results, as this can be seen as a bridge between the FPK and HJB worlds.

Returning to MFGs, in the past one and half decades many authors have contributed to the understanding of asymptotic behavior of MFG systems and the corresponding master equations. The main actor of this manuscript is the MFG system 

\begin{equation}\label{eq:MFG_intro}
\left\{
\begin{array}{ll}
-\partial_t u  -\b\Delta u + H(x,-D_{x}u,\rho) = 0, & \text{in } (0,T)\times \R^d,\\[5pt]
\partial_t\rho -\b\Delta \rho + \nabla\cdot (\rho D_pH(x,-D_{x}u,\rho)) = 0, & \text{in } (0,T)\times \R^d,\\[5pt]
\rho(0,\cdot) = \rho_0; \ u(T,\cdot) = g(\cdot,\rho_T), & \text{in } \R^d.
\end{array}
\right.
\end{equation}
Here the data consists of the Hamiltonian $H:\R^d\times\R^d\times\sP_2(\R^d)\to\R$, the final cost function $g:\R^d\times\sP_2(\R^d)\to\R$, the initial agent distribution $\rho_0\in\sP_2(\R^d)$, the time horizon $T>0$ and the intensity of the idiosyncratic noise $\b\ge 0$ which is allowed to be $0$ throughout the text. A solution to this system is a pair, given by the value function $u:[0,T]\times\R^d\to\R$ of any typical agent, and $(\rho_s)_{s\in[0,T]}$ the actual MFG Nash equilibrium, a flow in $\sP_2(\R^d)$ representing the evolution of the distribution of the agents. To emphasize the dependence on the time horizon $T>0$, we sometimes refer to $(u,\rho)$ as $(u^T,\rho^T)$.

Our goal in this paper is to present a class of sufficient assumptions on the data $H$ and $g$ which allow the study of the asymptotic behavior of $(u^T,\rho^T)$ as $T\to+\infty$, in a suitable sense. Our aim is to obtain results which are independent of $\b\ge 0$, allow a general class of Hamiltonians, which are {\it non-separable} and are valid for general $\rho_0\in\sP_2(\R^d)$.

\medskip

\noindent {\bf A brief history on the asymptotic behavior of $(u^T,\rho^T)$ as $T\to+\infty$.}

\medskip

This direction and some of the corresponding ideas were first mentioned by Lions in \cite{Lions}. A relatively up to date account on the development of this line of research is given in \cite[Section 1.3.6]{CarPor20}. In what follows we briefly describe this evolution in various contexts. In \cite{GomMohSou} the limiting properties were studied for discrete time discrete space MFG. In \cite{CarLasLioPor:12} and \cite{CarLasLioPor:13} the authors studied the long time convergence problem in the case of purely quadratic and separable Hamiltonians, uniformly parabolic setting (i.e. $\b>0$) both in the case of Hamiltonians depending locally and non-locally on the measure variable. \cite{Car} was concerned with the limiting behavior of first order nonlocal models in connection with weak KAM theory, while in \cite{CarGra} the authors relied on variational techniques for first order MFG involving Hamiltonians which are local in the measure variable. Similar results to those in \cite{CarLasLioPor:13} were developed in the context of the master equation in the work \cite{CarPor:19}.

All these results described above used crucially two facts: (i) the separability of the Hamiltonian, i.e. that  $H$ decomposes as 
\begin{align}\label{def:sep}
H(x,p,\rho) = H_0(x,p) - f(x,\rho),
\end{align}
for some $H_0:\R^d\times\R^d\to\R$ and $f:\R^d\times\sP_2(\R^d)\to\R$; and (ii)  the convexity of $H_0$ with respect to $p$ and the (strong/strict) {Lasry--Lions monotonicity} (LL-monotonicity) of the coupling functions $f$ and $g$. This means that 
\begin{equation}\label{LL-mon}
\int_{\R^d}\left[f(x,\rho_1) - f(x,\rho_2) \right]\dd(\rho_1-\rho_2)(x)\ge c\int_{\R^d} \left[f(x,\rho_1)-f(x,\rho_2)\right]^2\dd x, \  \forall \rho_{1},\rho_{2}\in\sP_{2}(\R^{d}),\tag{LL$_c$}
\end{equation}
for some $c\ge0$. For similar monotonicity conditions we refer also to \cite{Por:18} and \cite{CirPor}. In this latter reference, the authors in particular were able to allow a bit of loss of LL-monotonicity on the price of increasing the noise intensity and show \textit{exponential turnpike} type property for solutions.

Such exponential turnpike property typically reads as follows. Let $(u^T,\rho^T)$ be the solution to \eqref{eq:MFG_intro}. Then there exist $C>0$ and $\omega>0$ (typically depending on the data, and monotonicity constants) such that
\begin{align}\label{ineq:turnpike1}
\|\rho^T(t) - \bar\rho\|_{\mathcal X} + \|D_{x}u^T(t) - D_{x}\bar u \|_{\mathcal Y} \le C\left[e^{-\omega t} + e^{-\omega(T-t)}\right],\ \ \forall t\in (1,T-1),
\end{align}
where $(\bar u,\bar\rho,\bar\l)$ is the solution to the stationary MFG system
\begin{equation}\label{eq:MFG_st}
\left\{
\begin{array}{ll}
\bar\l  -\b\Delta \bar u + H(x,-D_{x}\bar u,\bar\rho) = 0, & \text{in } \R^d,\\[5pt]
 -\b\Delta \bar \rho + \nabla\cdot (\bar \rho D_pH(x,-D_{x}\bar u,\bar \rho)) = 0, & \text{in } \R^d,\\[5pt]
\bar\rho\ge 0; \ \int_{\T^d}\bar u\dd x=0;\ \ \int_{\T^d}\bar\rho\dd x=1, & \text{in } \R^d.
\end{array}
\right.
\end{equation}
Depending on the concrete settings, the function spaces $\mathcal X$ and $\mathcal Y$ have to be chosen suitably. These typically range mostly from $L^p$ or $L^\infty$ to $C^{k,\alpha}$ type spaces in the literature. Thus, stability manifests itself as the presence of a stationary (ergodic) state that attracts finite-horizon equilibria as $T \to \infty$. Another system that has been used to describe the long time asymptotic properties of \eqref{eq:MFG_intro} is the infinite horizon system
\begin{equation}\label{eq:MFG_infinite}
\left\{
\begin{array}{ll}
-\partial_t \tilde u +\tilde\l -\b\Delta \tilde u + H(x,-D_{x}\tilde u,\tilde\rho) = 0, & \text{in } (0,+\infty)\times \R^d,\\[5pt]
\partial_t\tilde\rho -\b\Delta \tilde\rho + \nabla\cdot (\tilde\rho D_pH(x,-D_{x}\tilde u,\tilde\rho)) = 0, & \text{in } (0,+\infty)\times \R^d,\\[5pt]
\tilde \rho(0,\cdot) = \rho_0, & \text{in } \R^d.
\end{array}
\right.
\end{equation}
With respect to the stationary system above, this one incorporates more information: while the former merely describes the stationary state, the latter also clarifies how equilibria starting from any initial state evolve into the stationary one (provided that one can prove that $\tilde u(t)$ approaches $\bar u$ as $t \to \infty$).
In fact, \cite{CirPor} shows that under their standing assumptions, there exists $(\tilde u,\tilde\rho,\tilde\l)$ solution to this system such that 
$$
u^{T}(t,x) - \tilde\l(T-t)\to \tilde u(t,x);\ \ \rho^{T}(t,x)\to\tilde\rho(t,x),\ \ {\rm{as}\ } T\to+\infty,
$$
locally uniformly in $(t,x)$,  and $\tilde u(t)$ itself converges to $\bar u$ as $t \to \infty$.

Departing completely from the LL-monotone regime, the recent paper \cite{BarKou} shows asymptotic characterization for first order MFGs in lack of LL-monotonicity (for a specific class of quadratic Hamiltonians), under particular assumption on the minima of the associated cost functions. Furthermore, in this direction in \cite{CarMas,Mas} the authors obtained weak KAM type results in the context of potential second order MFGs, for separable nonlocal general class of Hamiltonians. These results can in some sense be seen as the second order versions of the results from \cite{GanTud:14} and \cite{GomNur}.

These references mentioned above (with the exception of \cite{BarKou}) consider always data functions which are $\Z^d$-periodic, and hence rely on the compactness of the space $\sP(\T^d)$. In \cite{BarKou}, even though set on the whole space $\R^d$, there seem to be a hidden compactness argument, which comes from the assumption on the location of minima of the cost functions.  There are further other interesting results on the long time behavior, asymptotic analysis on discounted MFGs, and different applications on these, see for instance \cite{PorRic,MimMun,CarTanZha, BayZha}. 

\medskip

When it comes to the long time asymptotic analysis of MFGs with $\R^{d}$ as a state space, the literature is sparse. The very recent manuscript \cite{CecConDurEic} considers the long time behavior of solutions to MFGs genuinely set on $\R^{d}$. This work relies mainly on probabilistic techniques (via the so-called `coupling approach') to obtain exponential turnpike properties. The standing assumptions therein are: a weak form of asymptotic monotonicity on the drift of the controlled dynamics and regularity and smallness conditions on the interaction terms. The main results from this paper are for separable Hamiltonians, and are in the spirit of \eqref{ineq:turnpike1}. These read informally as
\begin{align}\label{ineq:turnpike2}
W_{1}\left(\rho^T(t), \bar\rho\right) + \|D_{x}u^T(t) - D_{x}\bar u \|_{L^{\infty}} \le C\left[e^{-\omega t} + e^{-\omega(T-t)}\right],\ \ \forall t\in (0,T-1),
\end{align}
where $(\bar u,\bar\rho,\bar\l)$ is the solution to the stationary MFG system \eqref{eq:MFG_st}. It is worth mentioning that this work is also completely outside of the LL-monotone regime, it allows a general uniformly elliptic smooth state dependent diffusion matrix, and under further assumptions on the data sometimes gives exponential turnpike properties also for the $D^{2}u^{T}$. 

\medskip
We finally mention that the general study of the long time behavior of MFG lacking of any monotone structure is a rather wide open field of research. In particular, in presence of a genuine multiplicity of stationary states, or more complicated dynamic patterns such as periodic solutions or traveling waves, only few specific models have been so far addressed (\cite{C19, CC21, CC24, GMP, KMfRb, PR}), and stability/instability properties of those patterns are just partially understood even in these special cases.

\medskip

\noindent {\bf Our contributions and the description of our main results.}

\medskip

As highlighted above, the literature on the long-time behavior of solutions to MFGs remains limited in the absence of LL-monotonicity (or semi-monotonicity) conditions on the data, particularly outside the compact regime of $\sP(\T^{d})$. In this manuscript, we address this gap by imposing {\it displacement monotonicity} (D-monotonicity; {see Assumptions \eqref{hyp:Hdis_strong} and \eqref{hyp:g_D-mon} below}) conditions on the data in the setting of the non-compact state space $\R^{d}$. More specifically, we impose \textit{strong} D-monotonicity assumptions, as long time stability fails in general if one requires only D-monotonicity (see Remark \ref{rmk:nonstrict}). Since D-monotonicity is generally incompatible with LL-monotonicity, our results significantly advance the understanding of the long-time asymptotic behavior of solutions beyond the LL-monotone regime in the non-compact setting of $\sP_{2}(\R^{d})$.

To recall, D-monotonicity turned out to be an instrumental sufficient condition for obtaining global in time well-posedness results for MFGs and the corresponding master equations beyond of the LL-monotone regime. This condition is particularly versatile, accommodating non-separable Hamiltonians and degenerate idiosyncratic noise. For a comprehensive overview of this research direction, we refer readers to the works \cite{Ahu, AhuRenYan, CD1, GanMes, GanMesMouZha, BanMesMou, BanMes:master, MesMou, JacTan}. Additionally, \cite{GraMes:23} provides a comparative analysis of various monotonicity conditions.

The key contributions of this work can be summarized as follows: 
\begin{itemize} 
\item Under displacement monotonicity conditions, we obtain the long-time stability in the {\it non-compact} state space $\R^d$. This requires in particular the derivation of suitable uniform in time second moment estimates.
\item We consider a broad class of {\it non-separable} Hamiltonians, i.e., the structural condition \eqref{def:sep} is not imposed at any point in the manuscript. 
\item We develop a robust approach capable of handling both deterministic models and models with non-degenerate idiosyncratic noise, allowing for the case $\b = 0$. 
\item We obtain exponential decay properties not only for $\rho^{T}, D_{x}u^{T}$ but also for the value function $u^{T}$.
\end{itemize}
 Non-separable Hamiltonians arise naturally in applications. For a simple example consider $H(x,p,\mu) := \frac{1}{2}|p-\alpha\overline{m}_\mu|^2 - \frac{\kappa}{2}|x|^2$, where $\alpha\in\R,\kappa\ge 0$ and $\overline{m}_\mu :=\int_{\R^d}y\dd\mu(y)$. This Hamiltonian  is D-monotone for $\kappa>\alpha^2/4$. A potential economical interpretation could be as follows: the underlying MFG is for a continuum of firms, where $x\in\R$ stands for the inventory level and $p\in\R$ stands for some investment rate. In this scenario the optimal investment decision (i.e. $D_pH(x,p,\mu) = p-\alpha\overline{m}_\mu$) depends on population's aggregate state. Further Hamiltonian examples satisfying our standing assumptions will be discussed later.

Our first set of main results can be summarized informally as follows. We refer to the precise statements in Theorem \ref{prop:pointwisedecay} and Theorem \ref{prop:Du_loc}.

\begin{theorem}\label{thm:intro1}
Let $H:\R^{d}\times\R^{d}\times\sP_{2}(\R^{d})\to\R$ be displacement $c_{0}$-monotone with $c_{0}>0$ and suppose that it satisfies our standing assumptions. Let $\left(u^{1,T}_{s},\rho^{1,T}_{s}\right)_{s\in[0,T]}$ and $\left(u^{2,T}_{s},\rho^{2,T}_{s}\right)_{s\in[0,T]}$ be two solutions to \eqref{eq:MFG_intro} with initial/final data $\left(\rho^{1}_{0},g^{1}\right)$ and $\left(\rho^{2}_{0},g^{2}\right)$, respectively. Suppose that $\rho^{1}_{0},\rho^{2}_{0}\in\sP_{2}(\R^{d})$ and $g^{1},g^{2}:\R^{d}\times\sP_{2}(\R^{d})\to\R$ are both D-monotone and satisfy our standing assumptions. Then, there exists $C>0$ depending on  $H,g^{1},g^{2}, \int_{\R^{d}}|x|^{2}\dd\rho^{1}_{0}(x), \int_{\R^{d}}|x|^{2}\dd\rho^{2}_{0}(x)$ and there exists $\d>0$ depending only on $c_{0}$ such that
\[
W_{2}\left(\rho^{1,T}_{s},\rho^{2,T}_{s}\right) \le C \left[e^{-\d s} + e^{-\d(T-s)}\right],\ \ s\in[0,T],
\]
{ and, if we assume in addition that $\rho^1_0 = \rho^2_0$,}
\[
\ds\sup_{x\in\R^d}\frac{|D_xu^{1,T}(s,x)-D_xu^{2,T}(s,x)|^2}{1+|x|^2} \le { C e^{-\d (T-s)}, \ \ s\in[0,T]} .
\]
\end{theorem}
As a consequence of this theorem, we can formulate our second set of main results. For the precise statement we refer to Theorem \ref{thm:valueconvergence} and Corollary \ref{prop:tildeulimit1}.

%
%

\begin{theorem}\label{thm:intro2} 
Let $(u^T,\rho^T)$ be the solution to \eqref{eq:MFG_intro} with initial datum $\rho_{0}$ and any final datum $g$ which satisfies our standing assumptions. Then there exist $\l\in\R$,  $u:[0,+\infty)\times\R^d\to\R$ which is of class $C_{\rm{loc}}^{1,1}$ in space and Lipschitz continuous in time and $\rho\in C([0,+\infty); (\sP_2(\R^d),W_2))$, such that 
\begin{align*}
\sup_{s\in [0,t]} & W^2_2(\rho^T_s,\rho_s) \le C e^{-\d(T-t)}, \quad \forall  t \in[0, T], \\
\sup_{t \in[0, T/8], x \in \R^d} & \frac{| u^T (t,x) - \lambda (T-t) - u(t,x) |}{1+|x|^2} + \frac{|D_xu^{T}(t,x)-D_xu(t,x)|}{1+|x|} \le Ce^{-\delta T}
\end{align*}
where $C>0$ depends only on $\rho_{0}, H$ and $\d>0$ depends on $c_{0}$.

 Moreover, the triple $(u, \lambda, \rho)$ is the unique solution to the infinite horizon system
\begin{align}\label{limitsys}
\left\{
\begin{array}{ll}
-\partial_t  u  -\b\Delta u + H(x,-D_{x}  u,\rho) + \l = 0, & {\rm{in}\ } (0,+\infty)\times \R^d,\\
\partial_t\rho -\b\Delta \rho + \nabla\cdot (\rho D_pH(x,-D_{x}u,\rho)) = 0, & {\rm{in}\ } (0,+\infty)\times \R^d,\\
\rho(0,\cdot) = \rho_0,\ \
\ds\sup_{t\in[0,+\infty)}\frac{| u(t,x)|}{1+|x|^2} < \infty,
\end{array}
\right.
\end{align}
where the first equation is satisfied in the viscosity sense, while the second equation is satisfied in the sense of distributions.

\end{theorem} 

\medskip
\medskip
We now outline the main approach that enabled us to establish our main results Theorems \ref{thm:intro1} and \ref{thm:intro2}, which are grounded in techniques based on D-monotonicity. Notably, our analysis does not rely at all on solutions to the stationary system of type \eqref{eq:MFG_st}. For instance, the analysis in \cite{CecConDurEic} depends on a stationary system, requiring the separate construction of a solution, which in turn necessitated additional (e.g., smallness-type) assumptions on the data.

In contrast, our approach is deeply rooted in the Pontryagin maximum principle and a variety of FBSDE systems that characterize both Nash equilibria and individual agent trajectories. These systems write as 
\begin{equation}\label{eq:FBSDE_intro}
\left\{
\begin{array}{l}
\ds X^{t,\xi}_s=\xi+\int_{t}^sD_pH(X^{t,\xi}_\tau,Y^{t,\xi}_\tau,\rho_\t)\dd \t+\sqrt{2\beta} B_s^t, \\
\ds Y^{t,\xi}_s=-D_xg(X^{t,\xi}_T,\rho_T)+\int_s^T D_xH(X^{t,\xi}_\t,Y^{t,\xi}_\t,\rho_\t)\dd \t-\sqrt{2\beta}\int_s^TZ^{t,\xi}_\t\dd B_\t^t.
\end{array}
\right.
\end{equation}
Here $\beta\ge 0$, $(B_\t)_{\t\in[0,T]}$ is a given Brownian motion on $\R^{d}$ and we have set $B^t_s:=B_s-B_t$, $s\in[t,T]$.  Whenever $t=0$ and $\cL(X_s)=\rho_s$ for all $s\in[0,T]$, the system gives a Nash equilibrium. Furthermore, once a Nash equilibrium $(\rho_s)_{s\in[t,T]}$ is fixed, the system describes an individual agent's trajectory distributed as $\xi$ at time $t$. The core of our analysis relies on three main ingredients: (i) uniform in time second moment estimates for the processes $\left(X^{t,\xi}_\tau
\right)_{\tau\in[t,T]}$ and $\left(Y^{t,\xi}_\tau\right)_{\tau\in[t,T]}$; (ii) quantified D-monotonicity propagation estimates and (iii) the analysis of the dissipation of the $W_{2}$-distance between distinct MFG Nash equilibria.

To achieve (i), we discover new {\it generalized confining properties} for non-separable Hamiltonians which will guarantee a sort of semi-convexity property for the curves $s\mapsto \E\left[\left|X^{t,\xi}_s\right|^{2}\right]$ and $s\mapsto \E\left[\left|Y^{t,\xi}_s\right|^{2}\right]$, which in turn will lead to the desired uniform second moment estimates. We demonstrate in several examples how these generalized confining properties go hand in hand with the D-monotonicity.

{We pause for a moment to comment on the confining properties imposed on the Hamiltonian and final costs to achieve the desired second moment estimates. While for {\it mechanical Hamiltonians} of form $H(x,p,\mu) = \frac12|p|^{2}-f(x,\mu)$, the sufficient confining properties are naturally fully ensured by the strong D-monotonicity of the coupling function $f$ (and of $g$; as we show this in Subsection \ref{subsec:mechanical}), the situation is more subtle in the case of general Hamiltonians. As the velocity field for the Kolmogorov--Fokker--Planck equation is given by $(s,x)\mapsto D_{p}H(x,-D_{x}u(s,x),\rho_{s})$, this itself would need to have strong confining properties in order to ensure uniform in time moment bounds for the measure flow $(\rho_{s})_{s\in[0,T]}$. However, in our setting $x\mapsto D_{x}u(\cdot,x)$ will grow naturally linearly at infinity, and despite being a monotone field because of D-monotonicity, its composition with $D_{p}H$ might not have confining properties for free. And so, moment estimates on $s\mapsto D_{x}u(s,\cdot)$ have to be established separately. This is why further sufficient assumptions on $H$, beyond strong D-monotonicity, would need to be imposed. These are the ones listed in \eqref{hyp:H_confining}-\eqref{hyp:H_confining3}. These extra assumptions on $H$ ensure the uniform bounds on $s\mapsto \int_{\R^{d}}|x|^{2}\dd\rho_{s}(x)$ and $s\mapsto \int_{\R^{d}} |D_{x}u(s,x)|^{2}\dd\rho_{s}(x)$, which are precisely obtained via the second moment estimates on the $(X_{s},Y_{s})_{s\in[0,T]}$ processes appearing in the FBSDE system. 

Our two sets of assumptions \eqref{hyp:H_confining}-\eqref{hyp:H_confining3} on precise confining assumptions on $H$ and $g$} are formulated in the most general form. Such general forms arise because we handle the potential presence of noise and non-separable feature at the same time. In Lemma \ref{lem:genH_example} we give a set of more transparent sufficient conditions, which are stronger than the original ones, but easier to verify in applications. These can be formulated as
\begin{align*}
- \E\left[ (D_{pp}^2H(X,Y,\rho)X)\cdot D_{x}H(X,Y,\rho)\right] \ge \frac{\tilde \d^{1}_H}{2}\E\left[|X|^2\right] + \tilde c^{1}_{H},
\end{align*}
and 
\begin{align*}
- \E\left[ (D_{xx}^2H(X,Y,\rho)Y)\cdot D_{p}H(X,Y,\rho)\right] &- \E\tilde\E\left[ (D_{x\mu}^2H(X,Y,\rho,\tilde X)Y)\cdot D_{p}H(\tilde X,\tilde Y,\rho)\right]\\ 
&\ge \frac{\tilde \d^{2}_H}{2}\E\left[|Y|^2\right] + \tilde c^{2}_{H},
\end{align*}
for all $X,Y\in L^{2}(\Omega,\cF,\P;\R^{d})$ and $\rho=\cL(X)$. Here the constants $\tilde \d^{1}_H, \tilde \d^{2}_H>0$ and $ \tilde c^{1}_{H}, \tilde c^{2}_{H}\in\R$ would need to depend on second and third order derivative bounds of $H$ in a suitable way.

Let $(u^{i,T},\rho^{i,T})$, $i=1,2$ be two solutions to \eqref{eq:MFG_intro}, and we consider
$(X^{i,t,\xi^i}_s,Y^{i,t,\xi^i}_s,Z^{i,t,\xi^i}_s)_{s\in[t,T]}$, $i=1,2$, the solutions to the associated FBSDE systems \eqref{eq:FBSDE_intro}. 

For steps (ii) and (iii) in our program,  we introduce the functions $\varphi :(t,T)\to \R$, $\Phi:(t,T)\to[0,+\infty)$, given by 
\begin{equation*}
\varphi(s):=\mathbb E\left[(X_s^{1,t,\xi^1}-X_s^{2,t,\xi^2}) \cdot \left(Y^{1,t,\xi^1}_s-Y^{2,t,\xi^2}_s \right) \right]
\end{equation*}
and
\begin{equation*}
\Phi(s) 
:= \mathbb E\left[\left|X_s^{1,t,\xi^1}-X_s^{2,t,\xi^2}\right|^2\right] +\mathbb E \left[\left|Y^{1,t,\xi^1}_s-Y^{2,t,\xi^2}_s\right|^2\right].
\end{equation*}  
A crucial observation is that $\varphi$ and $\Phi$ satisfy a joint differential inequality, namely
\begin{align*}
c_0 |\varphi(s)| \le \frac {c_0} 2 \Phi(s) \le \varphi'(s) +\frac{C}{2c_0} W_1^2(\rho_s^1,\rho^2_s),
\end{align*}
where $c_{0}>0$ is the strong D-monotonicity constant of $H$ and $C>0$ depends only on the data. This inequality will then unfold a series of important consequences, which eventually lead to Theorem \ref{thm:intro1}. Along the way, we are also using the fact that as an implication of the growth conditions imposed for $DH$ and $Dg$, we have that $D_{x}u$ grow at most linearly at infinity, and so the decay in Theorem \ref{thm:intro1}(2) arises naturally.

\medskip

While $D_{x}u$ is very naturally connected to the $Y^{t,\xi}_s$ variable in \eqref{eq:FBSDE_intro}, the value function $u$ itself cannot be directly recovered from \eqref{eq:FBSDE_intro}. Therefore, its long time behavior cannot be deduced from Theorem \ref{thm:intro1} in a straightforward manner. So, in order to establish Theorem \ref{thm:intro2}, we rely on the Lagrangian representation formulas along Nash equilibria and optimal stochastic paths. For this, we would need to impose two additional hypotheses on the growth of $D_\mu H$ and $D_\mu g$, which are collected in \eqref{hyp:add_H} and \eqref{hyp:add_g}.

A very important role in this analysis is played by a suitably chosen average of a (partial) Lagrangian action, defined as 
\begin{equation*}
\lambda^T :=\mathbb{E}\left\{ \int_{T/2}^{T/2+1} L(X^{0,\xi,T}_s,D_pH(X^{0,\xi,T}_s,Y^{0,\xi,T}_s,\rho^T_s),\rho^T_s)\dd s\right\}.
\end{equation*}
It turns out that $\left(\lambda^T\right)_{T>0}$ becomes a Cauchy family of real numbers, whose limit will appear as $\l$ in the infinite horizon system \eqref{limitsys} in Theorem \ref{thm:intro2}. Our construction provides in fact the existence of solutions to such a system. We also show in Theorem \ref{uniquelim} that solutions to \eqref{limitsys} are (almost) unique: $u$ is defined only up to translations. The reader may observe that no final condition $g$ appears in  \eqref{limitsys}, hence the long-time behavior of $u^T(0,\cdot)$ depends on $\rho_0$ only, as one expects. We stress again that no stationary behavior in terms of $\bar u$ nor invariant measures $\bar \rho$ are exploited here, though we believe that these objects could be a posteriori reconstructed.

We also observe that the method proposed here, i.e. not relying directly on the stationary structure of the system (that is, the presence of a sole stationary state), could be used to study problems involving time dependent Hamiltonians, provided that the D-monotonicity is satisfied for all $t$. To the best of our knowledge, the investigation of the stability of nonautonomous systems (such as those with time-periodic data) is widely open within the context of MFG.

Note finally that the estimate on the rate of convergence of $u^T$ to $u$ is restricted to the time horizon $[0,T/8]$. In fact, this estimate could be extended to any interval of the form $[0,\zeta T]$, with $\zeta <1$, but not to the whole $[0,T]$ as one expects $u^T$ to deviate from $u$ in order to achieve the final condition $g$. This fact is also evident in the form of the estimate for the distance between $\rho^T$ and $\rho$. To describe precisely what happens at times close to $T$, we believe that one should study the convergence of $u^T$ to (time shifts of) solutions of the infinite horizon system
\begin{align*}
\left\{
\begin{array}{ll}
-\partial_t  \hat u  -\b\Delta \hat u + H(x,-D_{x}  \hat u,\hat \rho) + \l = 0, & {\rm{in}\ } (-\infty,0)\times \R^d,\\
\partial_t \hat \rho -\b\Delta \hat \rho + \nabla\cdot (\hat \rho D_pH(x,-D_{x}\hat u,\hat \rho)) = 0, & {\rm{in}\ } (-\infty,0)\times \R^d,\\
\hat u(0,\cdot) = g(\cdot,\hat \rho_0),\ \
\ds\sup_{t\in(-\infty,0)} M_2(\hat \rho_t) < \infty.
\end{array}
\right.
\end{align*}

\medskip

%

The structure of the rest of the paper is as follows. In Section \ref{sec:tech} we present the setting, recall the definitions of displacement monotonicity, and we collect a first class of assumptions on the data $H$ and $g$. Section \ref{sec:2nd_mom_Nash} is dedicated to the uniform in $T$ second moment estimates on the processes $\left(X^{t,\xi}_\tau\right)_{\tau\in[t,T]}$ and $\left(Y^{t,\xi}_\tau\right)_{\tau\in[t,T]}$ from \eqref{eq:FBSDE_intro}. This is where we present our proposed generalized confining properties on $H$ and $g$ and we construct classes of examples of Hamiltonians satisfying these new assumptions. 
Section \ref{sec:main} can be seen as the main one of this paper. Based on the results from the previous sections, this is where we essentially give all the remaining ingredients for proving Theorem \ref{thm:intro1}. Section \ref{sec:u} concerns the analysis of the value function $u$, using the Lagrangian representation formula. This is where we also construct the solution to the infinite horizon system and prove Theorem \ref{thm:intro2}. Finally, in Appendix \ref{sec:app} we have proven a global in time semi-concavity estimate for value functions arising in stochastic control problems. Such a result is probably well-known for experts, but in lack of a precise reference (which would be applicable both for uniformly parabolic and degenerate, first order problems, at the same time), we decided to give the details on this for the convenience of the reader. Then, we end the paper with Appendix \ref{app:ito}, where we have collected some well known It\^o lemmas along measure flows, and their applications in our setting.

\medskip
\medskip

\noindent {\bf Acknowledgements.} Both authors would like to thank the support by the King Abdullah University of Science and Technology Research Funding (KRF) under award no. ORA-2021-CRG10-4674.2; ARM has also been supported by the EPSRC New Investigator Award ``Mean Field Games and Master equations'' under award no. EP/X020320/1. MC has been partially funded by the EuropeanUnion-NextGenerationEU under the National Recovery and Resilience Plan (NRRP), Mission 4 Component 2 Investment 1.1 - Call PRIN 2022 No. 104 of February 2, 2022 of Italian Ministry of University and Research; Project 2022W58BJ5 (subject area: PE - Physical Sciences and Engineering) ``PDEs and optimal control methods in mean field games, population dynamics and multi-agent models'', and he acknowledges the support by the Gruppo Nazionale per l'Analisi Matematica, la Probabilit\`a e le loro Applicazioni (GNAMPA) of the Istituto Nazionale di Alta Matematica (INdAM), Italy. We thank two anonymous referees for their constructive suggestions which helped us improving the presentation of our results.

\section{Technical setting and assumptions}\label{sec:tech}

We denote by $\sP(\R^d)$ the set of Borel probability measures supported on $\R^d$. For $p\ge 1$ we set $\sP_p(\R^d):=\left\{\mu\in\sP(\R^d):\int_{\R^d}|x|^p\dd\mu(x)<+\infty\right\}$, to denote the set of probability measures with finite $p$-moment. For $\mu\in\sP_p(\R^d)$, we use the notation $M_p(\mu):=\left(\int_{\R^d}|x|^p\dd\mu(x)\right)^{\frac1p}$ to denote the $p$-moment of $\mu$. We equip $\sP_p(\R^d)$ with the classical Wasserstein distance $W_p$, defined as 
$$W_p(\mu,\nu)=\inf\left\{\iint_{\R^d\times\R^d}|x-y|^p\dd\pi(x,y):\ \pi\in\Pi(\mu,\nu)\right\}^{\frac1p},$$ 
where $\Pi(\mu,\nu)$ stands for the transport plans between $\mu$ and $\nu$, i.e. $\Pi(\mu,\nu):=\{\pi\in\sP_2(\R^d\times\R^d): (p^x)_\sharp\pi=\mu,\ (p^y)_\sharp\pi=\nu\}$. Here $p^x,p^y$ denote the canonical projections from $\R^d\times\R^d$ onto $\R^d$ and $\sharp$ stands for the pushforward operator.

We consider a complete probability space $(\Omega,\cF,\P)$, equipped with a filtration $\mathbb{F} = \left(\cF_t\right)_{t\ge 0}$, satisfying the usual assumptions. We consider furthermore an $\R^d$-valued $\mathbb{F}$-Brownian motion $(B_t)_{t\ge 0}$. We suppose that $(\Om,\cF_0,\P)$ is atomless, and so, for any $\mu\in\sP_2(\R^d)$ there exists $X\in L^2(\Omega,\cF_0,\P;\R^d)$ such that $\cL(X)=\mu$. Here, $\cL:L^2(\Omega,\cF,\P;\R^d)\to\sP_2(\R^d)$ stands for the standard law operator. 

It is well-known that the $W_2$ distance can be also formulated as 
$$W_2(\mu,\nu)=\inf\left\{\E\left[|X-Y|^2\right]:\ \cL(X)=\mu,\cL(Y)=\nu\right\}^{\frac12}.$$
Based on \cite{AGS}, for a function $f:\sP(\R^{d})\to\R$ we say that it is {\it Wasserstein differentiable} (or $W_{2}$-differentiable, or simply differentiable) at $\mu\in\sP_{2}(\R^{d})$, if there exists $\xi\in L^{2}(\mu;\R^{d})$ such that
\begin{align}\label{tool:diff}
f(\nu) = f(\mu) + \iint_{\R^{d}\times\R^{d}}\xi(x)\cdot (y-x)\dd\pi(x,y) + o(W_{2}(\mu,\nu)),\ \forall \pi\in\Pi(\mu,\nu),
\end{align}
as $W_{2}(\mu,\nu)\to 0$. There is a unique solution $\bar\xi\in L^{2}(\mu;\R^{d})$ of the variational problem 
$$\inf\left\{\|\xi\|_{L^{2}(\mu;\R^{d})}:\ \xi\ {\rm{satisfies}}\ \eqref{tool:diff}\right\},$$ 
and this vector field is referred to as the $W_{2}$-derivative of $f$ at $\mu$, that we denote as $D_{\mu}f(\mu,\cdot):=\bar\xi(\cdot)$. It is worth noting that a priori $D_{\mu}f(\mu,\cdot)$ is defined on $\spt(\mu)$. If $f$ is differentiable at any $\mu\in\sP_{2}(\R^{d})$ and the vector field $D_{\mu}f$ has a unique jointly continuous extension, i.e. $D_{\mu}f:\sP_{2}(\R^{d})\times\R^{d}\to\R^{d}$, then (referring to \cite[Chapter 5]{CD1}) we call $f$ to be {\it fully $C^{1}$}. Similarly, for $k>1$, one can define the class of fully $C^{k}$ functions over the Wasserstein space $(\sP_{2}(\R^{d}),W_{2})$. A deep result from \cite{GanTud:19} implies that $f$ is differentiable at $\mu\in\sP_{2}(\R^{d})$ if and only if there exists $\xi\in L^{2}(\mu;\R^{d})$ such that
$$
f(\nu) = f(\mu) + \E\left[\xi(X) \cdot (Y-X)\right] + o\left(\|X-Y\|_{L^{2}}\right),
$$ 
for any $X,Y\in L^{2}(\Om,\cF_{0},\P;\R^{d})$ with $\cL(X)=\mu$, $\cL(Y)=\nu$ and $\|Y- X\|_{L^{2}}\to 0.$

\subsection{First set of assumptions}

We suppose that the non-separable Hamiltonian $H:\R^d\times\R^d\times\sP_2(\R^d)\to\R$ is strongly convex in the $p$-variable and it is sufficiently regular. More precisely, 
\begin{align}\label{hyp:H_reg}\tag{H\arabic{hyp}}
&\bullet H(\cdot,\cdot,\mu)\in C^3(\R^d\times\R^d),\ \ \forall \mu\in\sP_2(\R^d);\\
&\bullet H, D_pH, D_{x}H\ {\rm{are\ fully\ }} C^1;\\ 
&\bullet D^2_{p\mu}H(x,p,\mu,\cdot), D^2_{x\mu}H(x,p,\mu,\cdot)\in C^1(\R^d),\ \forall (x,p,\mu)\in\R^d\times\R^d\times\sP_2(\R^d);\\
&\bullet D^2_{pp}H, D^{2}_{xx}H, D^2_{xp}H, D^{2}_{x\mu}H, D^2_{p\mu}H \ {\rm{are\ uniformly\ bounded}};\\
&\bullet  D^3_{ppx}H, D^3_{pxx}H, D^3_{ppp}H, D^3_{xxx}H, D^3_{p\mu\tilde x}H, D^3_{x\mu\tilde x}H \ {\rm{are\ uniformly\ bounded}}.
\end{align}
\stepcounter{hyp}
We immediately notice that the bounds on the second derivatives in \eqref{hyp:H_reg} imply that 
\begin{align}\label{hyp:DH_Lip}
D_pH {\rm{\ and\ }} D_{x}H {\rm{\ are\ globally\ Lipschitz\ continuous}}, 
\end{align}
when the Lipschitz continuity is taken with respect to $W_1$ in the measure variable. This further implies that 
\begin{align}\label{growth:DH}
&|D_{p}H(x,p,\mu)|\le C\left(1+ M_{1}(\mu)+|x|+|p|\right), \ \forall (x,p,\mu)\in \R^{d}\times\R^{d}\times\sP_{2}(\R^{d}),\\
&|D_{x}H(x,p,\mu)|\le C\left(1+ M_{1}(\mu)+|x|+|p|\right), \ \forall (x,p,\mu)\in \R^{d}\times\R^{d}\times\sP_{2}(\R^{d}),
\end{align}
where $C>0$ is a constant depending only on the uniform bounds on $D^{2}_{xx}H,D^{2}_{pp}H,D^{2}_{px}H$, $D^{2}_{x\mu}H,D^{2}_{p\mu}H$ and on $|D_xH(0,0,\d_0)|$ and $|D_pH(0,0,\d_0)|$.

Assume furthermore that $H$ jointly {\it strongly} displacement monotone, i.e. there exists $c_0>0$ such that $\forall X^1,X^2,P^1,P^2\in L^2(\Omega,\cF,\P;\R^d)$ and $\mu_1,\mu_2\in\sP_2(\R^d)$ such that $\cL(X^1)=\mu_1$, $\cL(X^2)=\mu_2$, we have
\begin{align}\label{hyp:Hdis_strong}\tag{H\arabic{hyp}}
\mathbb{E}\Big[\big(-D_xH(X^1,P^1,\mu_1)&+D_xH(X^2,P^2,\mu_2)\big)\cdot (X^1-X^2)\Big]\\
&\nonumber+\mathbb{E}\Big[\left(D_pH(X^1,P^1,\mu^1)-D_pH(X^2,P^2,\mu^2)\right)\cdot\left(P^1-P^2\right)\Big]\\
&\nonumber\ge c_0 \E\left[|X^1-X^2|^2+|P^1-P^2|^2\right].
\end{align}
\stepcounter{hyp}

We notice that the strong monotonicity assumption \eqref{hyp:Hdis_strong} implies that 
\begin{align}\label{hyp:Hdis_strong_2}
\big[-D_xH(x^1,p^1,\mu)&+D_xH(x^2,p^2,\mu)\big]\cdot (x^1-x^2)\\
&\nonumber+\left[D_pH(x^1,p^1,\mu)-D_pH(x^2,p^2,\mu)\right]\cdot\left(p^1-p^2\right)\\
&\nonumber\ge c_0 \left[|x^1-x^2|^2+|p^1-p^2|^2\right],
\end{align}
for all $x^1,x^2, p^1, p^2\in\R^d$ and for all $\mu\in\sP_2(\R^d).$ Indeed, this implication could be deduced from a straightforward adaptation of \cite[Lemma 2.5]{MesMou}.

\medskip

For the final cost functions $g$ we assume
\begin{align}\label{hyp:g}\tag{H\arabic{hyp}}
&\bullet g\in C^1(\R^d\times\sP_2(\R^d))\\
&\bullet D^2_{xx}g, D^2_{x\mu}g {\rm{\ are\ uniformly\ bounded.}}
\end{align}
\stepcounter{hyp}

Finally, we assume the displacement monotonicity condition on $g$, i.e. for all $\forall X^1,X^2\in L^2(\Omega,\cF,\P;\R^d)$ and $\mu_1,\mu_2\in\sP_2(\R^d)$ such that $\cL(X^1)=\mu_1$, $\cL(X^2)=\mu_2$, we have
\begin{align}\label{hyp:g_D-mon}\tag{H\arabic{hyp}}
\E\left[\left(D_{x}g(X^{1},\mu^{1})-D_{x}g(X^{2},\mu^{2})\right)\cdot \left(X^{1}-X^{2}\right)\right]\ge 0.
\end{align}
\stepcounter{hyp}

This implies that $g(\cdot, \mu)$ is convex for every $\mu$ (cf. \cite[Lemma 2.3]{MesMou}). We also note that the bounds on the second derivatives imply that
\begin{align}
D_{x}g \text{ is globally Lipschitz continuous},
\end{align}
with respect to $W_{1}$ in the measure variable, and
\begin{align}\label{cons:Dxg_growth}
|D_{x}g(x,\mu)|\le C(1+ M_{1}(\mu)+|x|),\ \forall (x,\mu)\in\R^{d}\times\sP_{2}(\R^{d}),
\end{align}
where $C>0$ depends only on the uniform bounds on $D^{2}_{xx}g$, $D^{2}_{x\mu}g$ and on $|D_{x}g(0,\d_{0})|$.


\section{Global second moment estimates along MFG Nash equilibria}\label{sec:2nd_mom_Nash}

The solution to the MFG can be fully characterized by the solution of an FBSDE system (which plays the role of the classical Hamiltonian system in the deterministic case). Let $t\in[0,T]$ and let $\xi$ be a random variable. Let moreover $\b\ge 0$ and let $(B_\t)_{\t\in[0,T]}$ be a given Brownian motion on $\R^{d}$ and set $B^t_s:=B_s-B_t$, $s\in[t,T].$ For a given flow of probability measures $(\rho_s)_{s\in[t,T]},$ consider the FBSDE system set on the time interval $(t,T)$
\begin{equation}\label{eq:FBSDE}
\left\{
\begin{array}{l}
\ds X_s=\xi+\int_{t}^sD_pH(X_\tau,Y_\tau,\rho_\t)\dd \t+\sqrt{2\beta} B_s^t, \\
\ds Y_s=-D_xg(X_T,\rho_T)+\int_s^T D_xH(X_\t,Y_\t,\rho_\t)\dd \t-\sqrt{2\beta}\int_s^TZ_\t\dd B_\t^t.
\end{array}
\right.
\end{equation}
{By a (strong) solution to this FBSDE system we mean an $\mathbb{F}$-progressively measurable, $\R^{d}\times\R^{d}\times\R^{d\times d}$-valued process $(X_s,Y_s, Z_s)_{s\in (t,T)}$ such that 
$$
\E\left[ \sup_{s\in[t,T]}\left(|X_{s}|^{2}+|Y_{s}|^{2}\right) + \int_{t}^{T}|Z_{s}|^{2}\dd s \right]<+\infty.
$$
Such systems, associated to the Pontryagin maximum principle in stochastic control, received a huge attention in the past three decades or so (see for instance the pioneering works \cite{PenWu, ParTan}). The so-called {\it probabilistic description} of MFG Nash equilibria via such FBSDEs, to the best of our knowledge, was initiated in \cite{CarDel:13}. This reference had shown in particular the connection between the solutions of \eqref{eq:MFG_intro} and \eqref{eq:FBSDE}, that we also describe later. For a thorough description of such links and further properties we refer to \cite[Chapter 4]{CD1}. In particular, \cite[Section 3.4.3, Section 4.5]{CD1} detail how to solve MFGs using the stochastic Pontryagin principle and the FBSDE system of type \eqref{eq:FBSDE}. We also refer to \cite[Section 5.2]{ChaCriDel22}, where the authors motivate various FBSDE formulations in connection with MFG. In particular, the convexity and monotonicity properties of the Hamiltonian and final condition make the use of the system \eqref{eq:FBSDE} natural. Our assumptions on the data, and in particular the displacement monotonicity conditions will let us to rely on that particular framework. We refer also to the works \cite{Ahu, MesMou, BanMesMou, JacTan} for further developments in the displacement monotone setting, where such FBSDEs played a crucial role.
}

\begin{definition}
We say that the flow of probability measures $(\rho_s)_{s\in[0,T]}$ is a {\emph{mean field game Nash equilibrium}}, if when solving \eqref{eq:FBSDE} over $[0,T]$ with input $(\rho_s)_{s\in[0,T]}$ and $\xi$ such that $\cL(\xi) = \rho_0$, we have $\cL(X_s)=\rho_s$ for all $s\in[0,T].$ 
\end{definition}

{
We note that as a consequence of the dynamic programming principle, if $t\in(0,T)$ and\\ $(X_s,Y_s,Z_s)_{s\in[t,T]}$ is a solution to \eqref{eq:FBSDE} with $\cL(\xi) = \rho_t$, then $\cL(X_s)=\rho_s$ for any $s\in(t,T)$. 
}

It is well-known that MFG Nash equilibria can be fully characterized by {sufficiently regular} solutions {(a regularity class present in our setting that we describe below)} of the MFG system: a coupled PDE system of a backward in time Hamilton--Jacobi--Bellman and a forward in time Kolmogorov--Fokker--Planck equation. This in our setting would write as the system given in \eqref{eq:MFG_intro}.

Furthermore, sufficiently regular solution to the HJB equation in \eqref{eq:MFG_intro} can be used to construct decoupling fields for \eqref{eq:FBSDE}: in particular, if $u\in C^{0,2}((0,T)\times\R^d)$, then $(\rho_s)_{s\in[0,T]}$ is a MFG Nash equilibrium if and only if $(X_s,Y_s,Z_s)_{s\in[0,T]}$ has the representation 
\begin{align}\label{form:decoupling}
Y_s = - D_x u(s,X_s)\ \ {\rm{and}}\ \ Z_s = - D^2_{xx} u(s,X_s).
\end{align} 

\begin{remark}
Under our standing assumptions on $H$ and $g$ imposed in this paper, both systems \eqref{eq:MFG_intro} and \eqref{eq:FBSDE} have unique solutions $(u,\rho)$ and $(X,Y,Z)$, respectively, for any $T>0$ fixed. In particular, the solvability of \eqref{eq:FBSDE} in the displacement monotone setting, under our standing assumptions, but in much more generality (involving also control of the volatility) has been addressed in \cite[Lemma 4]{JacTan}.

Recalling the main results from \cite{MesMou} and Lemma \ref{lem:semi_concave}, we have that $u\in C^{1}((0,T)\times\R^{d})$, and $D^2_{xx}u$ is uniformly bounded (independently of the value of $\beta$ and $T$). Furthermore $t\mapsto\rho_{t}$ is a continuous flow of probability measures (w.r.t. the $W_{1}$ distance). In addition, if $\beta>0$, we have that under our standing assumptions $u\in C^{1+\alpha/2, 2+\alpha}_{\rm{loc}}$. Indeed, since $H$ and $g$ are of class $C^1$, a standard bootstrap procedure involving Schauder estimates applies (for instance, as in the proof of Theorem 1.4 in \cite{CarPor20}).

The equivalence between the two MFG formulations \eqref{eq:MFG_intro} and \eqref{eq:FBSDE} are justified in \cite[Theorem 4.53]{CD1}. Under the regularity properties on $u$ that we have at hand, the representation \eqref{form:decoupling} was justified in \cite[Theorem 4.1]{MesMou}. We also note that because of the representation, the process $(Z_{s})_{s\in (t,T)}$ is valued in $\R^{d\times d}_{\sym}$ (real symmetric matrices).
\end{remark}

\subsection{Second moment estimates} A crucial part of the analysis that follows is based on uniform in time second moment estimates for MFG Nash equilibria. Depending on the Hamiltonian $H$, it will be convenient to introduce the functions $Q^{1}_H,Q^{2}_H:\left(L^2(\Om,\cF,\P;\R^d)\right)^2\times L^2(\Om,\cF,\P;\R^{d\times d}_{\sym})\times\sP_2(\R^d)\to\R$, defined as
\begin{align}\label{def:Q_H}
Q^{1}_H(X,Y,Z,\rho)&:= \E\left[|D_{p}H(X,Y,\rho)|^2 +(D_{px}^2H(X,Y,\rho)X)\cdot D_{p}H(X,Y,\rho) \right]\\
& - \E\left[ (D_{pp}^2H(X,Y,\rho)X)\cdot D_{x}H(X,Y,\rho)\right]\\
& +\beta \E\left\{ D_p{\rm{trace}}[D^2_{xx}H(X,Y,\rho)]\cdot X + 2{\rm{trace}}[D^2_{px}H(X,Y,\rho)] \right\}\\
& +\beta \E\left\{ {\rm{trace}}\left[(D^3_{ppp}H(X,Y,\rho)X){Z Z^\top} \right]\right\}\\
& +{2}\beta\, \E\, {\rm{trace}}\left\{\left[D^{3}_{ppx}H(X,Y,\rho)X+D^2_{pp}H(X,Y,\rho)\right]Z\right\}\\
& + \E\tilde\E \left[(D^2_{p\mu}H(X,Y,\rho,\tilde X)X)\cdot D_pH(\tilde X,\tilde Y,\rho)\right] + \beta \E\tilde\E \left\{ {\rm{trace}}\left[D^{3}_{p\mu \tilde x}H(X,Y,\rho,\tilde X)X\right]\right\}.
\end{align}
and
\begin{align}\label{def:Q2_H}
\\
Q^{2}_H(X,Y,Z,\rho)&:= \E\left[-(D^{2}_{xx}H(X,Y,\rho)Y)\cdot D_{p}H(X,Y,\rho) +|D_{x}H(X,Y,\rho)|^{2} \right]\\
& + \E \left[ (D_{xp}^2H(X,Y,\rho)Y)\cdot D_{x}H(X,Y,\rho)\right]\\
& -\beta\E \left\{ D_x{\rm{trace}}\left[D^2_{xx}H(X,Y,\rho\right]\cdot Y\right\}\\
& -\beta\E \left\{ {\rm{trace}}[(D^2_{xpp}H(X,Y,\rho)Y + 2 D^2_{px}H(X,Y,\rho)){Z Z^{\top}}] \right\}\\
& - {2}\beta\, \E{\rm{trace}}\left\{\left[D^{3}_{xxp}H(X,Y,\rho)Y+D^2_{xx}H(X,Y,\rho)\right]Z\right\}\\
& - \E\tilde\E \left[D^{2}_{x\mu}H(X,Y,\rho,\tilde X)Y\cdot D_pH(\tilde X,\tilde Y,\rho)\right] - \beta \E\tilde\E \left\{ {\rm{trace}}\left[D^{3}_{x\mu \tilde x}H(X,Y,\rho,\tilde X)Y\right]\right\},
\end{align}
{where $\tilde X$ and $\tilde Y$ stand for independent copies of the random variables $X$ and $Y$, respectively, and $\tilde\E$ emphasizes the expectation taken with respect to the corresponding $\sigma$-algebra.}

\begin{remark}
In fact, both $Q^{1}_{H}$ and $Q^{1}_{H}$ will always be evaluated as $Q^{i}_{H}(X,Y,Z,\cL(X))$, therefore we in fact just emphasize the concrete dependence on the measure in the last variable for these operators. {Similarly, these functionals will depend only on symmetric matrix valued random variables in the $Z$-coordinate, and so the terms involving $Z Z^\top$ will be replaced by $Z^2$.}
\end{remark}

\medskip

\subsubsection{Second set of standing assumptions}
Since our problems are set on the whole space $\R^d$, in order to have uniform in time control on second moments of MFG Nash equilibria, it will be crucial to have suitable `confining' properties on the corresponding optimal velocity fields (i.e. the drift term appearing in the Kolmogorov--Fokker--Planck equation). We can achieve such desired property by imposing general structural assumptions on the Hamiltonian $H$ and final cost function $g$. These will eventually enter into our list of main assumptions. 

\begin{align}\label{hyp:H_confining}\tag{H\arabic{hyp}}
&\forall Z\in L^\infty(\Om,\cF,\P;\R^{d\times d})\ \exists c^1_H = c^1_H (\|Z\|_{L^{\infty}})\in\R {\rm{\ and\ }} \delta^1_H > 0 {\rm{\ such\ that}}\\
&Q^{1}_H(X,Y,Z,\mu)\ge c^1_H+\frac{\d^1_H}{2}\E \left[|X|^2\right],\ \ \forall X,Y\in L^2(\Om,\cF,\P;\R^{d}),\ \cL(X)=\mu.
\end{align}
\setcounter{hyp_prime}{\arabic{hyp}}
\stepcounter{hyp}

\begin{align}\label{hyp:H_confining2}\tag{H\arabic{hyp}}
&\\
&\forall Z\in L^\infty(\Om,\cF,\P;\R^{d\times d}), \ \forall X\in L^2(\Om,\cF,\P;\R^{d}),\ \cL(X)=\mu,\ \ \exists c^2_H=c^2_H(\|Z\|_{L^{\infty}},\|X\|_{L^{2}})\in\R\\
& {\rm{\ and\ }} \delta^2_H = \d^2_H(\|Z\|_{L^{\infty}},\|X\|_{L^{2}}) > 0 {\rm{\ such\ that}}\\
&Q^{2}_H(X,Y,Z,\mu)\ge c^2_H+\frac{\d^2_H}{2}\E\left[ |Y|^2\right],\ \ \forall Y\in L^2(\Om,\cF,\P;\R^{d}).
\end{align}
\stepcounter{hyp}

\begin{align}\label{hyp:H_confining3}\tag{H\arabic{hyp}}
&\forall X\in L^2(\Om,\cF,\P;\R^{d}),\ \cL(X)=\mu, \forall\ c>0\ \ \exists c^3_H = c^3_H (\|X\|_{L^{2}},c)\in\R\\ 
&{\rm{\ and\ }} \delta^3_H = \d^3_H (c) > 0 {\rm{\ such\ that}}\\
&\E\left[c|Y|^{2} - Y\cdot D_{x}H(X,Y,\mu)\right]\ge \frac{\delta^{3}_{H}}{2}\E\left[|Y|^2\right] +c^3_H,\ \ \forall Y\in L^2(\Om,\cF,\P;\R^{d}).
\end{align}
\stepcounter{hyp}

\begin{align}\label{hyp:g_confining}\tag{H\arabic{hyp}}
&\exists c_g\in\R {\rm{\ and\ }} \delta_g > 0 {\rm{\ such\ that}}\\
&\E[X\cdot D_{p}H(X,-D_xg(X,\mu),\mu)]\le c_g - \frac{\d_g}{2}\E \left[|X|^2\right], \ \ \forall X \in L^2(\Om,\cF,\P;\R^{d}),\ \cL(X)=\mu.
\end{align}
\stepcounter{hyp}

\begin{remark}
Beside their dependence on the data $H$, we have emphasized the potential dependence on  $\|Z\|_{L^{\infty}}$ or $\|X\|_{L^{2}}$ of the constants $c^{i}_{H},\d^{i}_{H},\ i=2,3,$ appearing in the previous assumptions.
\end{remark}

\begin{remark}
$\ $

\begin{itemize}
\item While the particular Assumptions \eqref{hyp:H_confining} and \eqref{hyp:H_confining2} might look slightly involved, they provide natural confining properties for non-separable Hamiltonians. Here we would like to to emphasize that these general assumptions are implied by slightly stronger sufficient conditions as we present this in the next lemma.
\item Assumption \eqref{hyp:H_confining} is responsible for the uniform in time second moment estimates for the $(X_{s})_{s}$ process in \eqref{eq:FBSDE} (i.e. for the measure flow $s\mapsto\rho_{s}$ in \eqref{eq:MFG_intro}), while Assumptions \eqref{hyp:H_confining2} and \eqref{hyp:H_confining3} are responsible for the uniform in time moment estimates for the $(Y_{s})_{s}$ process in \eqref{eq:FBSDE} (i.e. for $s\mapsto - D_{x}u(s,X_s)$).
\item We would like to also note that long time stability of MFG Nash equilibria is in general more sophisticated than the stabilization of measure flows associated to solutions of Kolmogorov--Fokker--Planck equations. A na\"ive way to see this is that the driving field in the KFP equation is $(s,x)\mapsto D_{p}H(x,-D_{x}u(s,x),\rho_{s})$, which is depending on $D_{x}u$, an unknown itself. Furthermore, in our setting we do not expect $D_{x}u$ to be uniformly bounded, instead it can growth linearly at infinity. Therefore, one must investigate the stabilization for the pair $s\mapsto (\rho_{s},D_{x}u(s,\cdot))$ or equivalently for $s\mapsto (X_{s},Y_{s})$. 
\end{itemize}
\end{remark}

In the following lemma we present a class of sufficient conditions which will imply our general confining assumptions. While these are stronger than our original assumptions, they are straightforward to check in applications.

\begin{lemma}\label{lem:genH_example}
Suppose that $Z\in L^2(\Om,\cF,\P;\R^{d\times d}_\sym)$ is uniformly bounded by a constant $C_Z>0$. Suppose that $H$ satisfies the standing regularity assumptions \eqref{hyp:H_reg}. { Let
\begin{align*}
&c_{px}:=\sup_{q,p,\rho}|D^{2}_{px}H(q,p,\rho)|, c_{p\mu}:=\sup_{q,p,\rho,\tilde x}|D^{2}_{p\mu}H(q,p,\rho,\tilde q)|, \ \ c_{p\mu\tilde x}:=\sup_{q,p,\rho,\tilde x}|D^{3}_{p\mu\tilde x}H(q,p,\rho,\tilde q)|,
\end{align*}
and define similarly $c_{pp}, c_{ppx},c_{ppp},$ etc.
}

\noindent (i) 
Suppose that there exists $\tilde \d^{1}_{H}>0$ and $\tilde c^{1}_{H}\in\R$ such that 
{
\begin{align}\label{suf:1-1}
\tilde \d^{1}_{H} >  5\beta + c^{2}_{px} + c^{2}_{p\mu},
\end{align}
and
}
\begin{align}\label{ineq:compensationX}
- \E\left[ (D_{pp}^2H(X,Y,\rho)X)\cdot D_{x}H(X,Y,\rho)\right] \ge \frac{\tilde \d^{1}_H}{2}\E\left[|X|^2\right] + \tilde c^{1}_{H},
\end{align}
for all $X,Y\in L^{2}(\Omega,\cF,\P;\R^{d})$ and $\rho=\cL(X)$. Then \eqref{hyp:H_confining} is satisfied {with $\d^{1}_{H}:=\tilde \d^{1}_{H} - (5\beta + c^{2}_{px} + c^{2}_{p\mu})$ and $c^{1}_{H}:= \tilde c^{1}_{H} - \beta\left( \frac{c_{ppx}^2}{2} + 2 c_{xp} + \frac{c_{ppp}^2 C_Z^4}{2} + c_{ppx}^2 C_{Z}^{2} + 2c_{pp}C_Z  +\frac12 c^{2}_{p\mu\tilde x} \right).$}

\medskip

\noindent (ii) Suppose that there exists $\tilde \d^{2}_{H}>0$ and $\tilde c^{2}_{H}\in \R$ such that
{
\begin{align}\label{suf:2-1}
\tilde \d^{2}_{H} >  5\beta + c^{2}_{px},
\end{align}
}
and
\begin{align}\label{ineq:compensationY}
&- \E\left[ (D_{xx}^2H(X,Y,\rho)Y)\cdot D_{p}H(X,Y,\rho)\right]- \E\tilde\E\left[ (D_{x\mu}^2H(X,Y,\rho,\tilde X)Y)\cdot D_{p}H(\tilde X,\tilde Y,\rho)\right] \ge \frac{\tilde \d^{2}_H}{2}\E\left[|Y|^2\right] + \tilde c^{2}_{H},
\end{align}
for all $X,Y\in L^{2}(\Omega,\cF,\P;\R^{d})$ and $\rho=\cL(X)$. Then \eqref{hyp:H_confining2} is satisfied {with $\d^{2}_{H}:= \tilde \d^{2}_{H} - (5\beta + c_{px}^{2})$ and $c^{2}_{H}:=\tilde c^{2}_{H} + \beta\left(- \frac{c_{xxx}^{2}}{2} - \frac{c_{xpp}^{2}C_{Z}^{4}}{2}  - 2 c_{xp}C_{Z}^{2} - c_{xxp}^{2}C_{Z}^{2} - 2c_{xx}C_{Z} -\frac{c_{x\mu\tilde x}^{2}}{2} \right)$}

\end{lemma}

\begin{proof}

\noindent (i) After recalling the definition of $Q^{1}_H$ from \eqref{def:Q_H}, we have
\begin{align*}
Q^{1}_H(X,Y,Z,\rho)&= \E\left[|D_{p}H(X,Y,\rho)|^2 +(D_{px}^2H(X,Y,\rho)X)\cdot D_{p}H(X,Y,\rho) \right]\\
& - \E\left[ (D_{pp}^2H(X,Y,\rho)X)\cdot D_{x}H(X,Y,\rho)\right]\\
& +\beta \E\left\{ D_p{\rm{trace}}[D^2_{xx}H(X,Y,\rho)]\cdot X + 2{\rm{trace}}[D^2_{px}H(X,Y,\rho)] \right\}\\
& +\beta \E\left\{{\rm{trace}}\left[(D^3_{ppp}H(X,Y,\rho)X){Z^2}\right]\right\}\\
& +{2}\beta\, \E\, {\rm{trace}}\left\{\left[D^{3}_{ppx}H(X,Y,\rho)X+D^2_{pp}H(X,Y,\rho)\right]Z\right\}\\
& + \E\tilde\E \left[(D^2_{p\mu}H(X,Y,\rho,\tilde X)X)\cdot D_pH(\tilde X,\tilde Y,\rho)\right] + \beta \E\tilde\E \left\{ {\rm{trace}}\left[D^{3}_{p\mu \tilde x}H(X,Y,\rho,\tilde X)X\right]\right\}.
\end{align*}
So, we have 
\begin{align*}
Q^{1}_H(X,Y,Z,\rho)&{\ge \E\left[  |D_{p}H(X,Y,\rho)|^2 -\frac12 |D_{p}H(X,Y,\rho)|^2 -\frac{c_{px}^{2}}{2}|X|^2 \right] }\\
&{- \E\left[ (D_{pp}^2H(X,Y,\rho)X)\cdot D_{x}H(X,Y,\rho)\right]}\\
&{+\beta \E\left\{ - \frac{c_{ppx}^2}{2} -\frac12 |X|^2 - 2 c_{xp} - \frac{c_{ppp}^2 C_Z^4}{2} - \frac{1}{2}|X|^2 -c_{ppx}^2 C_{Z}^{2} - |X|^2 - 2c_{pp}C_Z\right\}}\\
&{+ \E\left[ -\frac{c_{p\mu}^{2}}{2}|X|^{2} - \frac12 |D_{p}H(X,Y,\rho)|^{2} - \frac\beta2 |X|^{2} - \frac12\beta c^{2}_{p\mu\tilde x} \right] }\\
&{= - \frac{5\beta + c_{px}^{2} + c_{p\mu}^{2}}{2}\E\left[|X|^{2} \right] - \E\left[ (D_{pp}^2H(X,Y,\rho)X)\cdot D_{x}H(X,Y,\rho)\right]}\\
&{ -\beta\left( \frac{c_{ppx}^2}{2} + 2 c_{xp} + \frac{c_{ppp}^2 C_Z^4}{2} + c_{ppx}^2 C_{Z}^{2} + 2c_{pp}C_Z  + \frac12 c^{2}_{p\mu\tilde x}\right)},\\
\end{align*}
{where we have used a series of Young's inequalities.}
The statement in (i) follows.
\medskip

\noindent (ii) The proof of the second claim follows similar lines of thought. We recall the definition of $Q^{2}_H$ from \eqref{def:Q2_H}, that we estimate below. 
\begin{align*}
Q^{2}_H(X,Y,Z,\rho)&:= \E\left[-(D^{2}_{xx}H(X,Y,\rho)Y)\cdot D_{p}H(X,Y,\rho) +|D_{x}H(X,Y,\rho)|^{2} \right]\\
& + \E \left[ (D_{xp}^2H(X,Y,\rho)Y)\cdot D_{x}H(X,Y,\rho)\right]\\
& -\beta\E \left\{ D_x{\rm{trace}}\left[D^2_{xx}H(X,Y,\rho)\right]\cdot Y\right\}\\
& -\beta\E \left\{ {\rm{trace}}[(D^2_{xpp}H(X,Y,\rho)Y + 2 D^2_{px}H(X,Y,\rho)){Z^{2}}] \right\}\\
& - {2}\beta\, \E{\rm{trace}}\left\{\left[D^{3}_{xxp}H(X,Y,\rho)Y+D^2_{xx}H(X,Y,\rho)\right]Z\right\}\\
& - \E\tilde\E \left[D^{2}_{x\mu}H(X,Y,\rho,\tilde X)Y\cdot D_pH(\tilde X,\tilde Y,\rho)\right] - \beta \E\tilde\E \left\{ {\rm{trace}}\left[D^{3}_{x\mu \tilde x}H(X,Y,\rho,\tilde X)Y\right]\right\}\\
\end{align*}
So, we have
\begin{align*}
Q^{2}_H(X,Y,Z,\rho)&{\ge - \E\left[ (D_{xx}^2H(X,Y,\rho)Y)\cdot D_{p}H(X,Y,\rho)\right] - \E\tilde\E\left[ (D_{x\mu}^2H(X,Y,\rho,\tilde X)Y)\cdot D_{p}H(\tilde X,\tilde Y,\rho)\right]}\\
&{ + \E\left[|D_{x}H(X,Y,\rho)|^{2} -\frac{c_{px}^{2}}{2}|Y|^{2} - \frac12|D_{x}H(X,Y,\rho)|^{2}  \right]}\\
&{ + \beta\left(- \frac{c_{xxx}^{2}}{2} -\frac12\E\left[|Y|^{2}\right] - \frac{c_{xpp}^{2}C_{Z}^{2}}{2} -\frac12\E\left[|Y|^{2}\right] - 2 c_{xp}C_{Z}^{2} \right)}\\
&{+\beta\left( - c_{xxp}^{2}C_{Z}^{2} - \E\left[|Y|^{2}\right] - 2c_{xx}C_{Z} -\frac{c_{x\mu\tilde x}^{2}}{2} - \frac{\beta}{2}\E\left[|Y|^{2}\right]\right)}\\
&{\ge - \E\left[ (D_{xx}^2H(X,Y,\rho)Y)\cdot D_{p}H(X,Y,\rho)\right] - \E\tilde\E\left[ (D_{x\mu}^2H(X,Y,\rho,\tilde X)Y)\cdot D_{p}H(\tilde X,\tilde Y,\rho)\right]}\\ 
& { -\frac{5\beta+c_{px}^{2}}{2}\E\left[|Y|^{2}\right] + \beta\left(- \frac{c_{xxx}^{2}}{2} - \frac{c_{xpp}^{2}C_{Z}^{4}}{2}  - 2 c_{xp}C_{Z}^{2} - c_{xxp}^{2}C_{Z}^{2} - 2c_{xx}C_{Z} -\frac{c_{x\mu\tilde x}^{2}}{2} \right), }
\end{align*}
where we have concluded by \eqref{ineq:compensationY}.
\end{proof}

\begin{remark}
$\ $

\begin{itemize}
\item From the previous lemma we see that the assumptions \eqref{suf:1-1}-\eqref{ineq:compensationX} and \eqref{suf:2-1}-\eqref{ineq:compensationY} could serve alternative sufficient conditions to \eqref{hyp:H_confining} and \eqref{hyp:H_confining2}, respectively.
\item From this lemma we also see that the confining assumptions impose in fact that the driving vector fields  $-D_xH$ and $D_pH$ in the FBSDE system are suitably confining in the $x$ and $p$-variables, respectively, when carefully weighted in an anisotropic way with respect to $D^2_{pp}H$ and $-D^2_{xx}H$ and $-D^2_{x\mu}H$. Such properties, in general cannot be deduced from the strong $D$-monotonicity properties. This is why they need to be imposed separately.
\end{itemize}
\end{remark}

\begin{remark}
It is important to notice that most of the terms appearing in the definitions of $Q^1_H$ and $Q^2_H$ is due to the potential presence of the idiosyncratic noise. Indeed, if $\beta = 0$, then those expressions simplify greatly, i.e. 
\begin{align*}
Q^{1}_H(X,Y,Z,\rho)&= \E\left[|D_{p}H(X,Y,\rho)|^2 +(D_{px}^2H(X,Y,\rho)X)\cdot D_{p}H(X,Y,\rho) \right]\\
& - \E\left[ (D_{pp}^2H(X,Y,\rho)X)\cdot D_{x}H(X,Y,\rho)\right]+ \E\tilde\E \left[(D^2_{p\mu}H(X,Y,\rho,\tilde X)X)\cdot D_pH(\tilde X,\tilde Y,\rho)\right]
\end{align*}
and
\begin{align*}
Q^{2}_H(X,Y,Z,\rho)&:= \E\left[-(D^{2}_{xx}H(X,Y,\rho)Y)\cdot D_{p}H(X,Y,\rho) +|D_{x}H(X,Y,\rho)|^{2} \right]\\
& + \E \left[ (D_{xp}^2H(X,Y,\rho)Y)\cdot D_{x}H(X,Y,\rho)\right] - \E\tilde\E \left[D^{2}_{x\mu}H(X,Y,\rho,\tilde X)Y\cdot D_pH(\tilde X,\tilde Y,\rho)\right].
\end{align*}

\end{remark}


\subsubsection{{Back to the moment estimates}}

The first main result of this subsection can be formulated as follows.

\begin{proposition}\label{prop:2nd_moment_bound}
Suppose that $H:\R^d\times\R^d\times\sP_2(\R^d)\to\R$ and $g:\R^d\times\sP_2(\R^d)\to\R$ satisfy \eqref{hyp:H_reg}, \eqref{hyp:Hdis_strong} and \eqref{hyp:g}, \eqref{hyp:g_D-mon} respectively. Suppose furthermore that the assumptions \eqref{hyp:H_confining} and \eqref{hyp:g_confining} are fulfilled.
Let $(\rho_s)_{s\in[t,T]}$ be a given MFG Nash equilibrium, corresponding to the solution $(X_s,Y_s,Z_s)_{s\in[t,T]}$ of \eqref{eq:FBSDE}. Then we have 
that there exists $C>0$ depending on the data, but independent of $T$, such that 
\begin{align*}
\E \left[|X_s|^2\right]\le \E\left[|X_t|^2\right] + C,\ \ \forall s\in[t,T].
\end{align*}
\end{proposition}

\begin{proof}
Since $(\rho_s)_{s\in[t,T]}$ is an MFG Nash equilibrium, we have in particular $\cL(X_s)=\rho_s$ for all $s\in[t,T]$. Let $(u,\rho)$ be the solutions to the MFG system \eqref{eq:MFG_intro}. As a consequence of displacement monotonicity imposed on $g$ and $L$, we have by \cite[Lemma 2.3, Lemma 2.5]{MesMou} $x\mapsto g(x,\mu)$ and $(x,v)\mapsto L(x,v,\mu)$ are convex, for all $\mu\in\sP_2(\R^d)$. Therefore, $x\mapsto u(s,x)$ is convex, uniformly with respect to $s\in [0,T]$ {(this fact is documented in \cite[Lemma 3.4]{MesMou}; see also Lemma \ref{lem:appdx})}. Furthermore, by Lemma \ref{lem:semi_concave}, we have that $x\mapsto u(s,x)$ is semi-concave with a semi-concavity constant independent of $T$ and $\beta$ (and the initial agent distribution $\rho_0$). Note that Lemma \ref{lem:semi_concave} applies since the Lagrangian $L(x,v,\rho_t)$ is convex in $v$ and semi-concave in $(x,v)$ by the uniform convexity of $H$ in $p$ and the boundedness assumptions on $D^2_{pp} H, D^2_{xx} H$. Therefore, there exists $C_u>0$ (independent of $T,\beta$ and $\rho_0$) such that 
\begin{equation}\label{eq:Cu}
\sup_{(s,x)\in[0,T]\times\R^d}\left|D^2_{xx}u(t,x)\right|\le C_u.
\end{equation}
As a consequence of this, from \eqref{form:decoupling} we obtain $|Z_s|\le C_u$ almost surely, for all $s\in [t,T]$.

\medskip

Let $\a>0$ whose value will be set later. We define $h_\alpha:\R\to\R$ as
$$h_\a(s):=\E\left[\frac\a2 |X_s|^2+X_s\cdot D_pH(X_s,Y_s,\rho_s) \right]$$ 
and compute
\begin{align*}
\frac{\dd}{\dd s} h_\a(s) &= \E\left[\a X_s\cdot D_pH(X_s,Y_s,\rho_s) \right] +\alpha\beta d + Q^{1}_H(X_s,Y_s,Z_s,\rho_s)\\
&\ge \a\E\left[X_s\cdot D_pH(X_s,Y_s,\rho_s) \right] + \alpha\beta d + c_H + \frac{\d_H}{2}\E|X_s|^2\\
& = \a h_\a(s) + \alpha\beta d + c_H +\frac{(\d_H-\a^2)}{2}|X_s|^2,
\end{align*}
where in the inequality above we have used \eqref{hyp:H_confining}. 

We recall that by our assumptions we have that $\d_H>0$ and $\d_g>0$. Then we choose $\a\in (0, \sqrt{\d_H}],$ and so, we deduce that 
$$
\frac{\dd}{\dd s} h_\a(s)\ge \a h_\a(s) + \alpha\beta d + c_H = \a(h_\a(s) + c_H/\a + \beta d).
$$ 
Gr\"onwall's inequality yields
\begin{align*}
h_\a(s)&\le e^{-\a (T-s)}(h_\a(T)+c_H/\a + \beta d) - c_H/\a - \beta d\\
&=e^{-\a (T-s)}\left[\E\left[\frac\a2 |X_T|^2+X_T\cdot D_pH(X_T,-D_xg(X_T,\rho_T),\rho_T)\right] +c_H/\a+ \beta d\right ] - c_H/\a - \beta d\\
&\le e^{-\a (T-s)}\left[\frac\a2 \E |X_T|^2+c_g - \frac{\d_g}{2}\E|X_T|^2 +c_H/\a + \beta d\right] - c_H/\a -\beta d,
\end{align*}
where in the last inequality we have used \eqref{hyp:g_confining}. If necessary, we decrease $\a$ further such that $\a\in (0,\d_g)$ and so, we obtain
$$
h_\a(s) \le e^{-\a (T-s)}\left[c_g  +c_H/\a + \beta d\right] - c_H/\a -\beta d \le C,
$$
where $C$ is a constant independent of the time variable, depending only on the data (in particular, on the constants $c_g, c_H, \d_g,\d_H$).
Using the definition of $h_\a$, with this choice of $\a$, we have just obtained
$$  
\E\left[\frac\a2 |X_s|^2+X_s\cdot D_pH(X_s,Y_s,\rho_s) \right]\le C.
$$
This further implies  that 
\begin{equation}\label{M2differentialineq}
\frac{\dd }{\dd s}\frac12\E\left[ |X_s|^2\right]\le C + \beta d - \a\frac12\E |X_s|^2,
\end{equation}
from where we obtain 
$$
\E\left[ |X_s|^2\right]\le e^{-\a (s-t)}\E\left[ |X_t|^2\right] +C,
$$
and so the claim follows.
\end{proof}

\begin{remark} Note that the displacement monotonicity assumption \eqref{hyp:Hdis_strong} is used in Proposition \ref{prop:2nd_moment_bound} only to guarantee uniform bounds on $D^2_{xx} u$. One can prove the same second moment bounds by assuming directly the control on second derivatives of $u$, which may hold beyond the D-monotone setting (for instance as a consequence of uniform parabolic estimates).
\end{remark}

We can formulate a similar property for the dual precess $(Y_{s})_{s\in[t,T]}$ from \eqref{eq:FBSDE}.

\begin{proposition}\label{prop:L2_Y}
Suppose that the assumptions from the statement of Proposition \ref{prop:2nd_moment_bound} and the additional assumptions \eqref{hyp:H_confining2} and \eqref{hyp:H_confining3} take place. Let $(\rho_s)_{s\in[t,T]}$ be a given MFG Nash equilibrium, corresponding to the solution $(X_s,Y_s,Z_s)_{s\in[t,T]}$ of \eqref{eq:FBSDE}.  Then we have that there exists $C>0$ depending on the data and $\E\left[|X_t|^2\right]$, but independent of $T$, such that 
\begin{align*}
\E|Y_s|^2\le C\left(1 + \E\left[|X_{t}|^{2}\right] \right),\ \ \forall s\in[t,T].
\end{align*}
\end{proposition}

\begin{proof}
The proof follows similar ideas as the one of Proposition \ref{prop:2nd_moment_bound}. Let $\a>0$, to be chosen later and consider $h_{\a}:[t,T]\to\R$ defined as
$$
h_{\a}(s):=\E\left[\frac{\a}{2}|Y_{s}|^{2} - Y_{s}\cdot D_{x}H(X_{s},Y_{s},\rho_{s})\right].
$$
We compute
\begin{align*}
\frac{\dd}{\dd s} h_\a(s) &= - \E\left[\a Y_s\cdot D_xH(X_s,Y_s,\rho_s) \right] + \alpha\beta {\E[{\rm trace} (Z^2)]} + Q^{2}_H(X_s,Y_s,Z_s,\rho_s)\\
&\ge -\a\E\left[Y_s\cdot D_xH(X_s,Y_s,\rho_s) \right] + c^{2}_{H} +\frac{\delta^{2}_{H}}{2}\E\left[|Y_{s}|^{2}\right]\\
& \ge \a h_\a(s) + \tilde C +\frac{(\d^{2}_H-\a^2)}{2}|Y_s|^2,
\end{align*}
where the penultimate inequality we have used \eqref{hyp:H_confining2} {and the fact that since $Z\in\R^{d\times d}_\sym$ we have $\alpha\beta \E[{\rm trace} (Z^2)]\ge 0$}, and have set $\tilde C:=  c^{2}_{H}$. Now, one chooses $\a<\sqrt{\d^{2}_H}$, and so, one obtains
$$
\frac{\dd}{\dd s} \left(h_\a(s)+\tilde C/\alpha\right) \ge \a \left(h_\a(s) + \tilde C/\alpha\right).
$$
Thus, $e^{-\a(T-s)}\left(h_{\a}(T)+\tilde C/\alpha\right)-\tilde C/\alpha\ge h_{\alpha}(s)$ for all $s\in[t,T]$. In addition to this,
\begin{align*}
h_{\a}(T) & = \E\left[\frac{\a}{2}|Y_{T}|^{2} - Y_{T}\cdot D_{x}H(X_{T},Y_{T},\rho_{T})\right]\\
&= \E\left[\frac{\a}{2}|D_{x}g(X_{T},\rho_{T})|^{2} + D_{x}g(X_{T},\rho_{T})\cdot D_{x}H(X_{T},-D_{x}g(X_{T},\rho_{T}),\rho_{T})\right]\\
&\le c + c\E\left[|X_{T}|^{2}\right],
\end{align*}
where we have used the growth properties of $D_{x}H$ and $D_{x}g$ from \eqref{growth:DH} and \eqref{cons:Dxg_growth}, respectively, and Young's inequality, and the constant $c>0$ depends only on $\alpha$, and the constants appearing in the growth inequalities. By Proposition \ref{prop:2nd_moment_bound}, we will be able to say that there is a constant $C>0$, independent of $T$ such that
$$
h_{\a}(s)\le C + C\E\left[|X_{t}|^{2}\right], \ \ \forall s\in[t,T].
$$
Thus, by the definition of $h_{\a}$ and \eqref{hyp:H_confining3} we obtain from here that 
\begin{align*}
\frac{\delta^{3}_{H}}{2}\E\left[|Y_{s}|^2\right] +c^3_H \le\E\left[\frac{\a}{2}|Y_{s}|^{2} - Y_{s}\cdot D_{x}H(X_{s},Y_{s},\rho_{s})\right] \le C + C\E\left[|X_{t}|^{2}\right],
\end{align*}
where the constants $\delta^{3}_{H}$ and $c^3_H$ depend on $\alpha$ and in addition, $c^3_H$ might depend also on $\E\left[|X_{t}|^{2}\right]$. We conclude by Proposition \ref{prop:2nd_moment_bound}.
\end{proof}

\begin{remark}
\begin{enumerate}
\item Let us emphasize that in general, one must have strict inequalities, i.e. in all $\d^1_H>0,$ $\d^2_H>0,$ $\d^{3}_{H}>0$ and $\d_g>0$ appearing in \eqref{hyp:H_confining}-\eqref{hyp:H_confining3}  and \eqref{hyp:g_confining}, in order to be able to have the previously described uniform in time second moment estimates. This is partly due to the presence of the idiosyncratic noise. Indeed, one would be able to relax slightly these assumptions, and allow $\d^{i}_H = 0$ or  $\d_g = 0$, as long as we would be able to have $c^{1}_H, c^{2}_{H}\ge0$ and $c_g\le 0$. However, in the presence of the idiosyncratic noise one has a $Z$ contribution in the expressions of $Q^{1}_H, Q^{2}_{H}$ (corresponding to $D^2_{xx}u$), and so it would be in general impossible to have the inequality $c^{i}_H\ge 0$ satisfied. In the case of deterministic models, i.e. when $\beta = 0$, we would be able to use such a slightly relaxed version of the assumptions \eqref{hyp:H_confining}, \eqref{hyp:H_confining2} and \eqref{hyp:g_confining}. 
\item In the same time, it is important to remark that these assumptions in fact guarantee a sort of `confining property' for the driving vector fields in the corresponding MFG models, in the absence of which in general, one cannot hope for uniform in time second moment estimates (since we work in the non-compact setting of $\R^d$). Therefore, it would in general not be possible to allow $\d^{i}_H<0$ or $\d_g<0$, even for deterministic models.
\end{enumerate}
\end{remark}

\begin{corollary}\label{cor:D_xu_lin}
Suppose that the assumptions of Proposition \ref{prop:L2_Y} are fulfilled. Then, there exists a constant $C>0$ depending only on the data and $M_2(\rho_0)$ (but independent of $T$) such that
\begin{align}\label{gradbound}
|D_xu(s,x)| \le C(1+|x|),\ \forall (s,x)\in [t,T]\times\R^d.
\end{align}
\end{corollary}

\begin{proof}
This is a direct consequence of Propositions \ref{prop:2nd_moment_bound} and \ref{prop:L2_Y}, and the semi-concavity estimates of Proposition \ref{lem:semi_concave}, which guarantee the uniform bound $\|D^2_{xx} u\|_\infty \le C_u $ (see \eqref{eq:Cu}). Indeed, for any $(s,x)\in [t,T]\times\R^d$, we have
\[
|D_x u(s,x)| \le |D_x u(s,X_s)| + C_u |x-X_s| = |Y_s| + C_u |x-X_s|,
\]
hence
\[
|D_x u(s,x)| \le \E[|Y_s|] + C_u |x| + C_u \E[|X_s|]
\]
and by the quoted propositions we obtain the result.
 \end{proof}

\subsection{Examples of $H,g$ satisfying the  `generalized confining' assumptions}\label{subsec:examples}
Now let us pause to provide natural examples for $H$ and $g$, which will satisfy the assumptions imposed in \eqref{hyp:H_confining}, \eqref{hyp:H_confining2}, \eqref{hyp:H_confining3} and \eqref{hyp:g_confining}.

\subsubsection{{ Mechanical Hamiltonians}}\label{subsec:mechanical}
The simplest possible example we can think of is when
\begin{align}\label{H:mechanical}
H(x,p,\mu)=\frac12|p|^2-f(x,\mu).
\end{align}
Though we are now focusing on \eqref{hyp:H_confining}--\eqref{hyp:g_confining}, the other assumptions \eqref{hyp:H_reg} and \eqref{hyp:Hdis_strong} can be also checked easily in this situation. We can formulate the following result.

\begin{lemma}\label{lem:mechanical}
Suppose that $f:\R^d\times\sP_2(\R^d)\to\R$ is such that $f(\cdot, \mu) \in C^3$, $D_xf$ is $C^1$ and $D^2_{x\mu}f(x, \mu, \cdot)$ is bounded and with bounded first-order derivatives. Assume that $f$ is strongly displacement monotone, that is, for some $c_0 > 0$ it holds
\begin{equation}\label{fdisplmon}
\E\left[D_{x}f(X^{1},\cL(X^{1}))-D_{x}f(X^{2},\cL(X^{2}))\cdot(X^{1}-X^{2})\right]\ge c_{0}\E\left[|X^{1}-X^{2}|^{2}\right],\ \forall X^{1},X^{2}\in L^{2}(\Omega,\cF,\P;\R^{d}).
\end{equation}
Then, \eqref{hyp:H_reg}, \eqref{hyp:Hdis_strong} and \eqref{hyp:H_confining}--\eqref{hyp:H_confining3} are verified.
%
%
%
%
\end{lemma}

Note that \eqref{fdisplmon} is the same as condition \eqref{hyp:g_D-mon} for $g$. The latter implies \eqref{hyp:g_confining} in this simple case, arguing exactly as in the proof below.

\begin{proof}[Proof of Lemma \ref{lem:mechanical}] First, the regularity assumption \eqref{hyp:H_reg} can be checked directly, while condition \eqref{hyp:Hdis_strong} is a straightforward consequence of \eqref{fdisplmon} when $H = |p|^2/2 - f$.

To check \eqref{hyp:H_confining}, note first that the separability of $H$ implies that $D^2_{p\mu}H$ vanishes everywhere. In this case we obtain
\begin{align*}
Q^{1}_H(X,Y,Z,\rho)&= \E\left[|Y|^2 \right] + \E\left[ X\cdot D_{x}f(X,\rho)\right] +\beta \E\left[{\rm{trace}}Z\right],
\end{align*}
and so, the assumption \eqref{hyp:H_confining} boils down (since $Z$ is supposed to be bounded) to the inequality
\begin{align}\label{hyp:Df_confining}
\E\left[ X\cdot D_{x}f(X,\rho)\right]\ge c_f + \frac{c'_f}{2}\E|X|^2, \ \ \forall X\in L^{2}(\Omega,\cF,\P;\R^{d}), \ \ \rho=\cL(X),
\end{align}
that is required to hold for for some $c_f\in\R$ and $c'_f>0$. This is in turn a consequence of \eqref{fdisplmon}, once we take for example $X_1 = X$ and  $\cL(X_2) = \delta_0$, and we apply Young's inequality.

Let us now proceed with \eqref{hyp:H_confining2}. Direct computation yields
\begin{align*}
Q^{2}_H(X,Y,Z,\rho)&= \E\left[(D^{2}_{xx}f(X,\rho)Y)\cdot Y +|D_{x}f(X,\rho)|^{2} \right]\\
& +\beta\E \left\{ D_x{\rm{trace}}\left[D^2_{xx}f(X,\rho)\right]\cdot Y\right\}\\
& +\beta\, \E{\rm{trace}}\left\{\left[D^2_{xx}f(X,\rho)\right]Z\right\}\\
& + \E\tilde\E \left[\left(D^{2}_{x\mu}f(X,\rho,\tilde X)Y\right)\cdot \tilde Y\right] + \beta \E\tilde\E \left\{ {\rm{trace}}\left[D^{3}_{x\mu \tilde x}f(X,\rho,\tilde X)Y\right]\right\}.
\end{align*}
Now, we notice that, using second order Taylor expansion, \eqref{fdisplmon} is in fact equivalent to
$$
\E\left[(D^{2}_{xx}f(X,\cL(X))\xi)\cdot \xi \right] +  \E\tilde\E \left[\left(D^{2}_{x\mu}f(X,\cL(X),\tilde X)\xi\right)\cdot \tilde \xi\right]\ge c_{0}\E\left[|\xi|^{2}\right],
$$
for any $X,\xi\in L^{2}(\Omega,\cF,\P;\R^{d}).$ We combine this strong monotonicity inequality with the growth property \eqref{growth:DH} for $D_{x}f$, with the uniform boundedness of second and third order derivatives of $f$, to conclude about \eqref{hyp:H_confining2} after multiple use of Young's inequality.

To check \eqref{hyp:H_confining3}, we compute
\begin{align*}
\E\left[c|Y|^{2} - Y\cdot D_{x}H(X,Y,\rho)\right] = \E\left[c|Y|^{2} +Y\cdot D_{x}f(X,\rho)\right].
\end{align*}
Again, using the growth property \eqref{growth:DH} for $D_{x}f$ together with Young's inequality, we obtain the desired property.
\end{proof}

\begin{remark}
{ The reader may have noticed that in this separable, quadratic example, the strong displacement monotonicity assumption \eqref{hyp:Hdis_strong} gives the confining properties \eqref{hyp:H_confining}--\eqref{hyp:H_confining3} for free. However, for general Hamiltonians there does not seem to be a direct relationship between the two set of assumptions, and one needs to check them separately. 

The reason is that the confining assumptions are used to produce useful differential inequalities such as the one on $s\mapsto\frac12\E|X_s|^2$ in \eqref{M2differentialineq}, along the solution $(X_s,Y_s,Z_s)_{s\in[t,T]}$ to \eqref{eq:FBSDE}. On the other hand, displacement monotonicity is useful to prove that $s \mapsto E \left[\left(Y^{1}_s  -Y^{2}_s\right)\cdot\left(X_s^{1}-X_s^{2}\right) \right]$ is a sort of Lyapunov function (see Lemma \ref{lem:Phi_building} below). Since the two arguments look at the evolution of different quantities, they are likely to need different assumptions. }
\end{remark}

\begin{remark} { One might be tempted to add a drift term $x\mapsto b(x)$ to this model problem, in which case the Hamiltonian would read $H(x,p,\mu)=\frac12|p|^2-b(x) \cdot p-f(x,\mu)$. This can be certainly done, but would fall into the ``perturbative'' situation described in the following subsection. In other words, even if one imposes some nice properties on $b$, such as dissipativity and/or confinement, these do not play a clear role in the verification of our assumptions (unless the very special case of Ornstein--Uhlenbeck drift is considered).
}
\end{remark}

\subsubsection{General classes of Hamiltonians}\label{sec:general_H}


{
\begin{remark}
Lemma \ref{lem:genH_example} gives sufficient conditions on the Hamiltonian, which are straightforward to verify. Prototypical examples of non-separable displacement monotone Hamiltonians, as discussed in \cite[Remark 2.8]{MesMou}, will satisfy the standing assumptions of this paper, as we demonstrate in the next corollary.
\end{remark}
}

\begin{corollary}\label{ex:non-sepH}
The Hamiltonians for the form
$$
H(x,p,\mu) = H_0(x,p,\mu) +\frac{C_0}{2}\left(|p|^2 - |x|^2\right),
$$ 
where $H_0:\R^d\times\R^d\times\sP_2(\R^d)\to\R$ is uniformly bounded, with uniformly bounded derivatives up to order 3, and $C_0>0$ is a sufficiently large constant (depending on the uniform derivative bounds of $H_0$) will satisfy all the generalized confining assumptions \eqref{hyp:H_confining}, \eqref{hyp:H_confining2} and \eqref{hyp:H_confining3}. 
\end{corollary}

{
\begin{proof}
We verify the condition \eqref{hyp:H_confining}, using Lemma \ref{lem:genH_example}. For this we directly compute
\begin{align*}
- \E\big[ & (D_{pp}^2H(X,Y,\rho)X)\cdot D_{x}H(X,Y,\rho)\big]\\
& = - \E\left[ (D_{pp}^2H_{0}(X,Y,\rho)X+C_{0}X)\cdot \left(D_{x}H_{0}(X,Y,\rho)-C_{0}X\right)\right]\\
& = - \E\left[ (D_{pp}^2H_{0}(X,Y,\rho)X)\cdot D_{x}H_{0}(X,Y,\rho)\right] + C_{0} \E\left[X^{\top}D_{pp}^2H_{0}(X,Y,\rho)X - X\cdot D_{x}H_{0}(X,Y,\rho)\right]\\
& + C_{0}^{2}\E\left[|X|^{2}\right]\\
&\ge -\frac12 (c_{0,pp}c_{0,x})^{2} -\frac12\E\left[|X|^{2}\right] - C_{0}\left(c_{0,pp}\E\left[|X|^{2}\right]+\frac12\E\left[|X|^{2}\right] +\frac12 c_{0,x}^{2}\right)+ C_{0}^{2}\E\left[|X|^{2}\right]\\
&= \E\left[|X|^{2}\right]\left(C_{0}^{2} - C_{0}c_{0,pp}-C_{0}/2-1/2\right) - (c_{0,pp}c_{0,x})^{2}/2 - C_{0}c_{0,x}^{2}/2,
\end{align*}
where $c_{0,pp}$ is the uniform upper bound on $D^{2}_{pp}H_{0}$, while $c_{0,x}$ is the uniform upper bound on $D_{x}H_{0}$. Now, it is clear that we can choose $C_{0}$ large enough such that both
$$
C_{0}^{2} - C_{0}c_{0,pp}-C_{0}/2-1/2>0
$$
and
$$
C_{0}>  5\beta + c^{2}_{px} + c^{2}_{p\mu} = 5\beta + c^{2}_{0,px} + c^{2}_{0,p\mu},
$$
where $c_{0,px}$ is the uniform upper bound of $D^{2}_{px}H_{0}$ and $c_{0,p\mu}$ is the uniform upper bound of $D^{2}_{p\mu}H_{0}$, and so the assumptions of Lemma \ref{lem:genH_example} (i) are fulfilled. 

\medskip

The verification of \eqref{hyp:H_confining2} is analogous, and so, we leave this to the reader.

\medskip

Finally, let us verify \eqref{hyp:H_confining3}. For $c>0$, we compute
\begin{align*}
\E\left[c|Y|^{2} - Y\cdot D_{x}H(X,Y,\mu)\right] &= \E\left[c|Y|^{2} - Y\cdot \left(D_{x}H_{0}(X,Y,\mu)-C_{0}X\right)\right]\\
&\ge c \E\left[|Y|^{2}\right] - \E\left[|Y| |D_{x}H_{0}(X,Y,\mu)| - C_{0}|Y| |X|\right]\\
&\ge c \E\left[|Y|^{2}\right] - \frac{\varepsilon}{2} \E\left[|Y|^{2}\right]  -\frac{c_{0,x}^{2}}{2\varepsilon} -\frac{\varepsilon}{2} \E\left[|Y|^{2}\right] - \frac{C_{0}^{2}}{2\varepsilon} \E\left[|X|^{2}\right],
\end{align*}
for any $\varepsilon>0$. Now, if we choose $\varepsilon:=\frac c2$, we see \eqref{hyp:H_confining3} satisfied with $\d^{3}_{H}:=c$ and 
$$c^{3}_{H}:= -\frac{c_{0,x}^{2}}{c} - \frac{C_{0}^{2}}{c} \E\left[|X|^{2}\right].$$
\end{proof}

\begin{remark}
Typical running costs (Lagrangians) could be easily obtained from Hamiltonians via the Legendre transform, i.e. 
$$
L(x,v,\mu):=\sup_{p\in\R^{d}}\left\{v\cdot p - H(x,p,\mu)\right\}.
$$
Similarly to the example for the Hamiltonian in Corollary \ref{ex:non-sepH}, a typical non-separable Lagrangian in our context would be
$$
L(x,v,\mu):=L_{0}(x,v,\mu) +\frac{C_{0}}{2}\left(|v|^{2}+|x|^{2}\right),
$$
where $L_{0}$ has uniformly bounded first, second and third order derivatives and $C_{0}>0$ is a large enough constant. 
\end{remark}
}

\section{Long-time behavior of $\rho$ and $D_x u$}\label{sec:main}


Let $(\rho_s)_{s\in[0,T]}$ be a given MFG Nash equilibrium. To characterize single agent trajectories starting at position $x\in\R^d$ at time $t\in[0,T]$, it is convenient to introduce the system
\begin{equation}\label{eq:FBSDEx}
\left\{
\begin{array}{l}
\ds X^{t,x}_s=x+\int_{t}^sD_pH(X^{t,x}_\tau,Y^{t,x}_\tau,\rho_\t)\dd \t+\sqrt{2\beta} B_s^t, \\
\ds Y^{t,x}_s=-D_xg(X^{t,x}_T,\rho_T)+\int_s^T D_xH(X^{t,x}_\t,Y^{t,x}_\t,\rho_\t)\dd \t-\sqrt{2\beta}\int_s^TZ^{t,x}_\t\dd B_\t^t.
\end{array}
\right.
\end{equation}
This system is a particular version of the more generic one
\begin{equation}\label{eq:FBSDExi}
\left\{
\begin{array}{l}
\ds X^{t,\xi}_s=\xi+\int_{t}^sD_pH(X^{t,\xi}_\tau,Y^{t,\xi}_\tau,\rho_\t)\dd \t+\sqrt{2\beta} B_s^t, \\
\ds Y^{t,\xi}_s=-D_xg(X^{t,\xi}_T,\rho_T)+\int_s^T D_xH(X^{t,\xi}_\t,Y^{t,\xi}_\t,\rho_\t)\dd \t-\sqrt{2\beta}\int_s^TZ^{t,\xi}_\t\dd B_\t^t,
\end{array}
\right.
\end{equation}
where $\xi\in L^2(\Om,\cF_{t},\P;\R^{d})$. Indeed, if we consider $\xi$ such that $\cL(\xi)=\d_x$, we obtain essentially \eqref{eq:FBSDEx}. It is important to note that the system \eqref{eq:FBSDExi} (and so \eqref{eq:FBSDEx} as well) in general does not describe MFG Nash equilibria (as $(\rho_s)_{s\in[0,T]}$ is considered to be an input), unless $\cL(\xi)=\rho_t$ and $\cL(X^{t,\xi}_s)=\rho_s$ for all $s\in[0,T]$. Moreover, just as in the case of \eqref{form:decoupling}, we also have

\begin{align}\label{form:decoupling_loc}
Y^{t,\xi}_s = - D_x u(s,X^{t,\xi}_s)\ \ {\rm{and}}\ \ Z^{t,\xi}_s = - D^2_{xx} u(s,X^{t,\xi}_s),
\end{align} 
{and in particular
\begin{align*}
Y^{t,x}_t = - D_x u(t,x)\ \ {\rm{and}}\ \ Z^{t,x}_t = - D^2_{xx} u(t,x).
\end{align*}

\subsection{Some preparatory results}

{We underline that whenever two solutions to the FBSDE systems are used, these are always driven by the same fixed Brownian motion. This is possible because of the strong solvability of the FBSDE systems, results which we have recalled from the literature earlier.}

\begin{definition}\label{def:phi_fcts}
Let $(u^1,\rho^1)$ and $(u^2,\rho^2)$ be two solutions to the mean field game system over $[0,T]$ with initial and final data $(\rho^1_0,g^1)$ and $(\rho^2_0,g^2)$, respectively. Let $t\in[0,T]$ be given and let\\ $(X^{i,t,\xi^i}_s,Y^{i,t,\xi^i}_s,Z^{i,t,\xi^i}_s)_{s\in[t,T]}$, $i=1,2$, be the solutions to the associated FBSDE systems \eqref{eq:FBSDExi}, with initial data $\xi^1,\xi^2\in L^2(\Om,\cF_{t},\P;\R^{d})$ for the first equation (we do not assume in general that $\cL(\xi^1) = \rho^1_t$ or $\cL(\xi^2) = \rho^2_t$). We define the functions $\varphi_{(\xi^1,\xi^2)} :(t,T)\to \R$, $\Phi_{(\xi^1,\xi^2)}:(t,T)\to[0,+\infty)$ given as 
\begin{equation}\label{def:phi}
\varphi_{(\xi^1,\xi^2)}(\t):=\mathbb E\left[\left(X_\t^{1,t,\xi^1}-X_\t^{2,t,\xi^2}\right) \cdot \left(Y^{1,t,\xi^1}_\t-Y^{2,t,\xi^2}_\t \right) \right]
\end{equation}
and
\begin{equation}\label{def:Phi}
\Phi_{(\xi^1,\xi^2)}(\t) 
:= \mathbb E\left[|X_\t^{1,t,\xi^1}-X_\t^{2,t,\xi^2}|^2\right] +\mathbb E \left[ |Y^{1,t,\xi^1}_\t-Y^{2,t,\xi^2}_\t|^2\right].
\end{equation}
When there is no ambiguity regarding the random variables $\xi^1,\xi^2$, we simply use the notation $\varphi, \Phi$ instead of $\varphi_{(\xi^1,\xi^2)}, \Phi_{(\xi^1,\xi^2)}$.
\end{definition}
These functions will be the essential tools to quantify various decay in time estimates for the MFG system. We can formulate the following result on $\varphi, \Phi$. 

\begin{lemma}\label{lem:Phi_building} Suppose that the assumptions \eqref{hyp:H_reg} and \eqref{hyp:Hdis_strong} take place.
Let $(u^1,\rho^1)$ and $(u^2,\rho^2)$ be two solutions to the mean field game system over $[0,T]$ with initial and final data $(\rho^1_0,g^1)$ and $(\rho^2_0,g^2)$, respectively. Let $t\in [0,T]$ and let $\xi^1,\xi^2\in L^2(\Om,\cF_{t},\P;\R^{d})$ be given and let  $(X^{i,t,\xi^i}_s,Y^{i,t,\xi^i}_s,Z^{i,t,\xi^i}_s)_{s\in[t,T]}$, $i=1,2$, be the solutions to the associated FBSDE systems \eqref{eq:FBSDExi}, with initial data $\xi^1,\xi^2$ and given MFG Nash equilibria and final data $(\rho^{1}_{s})_{s\in[t,T]}$ and $g^{1}$ and $(\rho^{2}_{s})_{s\in[t,T]}$ and $g^{2}$, respectively . Let $\varphi=\varphi_{(\xi^1,\xi^2)}$ and $\Phi=\Phi_{(\xi^1,\xi^2)}$ be defined as in \eqref{def:Phi}. 
\begin{enumerate}
\item Then, for any $s \in [t,T]$ 
 we have
\begin{align*}
c_0 |\varphi(s)| \le \frac {c_0} 2 \Phi(s) \le \varphi'(s) +\frac{C}{2c_0} W_1^2(\rho_s^1,\rho^2_s),
\end{align*}
where $C>0$ depends only on the Lipschitz constant of $D_{x}H$ and $D_{p}H$ in \eqref{hyp:DH_Lip} and $c_0$ is the strong monotonicity constant in \eqref{hyp:Hdis_strong_2}.
\item Assume that $\rho^1_{s} = \rho^2_{s}$ for all $s\in[t,T]$ (i.e. we are considering only one MFG Nash equilibrium, but two different single agent trajectories), or that $\xi^1,\xi^2$ is such that $\cL(\xi^i)=\rho^i_t$, $i=1,2$. Then, for any $s \in [t,T]$ we have
\begin{align*}
2 c_0 |\varphi(s)| \le c_0 \Phi(s) \le \varphi'(s).
\end{align*}
\end{enumerate}
\end{lemma}

\begin{proof}
A direct computation yields
\begin{align*}
 \frac{\dd}{\dd s}\mathbb E & \left[\left(Y^{1,t,\xi^1}_s  -Y^{2,t,\xi^2}_s\right)\cdot\left(X_s^{1,t,\xi^1}-X_s^{2,t,\xi^2}\right) \right]\\
=\mathbb E & \Big[\left(-D_xH(X_s^{1,t,\xi^1},Y^{1,t,\xi^1}_s,\rho_s^1)+D_xH(X_s^{2,t,\xi^2},Y^{2,t,\xi^2}_s,\rho_s^2)\right)\cdot \left(X^{1,t,\xi^1}_s-X^{2,t,\xi^2}_s\right)\Big]\\
&+\mathbb E\Big[\left(Y_s^{1,t,\xi^1}-Y_s^{2,t,\xi^2}\right)\cdot\left(D_pH(X_s^{1,t,\xi^1},Y^{1,t,\xi^1}_s,\rho_s^1)-D_pH(X_s^{2,t,\xi^2},Y^{2,t,\xi^2}_s,\rho_s^2)\right)\Big]\\
=\mathbb E & \Big[\left(-D_xH(X_s^{1,t,\xi^1},Y^{1,t,\xi^1}_s,\rho_s^1)+D_xH(X_s^{2,t,\xi^2},Y^{2,t,\xi^2}_s,\rho_s^1)\right)\cdot \left(X^{1,t,\xi^1}_s-X^{2,t,\xi^2}_s\right)\Big]\\
&+\mathbb E\Big[\left(-D_xH(X_s^{2,t,\xi^2},Y^{2,t,\xi^2}_s,\rho_s^1)+D_xH(X_s^{2,t,\xi^2},Y^{2,t,\xi^2}_s,\rho_s^2)\right)\cdot \left(X^{1,t,\xi^1}_s-X^{2,t,\xi^2}_s\right)\Big]\\
&+\mathbb E\Big[\left(Y_s^{1,t,\xi^1}-Y_s^{2,t,\xi^2}\right)\cdot\left(D_pH(X_s^{1,t,\xi^1},Y^{1,t,\xi^1}_s,\rho_s^1)-D_pH(X_s^{2,t,\xi^2},Y^{2,t,\xi^2}_s,\rho_s^1)\right)\Big]\\
&+\mathbb E\Big[\left(Y_s^{1,t,\xi^1}-Y_s^{2,t,\xi^2}\right)\cdot\left(D_pH(X_s^{2,t,\xi^2},Y^{2,t,\xi^2}_s,\rho_s^1)-D_pH(X_s^{2,t,\xi^2},Y^{2,t,\xi^2}_s,\rho_s^2)\right)\Big]\\
\end{align*}
The assumptions \eqref{hyp:Hdis_strong_2} and \eqref{hyp:DH_Lip}, combined with Young's inequality further imply
\begin{align*}
\varphi'(s) = \frac{\dd}{\dd s}\mathbb E\Big[\big(Y^{1,t,\xi^1}_s-Y^{2,t,\xi^2}_s\big)&\cdot \big(X_s^{1,t,\xi^1}-X_s^{2,t,\xi^2}\big) \Big]\\
\ge & c_0\mathbb E\Big[|X_s^{1,t,\xi^1}-X_s^{2,t,\xi^2}|^2+|Y^{1,t,\xi^1}_s - Y^{2,t,\xi^2}_s|^2\Big]\\
&-\frac{C}{2c_0}W_1^2(\rho_s^1,\rho^2_s)\\
&-\frac{c_0}{2}\mathbb E\Big[|X_s^{1,t,\xi^1}-X_s^{2,t,\xi^2}|^2+|Y^{1,t,\xi^1}_s - Y^{2,t,\xi^2}_s|^2\Big],
\end{align*}
where $C>0$ depends only on the Lipschitz constant of $D_{x}H$ and $D_{p}H$. Rearranging and using \eqref{phivarphi} shows (1). Notice that if $\rho^1_s = \rho^2_s$ for all $s\in [t,T]$, then there is no need to use Young's inequality, and we have directly by \eqref{hyp:Hdis_strong_2} that
\[
 \frac{\dd}{\dd s}\mathbb E \left\{\left[\left(Y^{1,t,\xi^1}_s  -Y^{2,t,\xi^2}_s\right)\cdot\left(X_s^{1,t,\xi^1}-X_s^{2,t,\xi^2}\right) \right]\right\} \ge c_0\E\left[\left|X^{1,t,\xi^1}_s-X^{2,t,\xi^2}_s\right|^2+\left|Y^{1,t,\xi^1}_s-Y^{2,t,\xi^2}_s\right|^2\right].
\]

For the rest of the proof of (2) we argue similarly. We observe that the assumption $\cL(\xi^i)=\rho^i_t$, $i=1,2$ will imply that $\cL(X^{i,t,\xi^i}_s) = \rho^i_s$ for all $s\in (t,T)$. Therefore, in the previous computation we can use the strong monotonicity assumption \eqref{hyp:Hdis_strong} directly. As a consequence, we have 
\begin{align*}
&\frac{\dd}{\dd s}\mathbb E\left[\left(Y^{1,t,\xi^1}_s-Y^{2,t,\xi^2}_s\right)\cdot\left(X_s^{1,t,\xi^1}-X_s^{2,t,\xi^2}\right) \right]\\
&=\mathbb E\Big[\left(-D_xH(X_s^{1,t,\xi^1},Y^{1,t,\xi^1}_s,\rho_s^1)+D_xH(X_s^{2,t,\xi^2},Y^{2,t,\xi^2}_s,\rho_s^2)\right)\cdot \left(X^{1,t,\xi^1}_s-X^{2,t,\xi^2}_s\right)\Big]\\
&+\mathbb E\Big[\left(Y_s^{1,t,\xi^1}-Y_s^{2,t,\xi^2}\right)\cdot\left(D_pH(X_s^{1,t,\xi^1},Y^{1,t,\xi^1}_s,\rho_s^1)-D_pH(X_s^{2,t,\xi^2},Y^{2,t,\xi^2}_s,\rho_s^2)\right)\Big]\\
&\ge c_0\E\Big[|X^{1,t,\xi^1}_s-X^{2,t,\xi^2}_s|^2+|Y^{1,t,\xi^1}_s-Y^{2,t,\xi^2}_s|^2\Big].
\end{align*}
%
\end{proof}

\begin{lemma}\label{lem:integral_ineq}
Let $t\in[0,T]$, let $(\rho^i_s)_{s\in[t,T]}$ ($i=1,2$) be two given flows of probability measures, let $\xi^1,\xi^2\in L^2(\Om,\cF_{t},\P;\R^{d})$ and let $(X^{i}_s,Y^i_s,Z^i_s)_{s\in [t,T]}=(X^{i,t,\xi^i}_s,Y^{i,t,\xi^i}_s,Z^{i,t,\xi^i}_s)_{s\in[t,T]}$ ($i=1,2$) stand for the corresponding solutions to the FBSDE system \eqref{eq:FBSDExi}, with $X^{i,t,\xi^i}_t=\xi^i$, where $(\rho^i_s)_{s\in[t,T]}$ are given. Suppose that $H$ satisfies \eqref{hyp:DH_Lip}.

Then there exists $C>0$ depending on $H$ such that for any $s_1,s_2\in [t,T]$, $s_1<s_2$ we have
\begin{itemize}
\item[(i)]
$
\ds\mathbb E\left[|X^1_{s_2}-X^2_{s_2}|^2\right] \le   \E \left[|X^1_{s_1}-X^2_{s_1}|^2\right] + C\int_{s_1}^{s_2}\left\{\mathbb E\left[ |X_\t^1-X_\t^2|^2+|Y^1_\t-Y^2_\t|^2\right] + W_2^2(\rho^1_\t,\rho^2_\t)\right\} \dd \t,
$
\medskip
\item[(ii)]
$
\ds\mathbb E\left[|X^1_{s_1}-X^2_{s_1}|^2\right] \le   \E\left[|X^1_{s_2}-X^2_{s_2}|^2\right] + C\int_{s_1}^{s_2}\left\{\mathbb E\left[|X_\t^1-X_\t^2|^2+|Y^1_\t-Y^2_\t|^2\right]+W_2^2(\rho^1_\t,\rho^2_\t)\right\}\dd \t,
$
\item[(iii)]
$
\ds\mathbb E\left[|Y^1_{s_1}-Y^2_{s_1}|^2\right] \le   \E\left[|Y^1_{s_2}-Y^2_{s_2}|^2\right] + C\int_{s_1}^{s_2}\left\{\mathbb E \left[|X_\t^1-X_\t^2|^2+|Y^1_\t-Y^2_\t|^2\right]+W_2^2(\rho^1_\t,\rho^2_\t)\right\}\dd \t,
$
\end{itemize}
\end{lemma}

\begin{proof}
Taking the difference of the equations from \eqref{eq:FBSDExi} written for $X^1$ and $X^2$, we compute $$\frac{\dd}{\dd s}\frac12\E\left[ |X^1_s-X^2_s|^2\right].$$ Then integrating the obtained expression from $s=s_1$ to $s=s_2$ we obtain
\begin{align*}
\mathbb E\left[|X^1_{s_2}-X^2_{s_2}|^2\right]  = \E\left[|X^1_{s_1}-X^2_{s_1}|^2\right] +2\int_{s_1}^{s_2}\E\left[ (D_pH(X_\t^1,Y^1_\t,\rho^1_\t)-D_pH(X_\t^2,Y^2_\t,\rho^2_\t))\cdot (X_\t^1-X_\t^2)\right]\dd \t
\end{align*} 
and equivalently
\begin{align*}
\mathbb E\left[|X^1_{s_1}-X^2_{s_1}|^2\right]  = \E\left[|X^1_{s_2}-X^2_{s_2}|^2\right] - 2\int_{s_1}^{s_2}\E\left[ (D_pH(X_\t^1,Y^1_\t,\rho^1_\t)-D_pH(X_\t^2,Y^2_\t,\rho^2_\t))\cdot (X_\t^1-X_\t^2)\right]\dd \t.
\end{align*} 
Taking absolute values of the right hand sides of the previous two equations, using the Lipschitz continuity of $D_pH$ and Young's inequality, we obtain (i) and (ii).

Similar computations for the $Y^i$ variables yield (iii). {Indeed, by It\^o's lemma we have
\begin{align*}
\dd \E\left[|Y^{1}_{s} - Y^{2}_{s}|^{2}\right] & = - 2\E\left[(Y^{1}_{s}-Y^{2}_{s})\cdot(D_{x}H(X^{1}_{s},Y^{1}_{s},\rho^{1}_{s})-D_{x}H(X^{1}_{s},Y^{1}_{s},\rho^{1}_{s}))\right]\dd s\\
&+2\beta\E\left[{\rm trace}\left((Z^{1}_{s}-Z^{2}_{s})^{2}\right)\right] \dd s.
\end{align*}
Integrating this expression between $s_{1}$ and $s_{2}$ we find after rearranging
\begin{align*}
\E\left[|Y^{1}_{s_{1}} - Y^{2}_{s_{1}}|^{2}\right] &= \E\left[|Y^{1}_{s_{2}} - Y^{2}_{s_{2}}|^{2}\right] + \int_{s_1}^{s_2}2\E\left[(Y^{1}_{s}-Y^{2}_{s})\cdot(D_{x}H(X^{1}_{s},Y^{1}_{s},\rho^{1}_{s})-D_{x}H(X^{1}_{s},Y^{1}_{s},\rho^{1}_{s}))\right]\dd s\\
& - 2\beta \int_{s_1}^{s_2} \E\left[{\rm trace}\left((Z^{1}_{s}-Z^{2}_{s})^{2}\right)\right] \dd s\\
&\le \E\left[|Y^{1}_{s_{2}} - Y^{2}_{s_{2}}|^{2}\right] + \int_{s_1}^{s_2}2\E\left[(Y^{1}_{s}-Y^{2}_{s})\cdot(D_{x}H(X^{1}_{s},Y^{1}_{s},\rho^{1}_{s})-D_{x}H(X^{1}_{s},Y^{1}_{s},\rho^{1}_{s}))\right]\dd s,
\end{align*}
where we have used the fact that $(Z^{i}_{s})$, $i=1,2$ are symmetric matrix valued processes and so ${\rm trace}\left((Z^{1}_{s}-Z^{2}_{s})^{2}\right)\ge 0$ almost surely. We conclude by the Lipschitz continuity of $D_{x}H$ and by Young's inequality.
}
\end{proof}

\begin{corollary}\label{cor:integral_ineq}
Under the assumptions of Lemma \ref{lem:integral_ineq} we have the following.
\begin{enumerate}
\item Using the definition of $\Phi$ from \eqref{def:Phi}, we observe that the previous lemma implies that there exists a constant $C>0$ (depending on the data), such that 
\begin{align*}
\Phi(s_1)\le\Phi(s_2)+ C\int_{s_1}^{s_2}\left[ \Phi(s) + W_2^2(\rho^1_s,\rho^2_s)\right]\dd s,
\end{align*}
for any $[s_1,s_2]\subseteq[t,T]$.
\item If in addition we have that $(\rho^i_s)_{s\in [t,T]}$, $i=1,2$ are MFG Nash equilibria and $(X^{i}_{\tau},Y^{i}_{\tau},Z^{i}_{\tau})$, $i=1,2$ are the solutions to the corresponding \eqref{eq:FBSDE}, we in particular have that $\cL({X^i_s})=\rho^i_s$ and so $W_2^2(\rho^1_s,\rho^2_s)\le \E\left[|X^1_{s}-X^2_{s}|^2\right]$, for $s\in [t,T]$. Therefore (by potentially increasing the constants $C>0$) the inequalities in the statement of the lemma further imply 
$$
\ds\mathbb E\left[|X^1_{s_2}-X^2_{s_2}|^2\right] \le \E\left[|X^1_{s_1}-X^2_{s_1}|^2\right] + C\int_{s_1}^{s_2}\mathbb E\left[ |X_\t^1-X_\t^2|^2+|Y^1_\t-Y^2_\t|^2\right]\dd \t,
$$
and similarly for all the {two} other inequalities as well.
\item Under the additional assumptions in (2), we have 
\begin{align*}
\Phi(s_1)\le\Phi(s_2)+ C\int_{s_1}^{s_2}\Phi(s)\dd s,\ \ 
\end{align*}
for any $[s_1,s_2]\subseteq[t,T]$.
\end{enumerate}
\end{corollary}

\begin{proposition}\label{prop:decay_Phi}
Let $(u^1,\rho^1)$ and $(u^2,\rho^2)$ be two solutions to the mean field game system over $[0,T]$ with initial and final data $(\rho^1_0,g^1)$ and $(\rho^2_0,g^2)$, respectively. Suppose that the assumptions of Lemma \ref{lem:Phi_building} are fulfilled.

Let $t\in[0,T]$ and $\xi^1,\xi^2\in L^2(\Om,\cF_{t},\P;\R^{d})$ be given and let $(X^{i,t,\xi^i}_s,Y^{i,t,\xi^i}_s,Z^{i,t,\xi^i}_s)_{s\in[t,T]}$, $i=1,2$ be the solutions to the associated FBSDE systems \eqref{eq:FBSDExi}, with $X^{i,t,\xi^i}_t = \xi^i$. Let $\varphi=\varphi_{(\xi^1,\xi^2)}, \Phi= \Phi_{(\xi^1,\xi^2)}$ defined as in \eqref{def:Phi}.

\begin{enumerate}
\item There exists $C>0$ depending on the data 
such that for any $t \le t_1 < t_2 \le T$
\begin{multline}
\frac {c_0} 2 \int_{t_1}^{t_2} \Phi(s) \dd s \le e^{-c_0 (T-t_2) } {\varphi(T)} {-} e^{-c_0 (t_1-t) } {\varphi(t)}  \\ +C \left(  
 \int_t^{t_1} e^{-c_0 (t_1-s)} W_2^2(\rho_s^1,\rho^2_s)\dd s
 + \int_{t_1}^{t_2} W_2^2(\rho_s^1,\rho^2_s) \dd s 
+  \int_{t_2}^T e^{-c_0 (s-t_2)} W_2^2(\rho_s^1,\rho^2_s) \dd s \right).
\end{multline}
\item Assuming that $\cL({\xi^i})=\rho^i_t$ ($i=1,2$) and that $\cL({X^{i,t,\xi^i}_s})=\rho^i_s$, for all $s\in[t,T]$, $i=1,2$, then
\[
{c_0}  \int_{t_1}^{t_2} \Phi(s) \dd s \le e^{-{2c_0} (T-t_2) } {\varphi(T) -} e^{-{2c_0} (t_1-t) } {\varphi(t)}
\]
\item Assuming $\cL({\xi^i})=\rho^i_t$, ($i=1,2$) or that $\rho^1_{s} = \rho^2_{s},\ \forall s\in[t,T]$, and that $g^1=g^2 = g$ and $g$ is displacement monotone, then
\[
{c_0}  \int_{t_1}^{t_2} \Phi(s) \dd s \le e^{-{2c_0} (t_1-t) } |\varphi(t)|
\]
\end{enumerate}


\end{proposition}

\begin{proof} 

Setting $$h(s) := \frac{C}{2c_0}W_2^2(\rho_s^1,\rho^2_s),$$ where
$C>0$ is given in Lemma \ref{lem:Phi_building} (1). This means that there exists a constant $C>0$, depending only on the data (that we do not relabel) such that
\begin{equation}\label{hestimate}
h(s) {=} CW_2^2(\rho_s^1,\rho^2_s).
\end{equation}

Lemma \ref{lem:Phi_building} (1) yields the following inequalities for any $s \in [t,T]$:
\begin{align}\label{eqPhih}
-c_0 \varphi(s) \le \frac {c_0} 2 \Phi(s) \le \varphi'(s) + h(s), \qquad c_0 \varphi(s) \le \frac {c_0} 2 \Phi(s) \le \varphi'(s) + h(s).
\end{align}
By integrating the first one on $(t, \tau)$ and the second one on $(\tau, T)$ we get
\[
e^{-c_0 (\tau-t) } \varphi(t) - e^{-c_0 \tau } \int_t^\tau e^{c_0 s} h(s) \dd s \le \varphi(\tau) \le e^{-c_0 (T-\tau) } \varphi(T) + e^{c_0 \tau } \int_\tau^T e^{-c_0 s} h(s) \dd s.
\]
Finally, by integrating \eqref{eqPhih} on $(t_1, t_2)$ we obtain
\[
\frac {c_0} 2 \int_{t_1}^{t_2} \Phi(s) \dd s \le \varphi(t_2)-\varphi(t_1) + \int_{t_1}^{t_2} h(s) \dd s.
\]
The two previous inequalities together with \eqref{hestimate} show (1).

The proof of (2) is almost identical, noticing that Lemma \ref{lem:Phi_building} (2) provides the same kind of inequality with $h \equiv 0$.

To prove (3) one may argue in the same way as for (2), noticing in addition that integration on $(t_1, t_2)$ gives ${c_0}  \int_{t_1}^{t_2} \Phi(s) \dd s \le \varphi(t_2)-\varphi(t_1)$, and the inequality  
$$
\varphi(t_2) =\mathbb E\left[(X_{t_2}^{1,t,\xi^1}-X_{t_2}^{2,t,\xi^2}) \cdot \left(Y^{1,t,\xi^1}_{t_2}-Y^{2,t,\xi^2}_{t_2} \right) \right] \le 0
$$
can be used in view of the propagation of displacement monotonicity (see \cite[Theorem 4.2]{MesMou}), which applies in view of \eqref{hyp:Hdis_strong} if $\cL(\xi^i)=\rho^i_t$, or in view of \eqref{hyp:Hdis_strong_2} if $\rho^1_{s}= \rho^2_{s}$, for all $s\in[t,T]$. {In particular, $\varphi(T)\le 0$.}
\end{proof}

We now show that second-moment bounds for the MFG Nash equilibrium can be transferred to analogous bounds on solutions to the FBSDE system \eqref{eq:FBSDExi}.

\begin{proposition}\label{prop:2nd_local}
Suppose that all the assumptions  \eqref{hyp:H_reg}--\eqref{hyp:g_confining} take place.
Let $(\rho_s)_{s\in[t,T]}$ be an MFG Nash equilibrium and let $\xi \in L^2(\Om,\cF_{t},\P;\R^{d})$ be given. Let $(X^{t,\xi}_s,Y^{t,\xi}_s,Z^{t,\xi}_s)_{s\in[t,T]}$ be the solution to \eqref{eq:FBSDExi}. Then, there exists a constant $C>0$ depending on the data (in particular also on $M_2(\rho_0)$), but independent of $T$ such that
\begin{align*}
\E\left[|X^{t,\xi}_s|^2\right] \le C(1+\E\left[|\xi|^2\right]),\ \forall s\in [t,T],
\end{align*}
and
\begin{align*}
\E\left[|Y^{t,\xi}_s|^2\right] \le C(1+\E\left[|\xi|^2\right]),\ \forall s\in [t,T].
\end{align*}
\end{proposition}

\begin{proof}
Let $(X^{1}_s,Y^{1}_s,Z^{1}_s)_{s\in[t,T]}=(X^{t,\xi}_s,Y^{t,\xi}_s,Z^{t,\xi}_s)_{s\in[t,T]}$ be the solution to \eqref{eq:FBSDExi}, while let us set $(X^{2}_s,Y^{2}_s,Z^{2}_s)_{s\in[t,T]} = (X^{t,\hat{\xi}}_s,Y^{t,\hat{\xi}}_s,Z^{t,\hat{\xi}}_s)_{s\in[t,T]}$ to be the Nash equilibrium, that is, $\mathcal L(\hat \xi) = \rho_t$. Note that the two FBSDE systems share the same measure flow $\rho_s = \rho^1_s = \rho^2_s$ and the same final datum $g = g^1 = g^2$. Hence, with $\varphi$ and $\Phi$ defined as above (see Definition \ref{def:phi_fcts}) we first have by Lemma \ref{lem:integral_ineq}(i) that for $s \ge t$
\begin{equation}\label{eq:2345}
\ds\mathbb E\left[|X^1_{s}-X^2_{s}|^2\right] \le   \E \left[|X^1_{t}-X^2_{t}|^2\right] + C\int_{t}^{s} \Phi(\t )\dd \t =  \E \left[|\xi - \hat \xi |^2\right] + C\int_{t}^{s}   \Phi(\t )\dd \t.
\end{equation}
Now, integrating the inequality of Lemma \ref{lem:Phi_building}(2) with respect to $\tau$ on $[t,T]$ gives
\begin{multline*}
c_0 \int_t^T \Phi(\tau) \dd \t \le \varphi(T) - \varphi(t) = \\
- \E\left[\left(X_T^{1}-X_T^{2}\right) \cdot \left(D_x g(X_T^{1}, \rho_T) - D_x g(X_T^{2}, \rho_T)\right) \right] + \E\left[\left(\xi- \hat \xi\right) \cdot \left(D_x u(t, \xi) - D_x u(t, \hat \xi)\right) \right].
\end{multline*}
While the first term is non-positive by \eqref{hyp:g_D-mon} (note that we are not using the full D-monotonicity, rather the convexity of $g(\cdot,\mu)$ for each measure $\mu$), the second one can be bounded by the uniform in time estimate on $|D^2_{xx} u|$ in \eqref{eq:Cu} , that is
\[
c_0 \int_t^T \Phi(\tau) \dd \t \le C_u \E\left[ |\xi - \hat \xi|^2\right].
\]
Plugging this back into \eqref{eq:2345} yields
\[
\E\left[|X^1_{s}-X^2_{s}|^2\right] \le C\E \left[|\xi - \hat \xi |^2\right] \le 2C  \E \left[|\xi|^2\right] +2C  \E \left[|\hat \xi|^2\right] .
\]
Together with the uniform bounds on $\E[|\hat \xi|^2]$ from Proposition \ref{prop:2nd_moment_bound} we then estimate $\E[X^1_{s}|^2$. Finally, again in view of the estimates on $|D^2_{xx} u|$,
\[
\E\left[|Y^1_{s}-Y^2_{s}|^2\right] = \E\left[|D_x u(s,X^1_{s})-D_x u(s,X^2_{s})|^2\right] \le C_u^2 \E\left[|X^1_{s}-X^2_{s}|^2\right]
\]
and the bound on $\E[Y^1_{s}|^2$ follows by Proposition \ref{prop:L2_Y}.
\end{proof}

\begin{lemma}
Recall the notations from Definition \ref{def:phi_fcts}. Suppose that the assumptions of Proposition \ref{prop:2nd_local} take place. Then there exists a constant $C>0$ depending on the data (but independent of $T$) such that 
\begin{equation}\label{phivarphiest}
\E\left[|X^{i,t,\xi^i}_\t|^2\right],  \, \E\left[|Y^{i,t,\xi^i}_\t|^2\right], \, |\varphi(\t)| , \, \Phi(\tau) \le C\left(1+\E\left[|\xi^1|^2+|\xi^2|^2\right]\right), \qquad \forall \t \in [t,T], \, i =1,2.
\end{equation}
\end{lemma}

\begin{proof}
By Young's inequality we have
\begin{equation}\label{phivarphi}
|\varphi(\t)| \le \frac 12 \Phi(\tau) \qquad \forall \t \in [t,T].
\end{equation}

Under the assumptions of Proposition \ref{prop:2nd_local}, we notice that $\E \left[|X^{i,t,\xi^1}_\t|^2+|Y^{i,t,\xi^1}_\t|^2\right]$ is uniformly bounded on $[t,T]$, by a constant depending only on the data and $\E\left[|\xi^{1}|^{2}+|\xi^{2}|^{2}\right]$, and this dependence is precisely in the form stated in \eqref{phivarphiest}. The result follows.

\end{proof}


\subsection{Decay results in time for $\rho$ and for $D_xu$}

We are now ready to prove the first main result of this work, on the distance between two Nash equilibria for large $T$. Throughout the rest of this section, we assume that  \eqref{hyp:H_reg}--\eqref{hyp:g_confining} are in force.

\begin{theorem}\label{prop:pointwisedecay}
Let $\left(u^1,\rho^1\right)$ and $\left(u^{2},\rho^{2}\right)$ be the unique solutions to two MFG systems with the same Hamiltonian and final/initial data $(g^1,\rho^1_0)$ and $(g^2,\rho^2_0)$, respectively, on a given time interval $[0,T]$. 

Let $\xi^1,\xi^2\in L^2(\Om,\cF_{0},\P;\R^{d})$ be given and let $(X^{i,0,\xi^i}_s,Y^{i,0,\xi^i}_s,Z^{i,0,\xi^i}_s)_{s\in[0,T]}$, $i=1,2$ be the solutions to the associated FBSDE systems \eqref{eq:FBSDE}, with $X^{i,0,\xi^i}_0 = \xi^i$. In particular, we suppose that $\cL(\xi^i) = \rho^i_0$ and $\cL(X^{i,0,\xi^i}_s) = \rho^i_s$ for all $s\in[0,T]$, $i=1,2$. Let $\varphi=\varphi_{(\xi^1,\xi^2)}, \Phi=\Phi_{(\xi^1,\xi^2)}$ defined as in \eqref{def:Phi}. 

Then we have the following. 

\begin{enumerate}
\item There exists $C > 0$ such that for any $t \in [0,T]$,
\[
\Phi(t) \le C\left(e^{-{2c_0} (T-t)} + e^{-{2c_0} t}\right).
\]
\item Assuming in addition that $\rho^1_0 = \rho^2_0$, for any $t \in [0,T]$ we have
\[
\Phi(t) \le Ce^{-{2c_0} (T-t)}.
\]
\item Assuming otherwise that $g^1=g^2 = g$ and $g$ is displacement monotone, for any $t \in [0,T]$ we have
\[
\Phi(t) \le C e^{-{2c_0} t}.
\]
\end{enumerate}
The constants $C$ above depend on the data and $\E\left[|\xi^1|^2+|\xi^2|^2\right]$.
\end{theorem}

\begin{proof} 

{For (1), suppose first that $t\in[0,T-1]$. }
We start from the inequality in Proposition \ref{prop:decay_Phi} (2) (applied to $t=0$, {$t_{1} = t$ and $t_{2}= t + 1$}), i.e.
\[
{c_0}  {\int_{t}^{t+1}} \Phi(s) \dd s \le e^{-{2c_0} (T - {t-1}) } |\varphi(T)| + e^{-{2c_0}{t} } |\varphi(0)|.
\]
Hence, in view of bounds \eqref{phivarphiest} 
 and the Mean Value Theorem, there exists $\zeta_t \in {[t,t+1]}$ and a constant $C$ (depending on the data and $\E\left[|\xi^1|^2+|\xi^2|^2\right]$, but independent of $T$) such that
\[
\Phi(\zeta_t) \le C\left( { e^{-2c_0 (T-t)} + e^{-2c_{0}t}}\right).
\]
Therefore, by Corollary \ref{cor:integral_ineq} (3) we have
\begin{align*}
\Phi(t)&\le\Phi(\zeta_t)+ C\int_{t}^{\zeta_t}\Phi(s)\dd s 
\le \Phi(\zeta_t)+ C\int_{t}^{ t+1}\Phi(s)\dd s\\
&\le C(e^{-{c_0} (T-t)} + e^{-{2c_0} t}).
\end{align*}

{If now $t\in [T-1,T]$, we simply have by Corollary \ref{cor:integral_ineq} (3) that
\begin{align*}
\Phi(t)\le \Phi(T) + C\int_{t}^{T}\Phi(s)\dd s \le \Phi(T) + C (|\varphi(T)|+e^{-2c_{0}t}|\varphi(0)|),
\end{align*}
and by \eqref{phivarphiest} we conclude the proof of (1).
}

To obtain (2), we note that the additional assumption implies $\varphi(0) = 0$. {First, here also suppose that $t\in[0,T-1]$.}

Hence, Proposition \ref{prop:decay_Phi} (2) gives
\[
{c_0}  {\int_{t}^{t+1}} \Phi(s) \dd s \le e^{-{2c_0} (T-t-1) } |\varphi(T)| 
\]
for any $t \in [0,{T-1}]$. We can pick now {$\zeta_t \in [t,t+1]$} such that $\Phi(\zeta_t) =  {\int_{t}^{t+1}} \Phi(s)\dd s \le C e^{-2c_0 {(T-t)}}$. By Corollary \ref{cor:integral_ineq} (3), 
\[
\Phi(t)\le\Phi(\zeta_t)+ C{\int_{t}^{\zeta_t}}\Phi(s) \dd s \le \Phi(\zeta_t)+ C\int_{t}^{t+1}\Phi(s) \dd s \le Ce^{-2c_0 (T-t)},
\]
and the inequality is easily extended to any {$t \in [T-1,T]$, just as in the case of (1)}.

To prove (3) we can argue as before. { Again, suppose first that $t\in [0,T-1]$.}

Starting from Proposition \ref{prop:decay_Phi} (3) {(applied for $t:=0$, $t_{1}:=t$ and $t_{2}:=t+1$)} and the existence of {$\zeta_t \in [t, t+1]$} such that
\[
c_0 \Phi(\zeta_t) = {c_0}  {\int_{t}^{t+1} } \Phi(s) \dd s \le e^{-{2c_0}{t} } |\varphi(0)| \le C e^{-2c_0 {t}}.
\]
{Then, we have by Corollary \ref{cor:integral_ineq} (3) that
$$
\Phi(t)\le \Phi(\zeta_{t})+\int_{t}^{\zeta_t}\Phi(s)\dd s \le \Phi(\zeta_{t})+\int_{t}^{t+1}\Phi(s)\dd s \le C e^{-2c_0 t}.
$$
Suppose now that $t\in [T-1,T]$. Using now Proposition \ref{prop:decay_Phi} (3) (applied for $t:=0$, $t_{1}:=T-1$ and $t_{2}:=T$), we find that there exists $\zeta_T\in [T-1,T]$ such that
\[
c_0 \Phi(\zeta_T) = {c_0} \int_{T-1}^{T}  \Phi(s) \dd s \le e^{-2c_0(T-1)}  |\varphi(0)| \le C e^{-2c_0 T}.
\]
By Corollary \ref{cor:integral_ineq} (2) we have for any $s\in[T-1,T]$ that
\begin{align}\label{new:1}
\E\left[\left|X^{1,0,\xi^{1}}_{s}-X^{2,0,\xi^{2}}_{s}\right|^{2} \right] \le \Phi(\zeta_T) + \int_{\min\{s,\zeta_T\}}^{\max\{s,\zeta_T\}}\Phi(s)\dd s \le \Phi(\zeta_T) + \int_{T-1}^{T}\Phi(s)\dd s \le C e^{-2c_0 T}.
\end{align}
Then, we have
\begin{align}\label{new:2}
\nonumber\E\left[|Y^{1,0,\xi^{1}}_{t}-Y^{2,0,\xi^{2}}_{t}|^{2} \right] & \le \E\left[\left|D_x g(X^{1,0,\xi^{1}}_{T},\rho^{1}_{T})-D_x g(X^{2,0,\xi^{2}}_{T},\rho^{2}_{T})\right|^{2}\right] + C\int_{t}^{T}\Phi(s)\dd s\\
\nonumber& \le C \E\left[\left|X^{1,0,\xi^{1}}_{T}-X^{2,0,\xi^{2}}_{T}\right|^{2}\right] + C\int_{T-1}^{T}\Phi(s)\dd s\\
& \le C e^{-2c_0 T}.
\end{align}
Combining the two inequalities \eqref{new:1} and \eqref{new:2}, we conclude by the last point also for the range $t\in[T-1,T]$.
}


\end{proof}

\begin{remark}\label{rmk:nonstrict} The previous result says that different equilibria have to ``collapse'' as $T \to \infty$, that is: two equilibria starting from different initial conditions $\rho^1_0,\rho^2_0$ and approaching different final values $g^1, g^2$ must be exponentially close in the sense of Theorem \ref{prop:pointwisedecay}(1). This is a consequence of strong D-monotonicity of the data. We show below that the same conclusion may not hold if one has only D-monotonicity. This is in contrast with the Lasry--Lions monotone setting, where it is known (see for instance \cite{CirPor}) that the presence of the diffusion guarantees this kind of long time stability even {in} the presence of some mild anti-monotonicity. A similar compensation phenomenon does not appear in the D-monotone setting. Consider indeed the following system 
\[
\left\{
\begin{array}{ll}
\ds-\partial_t u  -\b\partial_{xx}^2 u + \frac12 |\partial_{x}u|^2 = \left(x- \int_{\R}y \rho(t,y)\dd y\right)^2, & {\rm{in}\ } (0,T)\times \R,\\[5pt]
\partial_t\rho -\b\partial_{xx}^2 \rho - \nabla\cdot (\rho \partial_{x}u) = 0, & {\rm{in\ }} (0,T)\times \R,\\[5pt]
\end{array}
\right.
\]
which is solved, for all $\beta > 0$ and $T > 0$, by the couple ($u, \rho$), where
\[
u(t,x) = \frac{\sqrt 2}2 x^2 -\beta (t-T) \sqrt 2 ,
\]
and $\rho$ is a normal distribution with zero mean and variance $\beta / \sqrt 2$. Note that the Hamiltonian is displacement monotone, but not strongly displacement monotone, in the sense that it satisfies \eqref{hyp:Hdis_strong} with $c_0 = 0$. Note also that any space translation $(u(t, \cdot + z), \rho(t, \cdot + z))$, for $z \in \R$ solves the same system of PDE, so we have a continuum of stationary equilibria (with the corresponding final condition $g(x, \rho) = \frac{\sqrt 2}2 (x+z)^2$, which satisfies our standing assumptions, and is in fact strongly D-monotone). We note also that for such MFG systems the measure component is always stationary, Gaussian centered at $-z$. Therefore, it is clear that any two such equilibria do not satisfy the stability property, as their (positive) $W_2$-distance (in the sense quantified by $\Phi(t)$) remains constant in $t$, uniformly with respect to the time horizon $T$.
\end{remark}

\begin{remark}\label{rmk:decay} Let $(X^{i,t,\xi^i}_s,Y^{i,t,\xi^i}_s,Z^{i,t,\xi^i}_s)_{s\in[t,T]}$, $i=1,2$ be solutions to the FBSDE system \eqref{eq:FBSDExi}, with $X^{i,t,\xi^i}_t = \xi^i\in L^{2}(\Om,\cF_{0},\P;\R^{d})$, with the same final data $g^{1}=g^{2}=g$ and with the same input flows of measures $\rho^{1}_{s}=\rho^{2}_{s}=\rho_{s}$ for $s\in[t,T]$, where $(\rho_{s})_{s\in[t,T]}$ is an MFG Nash equilibrium. It is immediate to see that the conclusion of Theorem \ref{prop:pointwisedecay} (3) holds true in this case, i.e. 
\[
\Phi(s) \le C e^{-{2c_0} (s-t)},\ s\in (t,T),
\]
where the constant $C>0$ depends on the data and $\E\left[|\xi^1|^2+|\xi^2|^2\right]$. Indeed, one may argue as in the previous proof, and use Proposition \ref{prop:decay_Phi} (3) and Corollary \ref{cor:integral_ineq} (1).

\end{remark}

\begin{corollary}\label{newcor}
As a consequence of Theorem \ref{prop:pointwisedecay}, under our standing assumptions we have that if $(\rho^i_t)_{t\in[0,T]}$, $i=1,2$ are two MFG Nash equilibria, then 
\begin{enumerate}
\item $W_2^2(\rho^1_t,\rho^2_t) \le C(e^{-{2c_0} (T-t)} + e^{-{2c_0} t}),$ for all $t\in[0,T]$.
\smallskip
\item If $\rho^1_0=\rho^2_0$, then $W_2^2(\rho^1_t,\rho^2_t) \le Ce^{-{2c_0} (T-t)},$ for all $t\in[0,T]$.
\smallskip
\item If $g^1=g^2$, then $W_2^2(\rho^1_t,\rho^2_t) \le C e^{-{2c_0} t},$ for all $t\in[0,T]$,
\end{enumerate}
where the constants $C>0$ above depend on the data and $M_2(\rho^1_0)+M_2(\rho^2_0).$
\end{corollary}

\subsubsection{Localization arguments for $D_x u$} We can see that Theorem \ref{prop:pointwisedecay} gives decay result for $D_xu$ only along the flow, i.e. for $D_x u (t,X^{0,\xi}_t) = -Y^{0,\xi}_t$, where $(X^{0,\xi}_t,Y^{0,\xi}_t,Z^{0,\xi}_t)_{t\in[0,T]}$ corresponds precisely to MFG Nash equilibria. In order to have a result which is pointwise in the spatial variable, we need to have some additional work. In particular, we need to `localize' the initial data $\xi^i$ in the FBSDE system, and hence consider \eqref{eq:FBSDExi} instead of \eqref{eq:FBSDE}, and like that the resulting flow is such that $\cL({X^{i,t,\xi^i}_{s}})$ does not correspond to $\rho^i_s$. 

We can now formulate the following result.

\begin{theorem}\label{prop:Du_loc}
Let $\left(u^1,\rho^1\right)$ and $\left(u^2,\rho^2\right)$ be the unique solutions to two mean field games systems with the same Hamiltonian and final/initial data $(g^1,\rho_0)$ and $(g^2,\rho_0)$, respectively, on a given time interval $[0,T]$. Then there exists $C>0$ depending on the data (and in particular on $M_2(\rho_0)$, but independent of $T$) such that, for every $t \in [0,T]$,
\begin{equation}\label{ineq:D_xu_final_ex}
\sup_{x\in\R^d}\frac{|D_xu^{1}(t,x)-D_xu^{2}(t,x)|^2}{1+|x|^2} \le C e^{-c_0 {(T-t)}}.
\end{equation}

\end{theorem}

\begin{proof}
First, let us notice that as $\rho^1_0=\rho^2_0=\rho_0$, Theorem \ref{prop:pointwisedecay} (2) implies that
\begin{align}\label{eq:W2_xx}
W_2^2(\rho^{1}_s,\rho^{2}_s)\le C_0 e^{-{2c_0} (T-s)},\ \  {\rm{for\ all}}\ s\in[0,T]
\end{align}
for some $C_0>0$ that depends on the second moment of $\rho_0$ (and on the data). Let $t \in [0, T/4]$ and let $\xi\in L^2(\Om,\cF_{t},\P;\R^{d})$ be a random variable. We are going to restrict the two MFGs to the time interval $[t,T]$. In particular, we consider the two FBSDE systems \eqref{eq:FBSDExi} with $\xi^1 = \xi^2 = \xi$ and $g^1, g^2$ as final conditions. 

{ Suppose first that $t\in [0,T-1]$.} According to Proposition \ref{prop:decay_Phi} (1), there exists a positive constant $C_\xi>0$ such that
\begin{multline*}
\frac {c_0} 2 \int_{{ t}}^{{ t+1}} \Phi(s) \dd s \le e^{-c_0 ({ T-t-1}) } |\varphi(T)| +  |\varphi(t)|  \\ +C_\xi \left(  
 \int_{{ t}}^{{ t+1}} W_2^2(\rho_s^1,\rho^2_s) \dd s 
 + \int_{{ t+1}}^T e^{-c_0 ({ s-t-1})} W_2^2(\rho_s^1,\rho^2_s) \dd s \right).
\end{multline*}
Since $\varphi(t) = 0$, by the estimate \eqref{eq:W2_xx} on $W_2^2(\rho^{1}_s,\rho^{2}_s)$, the uniform bounds on $|\varphi(T)|$ and the Mean Value Theorem we get
\begin{multline*}
\Phi(\zeta_t) = \int_{{ t}}^{{ t+1}} \Phi(s) \dd s \le C_\xi \Big(e^{-c_0 { (T-t)}} + \int_{{ t}}^{{ t+1}} e^{-2c_0 (T-s)} \dd s + \int_{{ t+1}}^T {  e^{-2c_0 (T-s)-c_0 (s-t-1)} } \dd s\Big) \\ \le C_\xi e^{-c_0 {(T-t)}}
\end{multline*}
for some $\zeta_t \in [{ t, t+1}]$. Note that $C_\xi$ may vary from line to line, but it always depends in an affine way on $\E \left[|\xi|^2\right]$. 
Let us now apply Corollary \ref{cor:integral_ineq} (1) and plug in the estimates obtained so far:
\[
\Phi(t)\le\Phi(\zeta_t)+ C_1\int_{t}^{\zeta_t}\left[ \Phi(s) + W_2^2(\rho^1_s,\rho^2_s)\right]\dd s \le C_\xi e^{-c_0 { (T-t)}}.
\]
For any $x \in \R^d$, choosing $\cL(\xi)=\d_{x}$ yields $\Phi(t) = \mathbb E \left[|Y^{1,t,\xi}_t-Y^{2,t,\xi}_t|^2\right] = |D_xu^{1}(t,x)-D_xu^{2}(t,x)|^2$. Since $C_\xi$ is an affine function of $|x|^2$, we obtain the assertion.

{ If $t \in [T-1, T]$, then estimate is a straightforward consequence of \eqref{gradbound}. }
\end{proof}



\section{Asymptotic behavior of the value function}\label{sec:u}

In this section our goal is to study the long time behavior of the value function. For this reason we consider the following objects. For a fixed time horizon $T>0$,  $g:\R^d\times\sP_2(\R^d)\to\R$ given final condition and $\rho_0\in\sP_2(\R^d)$ initial measure, let $(\rho^T_s)_{s\in(0,T)}$ stand for the MFG Nash equilibrium. In particular, if $(X^{0,\xi}_s,Y^{0,\xi}_s, Z^{0,\xi}_s)_{s\in(0,T)}$ is the solution to \eqref{eq:FBSDE} with $\xi\in L^2(\Om,\cF_{0},\P;\R^{d})$ such that $\cL(\xi) = \rho_0$, then we have $\rho^T_s=\cL({X^{0,\xi}_s})$ for all $s\in[0,T]$.

In a similar way, we define the associated value function $u^T:[0,T]\times\R^d\to\R$ as
\begin{equation*}
u^T(t,x):=\mathbb{E}\left\{\int_t^T L(X^{t,x}_s,D_pH(X^{t,x}_s,Y^{t,x}_s,\rho^T_s), \rho^T_s)\dd s + g(X^{t,x}_T,\rho^T_T)\right\},
\end{equation*}
where $(X^{t,x}_s,Y^{t,x}_s, Z^{t,x}_s)_{s\in(t,T)}$ is the solution to \eqref{eq:FBSDExi} with $\cL(\xi) = \d_x$, and in the definition of $u^T$ the Nash equilibrium $(\rho^T_s)_{s\in(0,T)}$ has been used.

{Here the Lagrangian $L:\R^{d}\times\R^{d}\times\sP_{2}(\R^{d})\to\R$ and Hamiltonian $H:\R^{d}\times\R^{d}\times\sP_{2}(\R^{d})\to\R$ are Legendre duals of each other with respect to the second variable, i.e.
$$
L(x,v,\mu) = \sup_{p\in\R^{d}}\left\{v\cdot p - H(x,p,\mu)\right\}.
$$
Using classical results from convex analysis, since $H$ is uniformly convex in the $p$ variable by \eqref{hyp:Hdis_strong_2},
the following identities hold:
\begin{align*}
D_{p}H(x,D_{v}L(x,v,\mu),\mu) = v,\ \ D_{v}L(x,D_{p}H(x,p,\mu),\mu) = p
\end{align*}
and
\begin{align*}
H(x,p,\mu) = p\cdot D_{p}H(x,p,\mu) - L(x,D_{p}H(x,p,\mu),\mu).
\end{align*}
These allow to transfer several properties of $H$ to $L$. In particular, we will need a suitable Lipschitz regularity of $L$, as described in Remark \ref{remL1L2} below.

The couple $(u^T,\rho^T)$ solves the MFG system

\begin{align*}
\left\{
\begin{array}{ll}
-\partial_t u^T  -\b\Delta u^T + H(x,-Du^T,\rho^T) = 0, & \text{in } (0,T)\times \R^d,\\
\partial_t\rho^T -\b\Delta \rho^T + \nabla\cdot (\rho^T D_pH(x,-Du^T,\rho^T)) = 0, & \text{in } (0,T)\times \R^d,\\
\rho^T_0 = \rho_0; \ u^T(T,\cdot) = g, & \text{in } \R^d.
\end{array}
\right.
\end{align*}

\medskip

In this section we need to assume further hypotheses on $H$ and $g$. We collect these ones here.

\begin{align}\label{hyp:add_H}\tag{H\arabic{hyp}}
\exists C>0,\ 
\left|D_\mu H\left(x, p,\mu,\tilde x \right)\right|\le C\left(1+ |\tilde x| + |x| + |p| \right),\  \forall (x,p,\mu,\tilde x)\in \R^{d}\times\R^{d}\times\sP_{2}(\R^{d})\times\R^{d}.
\end{align}
\stepcounter{hyp}

\begin{align}\label{hyp:add_g}\tag{H\arabic{hyp}}
\exists C>0,& \ \ {\rm{such\ that}}\\
&\bullet \left|g(x,\mu)\right|\le C\left(1+ M_{2}^{2}(\mu) + |x|^{2} \right),\ \ \forall (x,\mu)\in \R^{d}\times\sP_{2}(\R^{d});\\
&\bullet |D_\mu g\left(x,\mu,\tilde x \right)|\le C\left(1+ |\tilde x| + |x|\right),\  \forall (x,\mu,\tilde x)\in \R^{d}\times\sP_{2}(\R^{d})\times\R^{d}.
\end{align}
\stepcounter{hyp}

\begin{remark}
The inequality in \eqref{hyp:add_H} could be achieved if in addition to the assumptions in \eqref{hyp:H_reg} one would impose uniform bounds on $D^{2}_{\mu\mu}H$ and $D^{2}_{\mu\tilde x}H$.

Similarly, a sufficient condition for the fulfillment of \eqref{hyp:add_g}, in addition to \eqref{hyp:g}, would be uniform bounds on $D^{2}_{\mu\mu}g$ and $D^{2}_{\mu\tilde x}g$. Both of these additional hypotheses are slightly weaker than imposing additional regularity and derivative bounds on the data.
\end{remark}

\medskip

\begin{remark}\label{remL1L2} In what follows, it will be frequent to estimate the difference
\[
\left| \E \int_{t_1}^{t_2} L(X^1_s,D_pH(X^1_s,Y^1_s,\rho^1_s), \rho^1_s)\dd s - \E \int_{t_1}^{t_2} L(X^2_s,D_pH(X^2_s,Y^2_s,\rho^2_s), \rho^2_s)\dd s\right|
\]
where $(\rho^i_s)_{s\in[0,T]}$ and the associated $(X^i_s,Y^i_s, Z^i_s)_{s\in[0,T]}$ are MFG Nash equilibria, in terms of the function $\Phi$ defined in \eqref{def:Phi}, i.e.
\[
\Phi(s) = \E\left[|X_s^{1}-X_s^{2}|^2\right] +\mathbb E\left[ |Y^{1}_s-Y^{2}_s|^2\right],
\]
having the information that $X^i_s$, $Y^i_s$ enjoy universal second moment bounds. By employing the Legendre--Fenchel duality
$$L(x,D_pH(x,p,\rho),\rho)=  p\cdot D_pH(x,p,\rho) - H(x,p,\rho),\ \ \forall (x,p,\rho)\in\R^d\times\R^d\times\sP_2(\R^d),$$ 
and the fact that $W^2_2(\rho^1_s, \rho^2_s)\le  \E \left[| X^1_s - X^2_s|^2\right]$, one checks that assumptions \eqref{hyp:H_reg} guarantee
\begin{multline}\label{L1L2bound}
\left| \E \int_{t_1}^{t_2} L(X^1_s,D_pH(X^1_s,Y^1_s,\rho^1_s), \rho^1_s)\dd s - \E \int_{t_1}^{t_2} L(X^2_s,D_pH(X^2_s,Y^2_s,\rho^2_s), \rho^2_s)\dd s\right| \\ \le  C\int_{t_1}^{t_2} \sqrt{\Phi(s)} \dd s,
\end{multline}

where $C > 0$ depends only on $H$ and on the second moment bounds on $X^i$, $Y^i$. For the sake of completeness we prove this fact in the following lemma.
\end{remark}

\begin{lemma}
Suppose that $H$ satisfies the regularity assumptions in \eqref{hyp:H_reg} and \eqref{hyp:add_H}. Then \eqref{L1L2bound} holds true for a constant $C>0$ depending only on second moment bounds and on $H$.
\end{lemma}

\begin{proof}
Let $(x^i,p^i,\rho^i)\in \R^d\times\R^d\times\sP_2(\R^d),$ $i=1,2$ be given. First let us estimate
\begin{align*}
\left|H(x^2,p^2,\rho^2) - H(x^1,p^1,\rho^1) \right|.
\end{align*}
For this, let $(\rho_s)_{s\in[0,1]}$ be a $W_2$-geodesic connecting $\rho^1$ to $\rho^2$, so, in particular $\rho_0 = \rho^1$ and $\rho_1 = \rho^2$. 

In particular, we recall that there exists a family of Borel vector fields $(v_t)_{s\in[0,1]}$ such that $\partial_s\rho_s +\nabla\cdot (v_s\rho_s) = 0$ in the sense of distributions, $W_2^2(\rho^1,\rho^2) = \int_0^1\int_{\R^d}|v_s|^2\dd\rho_s\dd s$ and for any $f:\sP_2(\R^d)\to\R$, $W_2$-differentiable, we have that 
\begin{align*}
\left|f(\rho^2) - f(\rho^1) \right|&= \left|\int_0^1\int_{\R^d}D_\mu f(\rho_s)(\tilde x)\cdot v_s(\tilde x)\dd\rho_s(\tilde x)\dd s\right|\\
& \le \left(\int_0^1\int_{\R^d}|D_\mu f(\rho_s)(\tilde x)|^2\dd\rho_s(\tilde x)\dd s\right)^{\frac12}\left(\int_0^1\int_{\R^d}|v_s(\tilde x)|^2\dd\rho_s(\tilde x)\dd s\right)^{\frac12}\\
& \le \left(\int_0^1\int_{\R^d}|D_\mu f(\rho_s)(\tilde x)|^2\dd\rho_s(\tilde x)\dd s\right)^{\frac12}W_2(\rho^1,\rho^2).
\end{align*}
Furthermore as the function $\rho\mapsto M^2_2(\rho)$ is displacement convex, it is well known that 
$$M^2_2(\rho_s)\le\max\left\{M^2_2(\rho^1),M^2_2(\rho^2)\right\}.$$

With these tools in hand, we can perform the following computation
\begin{align*}
&\left|H(x^2,p^2,\rho^2) - H(x^1,p^1,\rho^1) \right| = \left|\int_0^1 \frac{\dd}{\dd s} H\left((1-s)x^1 + sx^2, (1-s)p^1 + sp^2,\rho_s \right)\dd s \right|\\
&\le |x^1-x^2|\int_0^1\left|D_x H\left((1-s)x^1 + sx^2, (1-s)p^1 + sp^2,\rho_s \right)\right|\dd s\\
&+ |p^1-p^2|\int_0^1\left|D_p H\left((1-s)x^1 + sx^2, (1-s)p^1 + sp^2,\rho_s \right)\right|\dd s\\
&+W_2(\rho^1,\rho^2)\left(\int_0^1\int_{\R^d}\left|D_\mu H\left((1-s)x^1 + sx^2, (1-s)p^1 + sp^2,\rho_s,\cdot \right)\right|^2\dd\rho_s\dd s\right)^{\frac12}
\end{align*}
By \eqref{growth:DH} and by the second moment estimates we have that there exists a constant $C>0$ such that
$$
\left|D_x H\left((1-s)x^1 + sx^2, (1-s)p^1 + sp^2,\rho_s \right)\right|\le C\left(1+\max\left\{M_2(\rho^1),M_2(\rho^2)\right\} + |x^1| + |x^2| + |p^1| + |p^2| \right)
$$
and 
$$
\left|D_x H\left((1-s)x^1 + sx^2, (1-s)p^1 + sp^2,\rho_s \right)\right|\le C\left(1+\max\left\{M_2(\rho^1),M_2(\rho^2)\right\} + |x^1| + |x^2| + |p^1| + |p^2| \right),
$$
for all $s\in[0,1]$. Furthermore, as $D^2_{x\mu}H, D^2_{p\mu}H$ are uniformly bounded and \eqref{hyp:add_H} is imposed, there exists a constant $C>0$ depending on these uniform bounds such that 
$$
\left|D_\mu H\left((1-s)x^1 + sx^2, (1-s)p^1 + sp^2,\rho_s,\tilde x \right)\right|\le C\left(1+ |\tilde x| + |x^1| + |x^2| + |p^1| + |p^2| \right),\ \ s\in[0,1].
$$

Therefore, all in all we can conclude that there exists a constant $C>0$ depending only on $H$ such that
\begin{align}\label{ineq:H_locLip}
&\left|H(x^2,p^2,\rho^2) - H(x^1,p^1,\rho^1) \right|\\ 
& \le C\left(1+M_2(\rho^1)+ M_2(\rho^2)  + |x^1| + |x^2| + |p^1| + |p^2|\right)\left(\left|x^1-x^2\right| + \left|p^1-p^2\right| + W_2(\rho^1,\rho^2)\right).
\end{align}

Then
\begin{align*}
&\left|L(x,D_pH(x^1,p^1,\rho^1),\rho^1) - L(x^2,D_pH(x^2,p^2,\rho^2),\rho^2)\right|\\[3pt]
&\le  \left|p^1\cdot D_pH(x^1,p^1,\rho^1) -  p^2\cdot D_pH(x^2,p^2,\rho^2)\right| + \left| H(x^1,p^1,\rho^1) - H(x^2,p^2,\rho^2)\right|\\[3pt]
&\le  \left|p^1\cdot D_pH(x^1,p^1,\rho^1) -  p^1\cdot D_pH(x^2,p^2,\rho^2)\right| +  \left|p^1\cdot D_pH(x^2,p^2,\rho^2) -  p^2\cdot D_pH(x^2,p^2,\rho^2)\right|\\[3pt] 
&+ \left| H(x^1,p^1,\rho^1) - H(x^2,p^2,\rho^2)\right|\\[3pt]
&\le C\left(|p^1|+ |D_pH(x^2,p^2,\rho^2)| \right)\left(\left|x^1-x^2\right| + \left|p^1-p^2\right| + W_2(\rho^1,\rho^2)\right)\\[3pt]
&+ \left| H(x^1,p^1,\rho^1) - H(x^2,p^2,\rho^2)\right|\\[3pt]
&\le C \left(1+|p^1|+ M_1(\rho^2) + |x^2| + |p^2| \right)\left(\left|x^1-x^2\right| + \left|p^1-p^2\right| + W_2(\rho^1,\rho^2)\right)\\[3pt]
&+ \left| H(x^1,p^1,\rho^1) - H(x^2,p^2,\rho^2)\right|\\[3pt]
&\le C \left(1+|x^1| + |p^1| + M_2(\rho^1) +M_2(\rho^2) + |x^2| + |p^2| \right)\left(\left|x^1-x^2\right| + \left|p^1-p^2\right| + W_2(\rho^1,\rho^2)\right) ,
\end{align*}
where we have used that \eqref{hyp:H_reg} implies \eqref{hyp:DH_Lip} and \eqref{growth:DH}, and in particular $C>0$ depends only on Lipschitz constants.
\end{proof}

Our main convergence result in this section reads as follows.

\begin{theorem}\label{thm:valueconvergence}  Assume that  \eqref{hyp:H_reg}--\eqref{hyp:add_g} are in force. 
Let $\rho_0\in\sP_2(\R^d)$ and $g:\R^d\times\sP_2(\R^d)$ be given, $(u^T,\rho^T)$ be as above. Then, there exists $\lambda \in \R$ such that the family of functions
\[
\left\{u^T (t,x) - \lambda (T-t)\right\}_{T \ge 0},
\]
converges locally uniformly on $[0,+\infty) \times \R^d$ as $T \to \infty$ to a function $u:[0,+\infty)\times\R^d\to\R$ which is $C_{\rm{loc}}^{1,1}$ in space and Lipschitz continuous in time, and $\rho^T$ converges to a flow of probability measures $\rho:[0,+\infty)\to\sP_2(\R^d)$ in $C([0,t]; (\sP_2(\R^d),W_2))$ for every $t > 0$. Moreover, the couple $(u,\rho)$ solves
\begin{align*}
\left\{
\begin{array}{ll}
-\partial_t  u  -\b\Delta u + H(x,-D  u,\rho) + \l = 0, & {\rm{in}\ } (0,+\infty)\times \R^d,\\
\partial_t\rho -\b\Delta \rho + \nabla\cdot (\rho D_pH(x,-Du,\rho)) = 0, & {\rm{in}\ } (0,+\infty)\times \R^d,\\
\rho(0,\cdot) = \rho_0,\ \
\ds\sup_{t\in[0,+\infty)}\frac{| u(t,x)|}{1+|x|^2} < \infty,
\end{array}
\right.
\end{align*}
where the fist equation is satisfied in the viscosity sense, while the second equation is satisfied in the sense of distributions.
\end{theorem}

\medskip

From now onwards, we suppose without loss of generality that $T \ge 2$. Let $\rho_0\in\sP_2(\R^d)$, $g:\R^d\times\sP_2(\R^d) \to \R$ be given. We introduce the following quantity
\begin{equation}\label{def:lambda_T}
\lambda^T :=\mathbb{E}\left\{ \int_{T/2}^{T/2+1} L(X^{0,\xi,T}_s,D_pH(X^{0,\xi,T}_s,Y^{0,\xi,T}_s,\rho^T_s),\rho^T_s)\dd s\right\}
\end{equation}
where $(X^{0,\xi,T}_s,Y^{0,\xi,T}_s, Z^{0,\xi,T}_s)_{s\in(0,T)}$ is the solution to \eqref{eq:FBSDE} with $\xi$ such that $\cL(\xi) = \rho_0$. 

We define moreover $\tilde u^T:[0,T]\times\R^d\to\R$ as
\begin{equation}\label{def:tilde_u}
\tilde u^T(t,x) := u^T (t,x) - \lambda^T (T-t), \ \ \forall\ (t,x)\in (0,T)\times\R^d
\end{equation}
and, for any $t\in [0,T]$ we introduce $\Psi^{T}:[0,T]\to\R$ given as
$$
\Psi^T(t) := \E \left[\tilde u^T(t,X^{0,\xi,T}_t)\right].
$$

\begin{lemma}\label{lem:bounds_in_T_2n}
Let $\rho_0\in\sP_2(\R^d)$, $g:\R^d\times\sP_2(\R^d) \to \R$ be given and let $(\tilde u^T,\rho^T,\lambda^T)$ be defined as above. Then there exists a constant $C>0$, depending on the data $\rho_0,H,g$, but independent of $T$, such that 
\begin{enumerate}
\item $\ds\sup_{t\in[0,T]}M_2(\rho^T_s)<C$,
\item $\ds\sup_{t\in[0,T]}|\tilde u^T(t,x)|\le C\left(1+|x|^2\right)$, 
\item for any $0 \le t_1 \le t_2 \le T$,
\[
\left|\lambda^T(t_1 - t_2) +  \mathbb{E} \left[\int_{t_1}^{t_2} L(X^{0,\xi,T}_s,D_pH(X^{0,\xi,T}_s,Y^{0,\xi,T}_s,\rho^T_s), \rho^T_s)\dd s\right]\right| \le C\left[e^{-c_0 T /4}+e^{-c_0 t_1} +e^{-c_0 (T-t_2)}\right].
\]
\end{enumerate}
\end{lemma}

\begin{proof}

First, by Proposition \ref{prop:2nd_moment_bound}, we find that $\sup_{s\in(0,T)}\E\left[|X^{0,\xi,T}_s|^2\right]$ is uniformly bounded, from where we conclude point (1) immediately. 

\medskip

Recall that 
$$
\E \left[u^T(t,X^{0,\xi,T}_t)\right] = \mathbb{E}\left\{\int_t^T L(X^{0,\xi,T}_s,D_pH(X^{0,\xi,T}_s,Y^{0,\xi,T}_s,\rho^T_s), \rho^T_s)\dd s + g(X^{0,\xi,T}_T,\rho^T_T)\right\}.
$$

Using the definition \eqref{def:tilde_u}, we have
\begin{align*}
\E \left[\tilde u^T(t,X^{0,\xi,T}_t)\right] + \lambda^T(T-t) &= \mathbb{E}\left\{\int_t^T L(X^{0,\xi,T}_s,D_pH(X^{0,\xi,T}_s,Y^{0,\xi,T}_s,\rho^T_s), \rho^T_s)\dd s + g(X^{0,\xi,T}_T,\rho^T_T)\right\},
\end{align*}
hence for any $0 \le t_1 \le t_2 \le T$,
\[
\Psi^T(t_1) - \Psi^T(t_2) =\lambda^T(t_1 - t_2) +  \mathbb{E} \int_{t_1}^{t_2} L(X^{0,\xi,T}_s,D_pH(X^{0,\xi,T}_s,Y^{0,\xi,T}_s,\rho^T_s), \rho^T_s)\dd s
\]
and so, for any $t \in [0,T]$
\begin{equation}\label{psi_equation}
\frac{\dd }{\dd t} \Psi^T(t) =  \lambda^T -  \mathbb{E} \big[ L(X^{0,\xi,T}_t,D_pH(X^{0,\xi,T}_t,Y^{0,\xi,T}_t,\rho^T_t), \rho^T_t) \big]
\end{equation} 

\medskip

\noindent {\bf Claim 1.} {\it There exists $C > 0$ such that for any $0 \le t \le T$ we have} 
$$\left|\frac{\dd }{\dd t} \Psi^T(t)\right| \le C\left[e^{-c_0 T /2}+e^{-c_0 t} +e^{-c_0 (T-t)}\right].$$ 

\medskip

\noindent {\bf Proof of Claim 1.} Set, for $0 \le t_1 \le t_2 \le T$ and $s\in[0,T-t_2+t_1]$,
$$(\bar X^{0,\xi,T}_s,\bar Y^{0,\xi,T}_s, \bar Z^{0,\xi,T}_s):= (X^{0,\xi,T}_{s+t_2-t_1},Y^{0,\xi,T}_{s+t_2-t_1}, Z^{0,\xi,T}_{s+t_2-t_1})\ \ {\rm{and}}\ \ \bar\rho^T_s:=\rho^T_{s+t_2-t_1}.$$
Then, we have
\begin{align*}
&\left|\frac{\dd }{\dd t} \Psi^T(t_1) - \frac{\dd }{\dd t} \Psi^T(t_2)\right| \\
&\le \E \left[\left| L(\bar X^{0,\xi,T}_{t_1},D_pH(\bar X^{0,\xi,T}_{t_1}, \bar Y^{0,\xi,T}_{t_1}, \bar \rho^T_{t_1}), \bar \rho^T_{t_1}) - L(X^{0,\xi,T}_{t_1},D_pH(X^{0,\xi,T}_{t_1},Y^{0,\xi,T}_{t_1},\rho^T_{t_1}), \rho^T_{t_1}) \right|\right] .
\end{align*}
Now, arguing as in \eqref{L1L2bound} we find
\begin{equation}\label{ineq:tilde_u2}
\left|\frac{\dd }{\dd t} \Psi^T(t_1) - \frac{\dd }{\dd t} \Psi^T(t_2)\right| \le C \sqrt{\Phi(t_1)} \le C\left[e^{-c_0(T-t_2)}+e^{-c_0 t_1}\right],
\end{equation}
where in the last inequality we have used Theorem \ref{prop:pointwisedecay} (1) (here, $C>0$ depends on $M_2(\rho_0)$, and in what follows it may increase from line to line). Indeed, it is crucial to mention that the results from Theorem \ref{prop:pointwisedecay} have been used for the two flows $( X^{0,\xi,T}_s, Y^{0,\xi,T}_s, Z^{0,\xi,T}_s)$ and $(\bar X^{0,\xi,T}_s,\bar Y^{0,\xi,T}_s, \bar Z^{0,\xi,T}_s)$ only on the time interval $[0,T-t_2+t_1]$. In particular, the inequality on $\Phi$ reads as 
$$
\Phi(s)\le C\left[e^{-2 c_0(T-t_2+t_1-s)}+e^{-2c_0s}\right],\ \ \forall s\in (0,T-t_2+t_1).
$$
Restricting now $t_1 \le T/2$ and integrating \eqref{ineq:tilde_u2} on $(T/2, T/2+1)$ with respect to $\dd t_2$ we deduce
$$
-C\left[e^{-c_0 T /2}+e^{-c_0 t_1}\right] \le  \frac{\dd }{\dd t} \Psi^T(t_1) - \Psi^T(T/2+1) + \Psi^T(T/2) \le C\left[e^{-c_0 T /2}+e^{-c_0 t_1}\right]
$$
but 
$$
\Psi^T(T/2) - \Psi^T(T/2+1)  = - \lambda^T +  \mathbb{E} \left[\int_{T/2}^{T/2+1} L(X^{0,\xi,T}_s,D_pH(X^{0,\xi,T}_s,Y^{0,\xi,T}_s,\rho^T_s), \rho^T_s)\dd s\right] = 0
$$
by \eqref{def:lambda_T}, so we get the claim for $t=t_1 \le T/2$. For $t \ge T/2$ we argue similarly, and integrate \eqref{ineq:tilde_u2} on $(T/2-1, T/2)$ with respect to $\dd t_1$, for $t_2\ge T/2$ fixed. The conclusion follows in this case by substituting $t_2=t$. 

We notice that in this case we need to also use that  
$$- \lambda^T +  \mathbb{E}\left[ \int_{T/2-1}^{T/2} L(X^{0,\xi,T}_s,D_pH(X^{0,\xi,T}_s,Y^{0,\xi,T}_s,\rho^T_s), \rho^T_s)\dd s\right]$$ 
is comparable to $e^{-c_0 T /2}$.  Indeed, we have
\begin{align*}
& \left|  -\lambda^T +  \mathbb{E} \int_{T/2-1}^{T/2} L(X^{0,\xi,T}_s,D_pH(X^{0,\xi,T}_s,Y^{0,\xi,T}_s,\rho^T_s), \rho^T_s)\dd s \right|\\
&=\Bigg| \mathbb{E} \left[\int_{T/2}^{T/2+1} L(X^{0,\xi,T}_s,D_pH(X^{0,\xi,T}_s,Y^{0,\xi,T}_s,\rho^T_s), \rho^T_s)\dd s\right]\\ 
&- \mathbb{E} \left[\int_{T/2-1}^{T/2} L(X^{0,\xi,T}_s,D_pH(X^{0,\xi,T}_s,Y^{0,\xi,T}_s,\rho^T_s), \rho^T_s)\dd s\right]  \Bigg|\\
&=\left| \mathbb{E} \left[\int_{T/2}^{T/2+1}\left[ L(X^{0,\xi,T}_s,D_pH(X^{0,\xi,T}_s,Y^{0,\xi,T}_s,\rho^T_s), \rho^T_s) -  L(\bar X^{0,\xi,T}_s,D_pH(\bar X^{0,\xi,T}_s,\bar Y^{0,\xi,T}_s,\bar \rho^T_s), \bar\rho^T_s)\right]\dd s\right]  \right|,\\
\end{align*}
where $\left(\bar X^{0,\xi,T}_s,\bar Y^{0,\xi,T}_s,\bar \rho^T_s\right):=\left(X^{0,\xi,T}_{s+1},Y^{0,\xi,T}_{s+1},\rho^T_{s+1}\right).$ Now, again, by \eqref{L1L2bound} and Theorem \ref{prop:pointwisedecay} (1) we find
\begin{align*}
 \Bigg|  -\lambda^T +  \mathbb{E}\bigg[ \int_{T/2-1}^{T/2} L(X^{0,\xi,T}_s,&D_pH(X^{0,\xi,T}_s,Y^{0,\xi,T}_s,\rho^T_s), \rho^T_s)\dd s\bigg] \Bigg|\\
&\le C\int_{T/2}^{T/2+1}\sqrt{\Phi(s)}\dd s \le C\int_{T/2}^{T/2+1}\left[e^{-{c_0} (T-s)} + e^{-{c_0} s}\right]\dd s\le Ce^{-c_{0}T/2},
\end{align*}
as desired.

\medskip

\noindent {\bf Claim 2.} {\it $\Psi^T(t)$ is uniformly bounded in $t, T$} . 

 For this, it is sufficient to observe that
 $$
 \Psi^T(t) = \Psi^T(T) - \int_{t}^T \frac{\dd }{\dd s} \Psi^T (s) \dd s,
 $$
 and one concludes by the estimate of Claim 1 and the fact that $\Psi^T(T)$ is uniformly bounded with respect to $T$ using the assumption \eqref{hyp:add_g} on $g$.
 
 \medskip

\noindent {\bf Claim 3.} {\it There exists a constant $C>0$ such that $|\tilde u^T(t,x)| \le C\left(1+|x|^2\right).$}

We have
\begin{align*}
\left|\tilde u^T(t,x)\right| &\le \E\left[ \left|\tilde u^T(t,x) - \tilde u^T(t,X^{0,\xi,T}_t)\right|\right] + \E\left[\left|\tilde u^T(t,X^{0,\xi,T}_t)\right|\right]\\
&\le \E \left[ \left|D_x\tilde u^T(t,y)\right| \left|x-X^{0,\xi,T}_t\right| \right] + \left|\Psi^T(t)\right| \\
& \le C \E \left[\left(1+|x|+\left|X^{0,\xi,T}_t\right|\right)\left|x-X^{0,\xi,T}_t\right|\right] + C \le C\left(1+|x|^2\right),
\end{align*}
where in the penultimate line $y$ is a vector on the line segment connecting $x$ to $X^{0,\xi,T}_t$ and in the last line we have used that $D_x\tilde u^T(t,\cdot)=D_xu^T(t,\cdot)$ grows at most linearly at infinity (cf. Corollary \ref{cor:D_xu_lin}). This concludes the proof of {\bf Claim 3} and shows point (2) in the statement of this lemma.

\medskip

To get the last point (3) of the lemma, recall that \eqref{psi_equation} reads
\begin{multline*}
\lambda^T(t_1 - t_2) +  \mathbb{E}\left[ \int_{t_1}^{t_2} L(X^{0,\xi,T}_s,D_pH(X^{0,\xi,T}_s,Y^{0,\xi,T}_s,\rho^T_s), \rho^T_s)\dd s\right] \\ = \Psi^T(t_1) - \Psi^T(t_2) = - \int_{t_1}^{t_2} \frac{\dd }{\dd s} \Psi^T (s) \dd s,
\end{multline*}
and again by the estimates on the derivative of $\Psi^T$ in {\bf Claim 1} one concludes.


\end{proof}

\begin{proposition}\label{prop:lambda_unique}
We suppose that we are in the setting of Lemma \ref{lem:bounds_in_T_2n}. Then the limit 
\[
\lim_{T\to \infty} \lambda^T = \lambda
\]
exists and it is finite, and it is independent of $\rho_0$ and $g$. Moreover, $|\lambda^T-\lambda| \le C e^{-c_0 T/2}$,  where $C$ depends on the data and the second moment of $\rho_0$.
\end{proposition}

\begin{proof}
Let us consider two MFG Nash equilibria $(\rho^T_s)_{s\in(0,T)}$ and $(\rho^{\hat T}_s)_{s\in(0,\hat T)}$ with data $(\rho^1_0,g^1)$ and $(\rho^2_0,g^2)$, set on time horizons of $T$ and $\hat T$, respectively. We set the Hamiltonian to be the same for both of them, and without loss of generality, we assume that $\hat T > T$. Recall the definition of $\lambda^T$ and $\lambda^{\hat T}$ from \eqref{def:lambda_T}. Let $(X^{0,\xi,T}_s,Y^{0,\xi,T}_s, Z^{0,\xi,T}_s)_{s\in(0,T)}$ and $(\hat X^{0,\hat \xi,\hat T}_s,\hat Y^{0,\hat\xi,\hat T}_s, \hat Z^{0,\hat\xi,\hat T}_s)_{s\in(0,\hat T)}$ be the corresponding solutions to \eqref{eq:FBSDE} with $\cL(\xi) = \rho^1_0$ and $\cL({\hat\xi}) = \rho^2_0$. Set
$$(\bar X^{0,\hat\xi,\hat T}_s,\bar Y^{0,\hat\xi,\hat T}_s, \bar Z^{0,\hat\xi,\hat T}_s):= (\hat X^{0,\hat\xi,\hat T}_{s+\hat T / 2- T/2},\hat Y^{0,\hat\xi,\hat T}_{s+\hat T / 2- T/2}, \hat Z^{0,\hat\xi,\hat T}_{s+\hat T / 2- T/2})\ \ {\rm{and}}\ \ \bar\rho^{\hat T}_s:=\rho^{\hat T}_{s+\hat T / 2- T/2},$$
for $s\in[0,T].$ We notice that if $s\in[0,T]$, then $s+\hat T / 2- T/2\in [\hat T/2 - T/2,T/2+\hat T/2]\subset [0,\hat T]$ , and so these new curves are well-defined.

Now, arguing as in \eqref{L1L2bound}, we can deduce
\begin{multline}\label{ineq:same_lambda}
\left|\lambda^T - \lambda^{\hat T}\right| = \\
\le \int_{T/2}^{T/2+1}\E \left[\left| L(X^{0,\xi,T}_s,D_pH(X^{0,\xi,T}_s,Y^{0,\xi,T}_s,\rho^T_s), \rho^T_s) - L(\bar X^{0,\hat\xi,\hat T}_s,D_pH(\bar X^{0,\hat\xi,\hat T}_s,\bar Y^{0,\hat\xi,\hat T}_s,\bar \rho^{\hat T}_s), \bar \rho^{\hat T}_s) \right|\right]\dd s\\
\le Ce^{-c_0 T/2},
\end{multline}
where in the last inequality we have used that since the flows are taken on the interval $[0,T]$, hence for $s\in[T/2,T/2+1]$ we have by Theorem \ref{prop:pointwisedecay} (1) 
$$
\Phi(s)\le C\left[e^{-2c_0(T-s)}+e^{-2 c_0 s}\right] = Ce^{-c_0 T}.
$$
Now, \eqref{ineq:same_lambda} shows that $\lambda^T$ is a Cauchy sequence, and yields the desired result. The above arguments show also that $\lambda$ does not depend on $\rho_0, g$.
\end{proof}

\begin{proposition}\label{prop:tildeulimit}
We suppose that we are in the setting of Lemma \ref{lem:bounds_in_T_2n}. Let $T^1 \le T^2$ and $t \le T^1/2$. Let $(X^{0,\xi,T^1}_s,Y^{0,\xi,T^1}_s, Z^{0,\xi,T^1}_s)_{s\in(0,T)}$ and $(X^{0,\xi,T^2}_s,Y^{0,\xi,T^2}_s, Z^{0,\xi,T^2}_s)_{s\in(0,T)}$ denote the corresponding solutions to \eqref{eq:FBSDE} with $\xi\in L^2(\Om,\cF_{0},\P;\R^{d})$ such that $\cL(\xi) = \rho_0$, with the same final datum $g$, on time horizons $(0,T^1)$ and $(0,T^2)$, respectively. Then, 
\[
\left|\E \left[\tilde u^{T^1}(t,X^{0,\xi,{T^1}}_{t})\right]  - \E \left[\tilde u^{T^2}(t,X^{0,\xi,{T^2}}_{t})\right]\right| \le C \left[e^{-c_0 t} + e^{-c_0 (T^1-t)} + e^{-c_0 T^1 /4}\right] 
\]
for some $C>0$ depending on the data $\rho_0,H,g$.
\end{proposition}

\begin{proof} Let us consider two MFG Nash equilibria $(\rho^{T^1}_s)_{s\in(0,{T^1})}$ and $(\rho^{T^2}_s)_{s\in(0,{T^2})}$  with data $(\rho_0,g)$ set on time horizons of $T^1$ and $T^2$. Recall that
\begin{multline}\label{tildeuT}
\E \left[\tilde u^{T^1}(t,X^{0,\xi,{T^1}}_{t})\right] + t \lambda^{T^1}  = \E \left\{\int_{t}^{T^1-t} L(X^{0,\xi,T^1}_s,D_pH(X^{0,\xi,T^1}_s,Y^{0,\xi,T^1}_s,\rho^{T^1}_s), \rho^{T^1}_s)\dd s\right\}\\
 - \left(T^1- 2t\right)\lambda^{T^1} \\
+ \E \left\{\int_{T^1-t}^{T^1} L(X^{0,\xi,T^1}_s,D_pH(X^{0,\xi,T^1}_s,Y^{0,\xi,T^1}_s,\rho^{T^1}_s), \rho^{T^1}_s)\dd s + g(X^{0,\xi,T^1}_{T^1}, \rho^{T^1}_{T^1})\right\}.
\end{multline}
and
\begin{multline}\label{tildeuhatT}
\E\left[ \tilde u^{T^2}(t,X^{0,\xi,{T^2}}_{t})\right] +t \lambda^{T^2}  = \E \left\{\int_{t}^{T^2-t} L(X^{0,\xi,T^2}_s,D_pH(X^{0,\xi,T^2}_s,Y^{0,\xi,T^2}_s,\rho^{T^2}_s), \rho^{T^2}_s)\dd s\right\}\\
 - \left(T^2- 2t\right)\lambda^{T^2} \\
+ \E \left\{\int_{T^2-t}^{T^2} L(X^{0,\xi,T^2}_s,D_pH(X^{0,\xi,T^2}_s,Y^{0,\xi,T^2}_s,\rho^{T^2}_s), \rho^{T^2}_s)\dd s + g(X^{0,\xi,T^2}_{T^2}, \rho^{T^2}_{T^2})\right\}.
\end{multline}
Notice first that using Lemma \ref{lem:bounds_in_T_2n} (3),
\[
\left| \E \left[\int_{t}^{T^1-t} L(X^{0,\xi,T^1}_s,D_pH(X^{0,\xi,T^1}_s,Y^{0,\xi,T^1}_s,\rho^{T^1}_s), \rho^{T^1}_s)\dd s\right] - \left(T^1- 2t\right)\lambda^{T^1} \right | \le  C\left[e^{-c_0 T^1 /4}+e^{-c_0 t}\right].
\]
and
\[
\left| \E \left[\int_{t}^{T^2-t} L(X^{0,\xi,T^2}_s,D_pH(X^{0,\xi,T^2}_s,Y^{0,\xi,T^2}_s,\rho^{T^2}_s), \rho^{T^2}_s)\dd s\right] - \left(T^2- 2t\right)\lambda^{T^2} \right | \le C\left[e^{-c_0 T^2 /4}+e^{-c_0 t}\right].
\]

Furthermore, by time shift $ s(\bar s) = s -T^2 + T^1$ we can compare
\begin{multline*}
\left| \E \int_{T^1-t}^{T^1} L(X^{T^1}_s,D_pH(X^{T^1}_s,Y^{T^1}_s,\rho^{T^1}_s), \rho^{T^1}_s)\dd s -  \E \int_{T^2-t}^{T^2} L(X^{T^2}_s,D_pH(X^{T^2}_s,Y^{0,\xi,T^2}_s,\rho^{T^2}_s), \rho^{T^2}_s)\dd s \right| = \\
\left| \E \int_{T^1-t}^{T^1} L(X^{T^1}_s,D_pH(X^{T^1}_s,Y^{T^1}_s,\rho^{T^1}_s), \rho^{T^1}_s)\dd s -  \E \int_{T^1-t}^{T^1} L(X^{T^2}_{s(\bar s)},D_pH(X^{T^2}_{s(\bar s)},Y^{0,\xi,T^2}_{s(\bar s)},\rho^{T^2}_{s(\bar s)}), \rho^{T^2}_{s(\bar s)})\dd \bar s \right| \\ \le Ce^{-c_0(T^1-t)} ,
\end{multline*}
arguing as in \eqref{L1L2bound} and applying Theorem \ref{prop:pointwisedecay} (3) (note that we are using the crucial fact that, after time-shift, the two MFG Nash equilibria enjoy the same final condition at time $T^{1}$). Similarly, we will have
\[
\left|\E \left[g(X^{0,\xi,T^{1}}_{T^{1}},\rho^{T^{1}}_{T^{1}})\right] - \E \left[g(X^{0,\xi,{T^{2}}}_{T^{2}},\rho^{T^{2}}_{T^{2}})\right]\right| \le Ce^{-c_0T^1}.
\]
Indeed, for $s\in (t,T^1)$ let us consider 
$$
\left(\bar X^{T^2}_{s},\bar Y^{0,\xi,T^2}_{s},\bar Z^{0,\xi,T^2}_{s},\bar\rho^{T^2}_{s}\right):=\left(X^{T^2}_{s-T_1+T_2},Y^{0,\xi,T^2}_{s-T_1+T_2},\bar Z^{0,\xi,T^2}_{s-T_1+T_2},\rho^{T^2}_{s-T_1+T_2}\right).
$$
With this choice, we have in particular that both flows, $\left( X^{T^1}_{s}, Y^{0,\xi,T^1}_{s}, Z^{0,\xi,T^1}_{s},\rho^{T^1}_{s}\right)_{s\in(t,T^1)}$ and $\left(\bar X^{T^2}_{s},\bar Y^{0,\xi,T^2}_{s},\bar Z^{0,\xi,T^2}_{s},\bar\rho^{T^2}_{s}\right)_{s\in(t,T^1)}$ are defined on the same time interval $(t,T^1)$ and both MFG Nash equilibria correspond to the same final condition $g$ (and possible different initial conditions). With this in mind we find that $\Phi$, when associated to these two flows will satisfy in particular (see again Theorem \ref{prop:pointwisedecay} (3))
$$
\Phi(T^1)\le Ce^{-2c_0 T^1}.
$$
Now, using \eqref{hyp:add_g} and the exact same arguments which led to \eqref{ineq:H_locLip}, we find that 
$$
\left|\E \left[g(X^{0,\xi,T^{1}}_{T^{1}},\rho^{T^{1}}_{T^{1}})\right] - \E \left[g(X^{0,\xi,{T^{2}}}_{T^{2}},\rho^{T^{2}}_{T^{2}})\right]\right| \le C \sqrt{\Phi(T^1)}\le Ce^{-c_0 T^1},
$$
where the constant $C>0$ depends on the data and on $M_2(\rho_0)$, and so our claim follows.

\medskip

Taking finally the difference between \eqref{tildeuT} and \eqref{tildeuhatT}, and plugging in all the previous inequalities we obtain
\[
\left|\E \left[\tilde u^{T^1}(t,X^{0,\xi,{T^1}}_{t})\right]- \E \left[\tilde u^{T^2}(t,X^{0,\xi,{T^2}}_{t})\right] \right| \le C \left[e^{-c_0 T^1 /4} + e^{-c_0 T^2 /4}+ e^{-c_0 t} + e^{-c_0(T^1-t)}\right],
\]
which shows the desired assertion.
\end{proof}

Now, let us discuss about the convergence of optimal trajectories.  We have the following result.

\begin{proposition}\label{prop:unique_X}
We suppose that we are in the setting of Lemma \ref{lem:bounds_in_T_2n}. For $T>0$ recall that $(X^{0,\xi,T}_s,Y^{0,\xi,T}_s, Z^{0,\xi,T}_s)_{s\in(0,T)}$ denotes the corresponding solutions to \eqref{eq:FBSDE} with $\xi\in L^2(\Om,\cF_{0},\P;\R^{d})$ such that $\cL(\xi) = \rho_0$. Then there exists a unique process $(X^\xi_s)_{s\in(0,+\infty)}$, independent of the final condition $g$, such that $X^\xi_0 = \xi$ and
$$
\lim_{T\to+\infty}\sup_{s\in[0,t]}\E \left[\left|X^{0,\xi,T}_s - X^\xi_s\right|^2\right] = 0,\ \ \forall t>0.
$$
Note that, since $\rho^{T}_s = \cL({X^{0,\xi,T}_s})$, we have that there exists a unique continuous curve $\rho:[0,+\infty)\to\sP_2(\R^d)$ starting at $\rho_0$, such that $\lim_{T \to+\infty}\sup_{s\in[0,t]} W_2(\rho^T_s,\rho_s) = 0$ for all $t>0$. More precisely, 
$$
W^2_2(\rho^T_s,\rho_s) \le C e^{-2c_0(T-t)}, \forall \ s\in [0,t].
$$
\end{proposition}

\begin{proof}
First, by Proposition \ref{prop:2nd_moment_bound}, we know that $\sup_{s\in(0,T)}\E\left[\left|X^{0,\xi,T}_s \right|^2\right]$ is uniformly bounded with respect to  $T>0$. Therefore, it is enough to prove that the family $(X^{0,\xi,T}_s)_{s\in(0,T)}$ is Cauchy with respect to $T$ in $C([0,t]; L^2(\Om,\cF,\P;\R^{d}))$, for any $t > 0$. For this, fix $t>0$ and let $T$ and $\hat T$ be given with the property $t<T<\hat T$. Let $\hat g:\R\times\sP_2(\R^d)\to\R$ be any final datum which satisfies our standing assumptions. Let us consider $(\hat X^{0,\xi,\hat T}_s,\hat Y^{0,\xi,\hat T}_s, \hat Z^{0,\xi,\hat T}_s)_{s\in(0,\hat T)}$ be the corresponding solutions to \eqref{eq:FBSDE} with $\xi\in L^2(\Om,\cF_{0},\P;\R^{d})$ such that $\cL(\xi) = \rho_0$ and $\hat g$ as a final datum.

Now, when restricting both triples to the time interval $(0,T)$, they will both describe MFG Nash equilibria, with the final data given by $g$ and $\hat u(T,\cdot)$, this being the value function associated to the second game, at time $T$.

Therefore, by Theorem \ref{prop:pointwisedecay}(2), using the definition of $\Phi$ from Definition \ref{def:phi_fcts}, we deduce that
\begin{equation}\label{rhoTconvergence}
\E\left[\left|X^{0,\xi,T}_s - \hat X^{0,\xi,\hat T}_s\right|^2\right] \le C e^{-2c_0(T-s)} \le C e^{-2c_0(T-t)}, \forall \ s\in (0,t).
\end{equation}

Then, the result follows.
\end{proof}

\begin{remark}
It is important to remark that the limit process $(X^\xi_s)_{s\in(0,+\infty)}$ that we obtain in Proposition \ref{prop:unique_X} could in general have very low regularity. Indeed, in particular a priori we do not even know if $(X^\xi_s(\omega))_{s\in(0,+\infty)}$ is a continuous path, for $\omega\in\Om$. However, by construction we must have that $X^\xi_s$ is $\cF_{s}$-measurable for all $s\ge 0$.

It remains an interesting open question to investigate whether the process $(X^\xi_s(\omega))_{s\in(0,+\infty)}$ could be related to an infinite horizon FBSDE system, such as the ones appearing in \cite{BayZha} for instance.
\end{remark}

\begin{proposition}\label{prop:tildeulimit0}
We suppose that we are in the setting of Lemma \ref{lem:bounds_in_T_2n}. Then, for any $\t \le T^1/4$,
\[
\left| \tilde u^{T^1}(\tau , 0)  - \tilde u^{T^2}(\tau , 0) \right| \le C\left[e^{-c_0 T^1 /4} +e^{-c_0(T^1/4 - \tau)} \right],
\]
for some $C>0$ depending on the data $\rho_0,H,g$.
\end{proposition}

\begin{proof} Let us denote by $(X^{i,0,\xi}_s,Y^{i,0,\xi}_s,Z^{i,0,\xi}_s)_{s\in[0,T^i]}$, $i=1,2$ the solutions to the FBSDE system \eqref{eq:FBSDE} associated with the Mean Field equilibria originating from $\rho_0$, on time horizons $T^1$ and $T^2$ respectively, and by $(\hat X^{i,\tau,\hat \xi}_s, \hat Y^{i,\tau,\hat \xi}_s, \hat Z^{i,\tau,\hat \xi}_s)_{s\in[\tau,T^i]}$, $i=1,2$ the solutions to the FBSDE system \eqref{eq:FBSDExi} with input $\rho_s^i = \mathcal L (X^{i,0,\xi^i}_s)$ and $\hat \xi = \delta_0$.

\medskip
{\bf Step 1.} Recall first that by \eqref{rhoTconvergence} we have that $W_2 (\rho_s^{T^1}, \rho_s^{T^2}) \le Ce^{-c_0 (T^1-s)}$ for any $s \in [0, T^1]$. Hence, denoting by
\[
\hat \Phi (s) :=  \mathbb E\left[\left| \hat X^{1,\tau,\hat \xi}_s - \hat  X^{2,\tau,\hat \xi}_s\right|^2\right] +\mathbb E\left[ \left| \hat Y^{1,\tau,\hat \xi}_s - \hat  Y^{2,\tau,\hat \xi}_s \right|^2\right],
\]
and the associated function $\hat\varphi$ along these flows (as defined in Definition \ref{def:phi_fcts}). First, we observe that by \eqref{phivarphiest} and \eqref{phivarphi} we have that $\hat\varphi(T^1)$ is uniformly bounded and $\hat\varphi(\tau) = 0$ (as we consider the same starting random variable $\hat\xi$ at time $\tau$).

Then, for any $t > \tau$, we get by Proposition \ref{prop:decay_Phi} (1) (by setting $(t,t_1,t_2,T)$ from that proposition as $(\tau,\tau,t,T^1)$)

\begin{equation}\label{hatphi_est}
C \int_{\t}^{t} \hat \Phi(s) \dd s \le e^{-c_0 (T^1-t) } +\left(  
  \int_{\t}^{t} e^{-c_0 (T^1-s) }\dd s + 
 \int_{t}^{T^1} e^{-2 c_0 (T^1-t)}  \dd s \right) \le C' e^{-c_0 (T^1-t) },
\end{equation}
for some $C'>0$ depending only on the data and $M_2(\rho_0).$
Recall that for $i=1,2$ it holds
\begin{align*}
\tilde u^{T^i}(\tau , 0) + \lambda^{T^i}(t-\tau) &= \mathbb{E}\left\{\int_\tau^t L(\hat X^{i, \tau,\hat \xi}_s,D_pH(\hat X^{i,\tau,\hat \xi},\hat Y^{i, \tau,\hat \xi},\rho^{T^i}_s), \rho^{T^i}_s)\dd s + \tilde u^{T^i}(t, \hat X^{i, \tau,\hat \xi}_t)\right\}.
\end{align*}
Therefore, taking differences and arguing as in \eqref{L1L2bound}, with the estimate in \eqref{hatphi_est} we get
\[
\left|\tilde u^{T^1}(\tau , 0) - \tilde u^{T^2}(\tau , 0)\right| \le \left|\lambda^{T^1}-\lambda^{T^2}\right|(t-\tau) + C  e^{-c_0 (T^1-t) /2} +
\left|\E \left[\tilde u^{T^1}(t, \hat X^{1, \tau,\hat \xi}_t)\right] - \E \left[\tilde u^{T^2}(t, \hat X^{2, \tau,\hat \xi}_t)\right]\right|.
\]
By the rate of convergence of $\lambda^T$ provided in Proposition \ref{prop:lambda_unique} we the conclude
\[
\left|\tilde u^{T^1}(\tau , 0) - \tilde u^{T^2}(\tau , 0)\right| \le Ce^{-c_0 T^1 /4} + C  e^{-c_0 (T^1-t) /2} +
\left|\E \left[\tilde u^{T^1}(t, \hat X^{1, \tau,\hat \xi}_t)\right] - \E \left[\tilde u^{T^2}(t, \hat X^{2, \tau,\hat \xi}_t)\right]\right|.
\]

\medskip

{\bf Step 2.} We now proceed by estimating the last term of the previous inequality. By the triangle inequality, that is controlled by 
\begin{multline*}
\left|\E \left[\tilde u^{T^1}(t, \hat X^{1, \tau,\hat \xi}_t)\right] - \E \left[\tilde u^{T^1}(t, X^{0,\xi,{T^1}}_{t})\right]\right| + \left|\E \left[\tilde u^{T^1}(t, X^{0,\xi,{T^1}}_{t})\right] - \E \left[\tilde u^{T^2}(t, X^{0,\xi,{T^2}}_{t})\right]\right| \\ + \left|\E \left[\tilde u^{T^2}(t, X^{0,\xi,{T^2}}_{t})\right] - \E \left[\tilde u^{T^2}(t, \hat X^{2, \tau,\hat \xi}_t)\right]\right|.
\end{multline*}
We start by the second term, that can be controlled by Proposition \ref{prop:tildeulimit} by 
$$\left|\E\left[ \tilde u^{T^1}(t, X^{0,\xi,{T^1}}_{t})\right] - \E \left[\tilde u^{T^2}(t, X^{0,\xi,{T^2}}_{t})\right]\right| \le C\left(e^{-c_0 t} + e^{-c_0 (T^1-t) } + e^{-c_0 T/4}\right).$$
To estimate the first one, define
\[
\Phi (s) =  \mathbb E\left[\left| \hat X^{1,\tau,\hat \xi}_s - X^{0,\xi,{T^1}}_s\right|^2\right] +\mathbb E \left[\left| \hat Y^{1,\tau,\hat \xi}_s - Y^{0,\xi,{T^1}}_s\right|^2\right],
\]
and apply Remark \ref{rmk:decay} to conclude that for any $s \in [\tau, T^1]$. 
\[
\Phi (s) \le C e^{-2c_0(s-\tau)}
\]
Therefore, employing the gradient bounds of Corollary \ref{cor:D_xu_lin} and uniform second moment bounds we get
\begin{multline*}
\left|\E \left[\tilde u^{T^1}(t, \hat X^{1, \tau,\hat \xi}_t)\right] - \E \left[\tilde u^{T^1}(t, X^{0,\xi,{T^1}}_{t})\right]\right| \le C \E \left[\left(1+ \left|\hat X^{1, \tau,\hat \xi}_t\right| + \left|X^{0,\xi,{T^1}}_{t}\right|\right) \left|\hat X^{1, \tau,\hat \xi}_t - X^{0,\xi,{T^1}}_{t}\right| \right] \\
\le C \sqrt{\Phi(t)} \le Ce^{-c_0(t-\tau)}.
\end{multline*}
Since the term $\left|\E \left[\tilde u^{T^2}(t, X^{0,\xi,{T^2}}_{t})\right] - \E \left[\tilde u^{T^2}(t, \hat X^{2, \tau,\hat \xi}_t)\right]\right|$ can be handled analogously, we conclude that
\[
\left|\tilde u^{T^1}(\tau , 0) - \tilde u^{T^2}(\tau , 0)\right| \le C\left[e^{-c_0 T^1 /4} +  e^{-c_0 (T^1-t) /2} + e^{-c_0 t}  +e^{-c_0(t-\tau)}\right].
\]
Choosing finally $t=T^1/4$ yields the desired assertion.

\end{proof}

\begin{corollary}\label{prop:tildeulimit1} 
We suppose that we are in the setting of Lemma \ref{lem:bounds_in_T_2n}. Let $\lambda$ be as in Proposition \ref{prop:lambda_unique}. Then, the family of functions
\[
\left\{u^T (t,x) - \lambda (T-t)\right\}_{T \ge 0}
\]
converges locally uniformly on $[0,+\infty) \times \R^d$ as $T \to \infty$ to a function $u$ which is $C^{1,1}_{\rm{loc}}$ is space and Lipschitz in time. In particular, for every $t \le T/8$, $x \in \R^d$,
\begin{equation}\label{convuDxu}
\begin{gathered}
\left| u^T (t,x) - \lambda (T-t) - u(t,x) \right| \le C e^{-c_0 T^1 /8} \left(1+|x|^2\right), \\
\left| D_x u^{T} (t,x) -D_x u (t,x) \big)\right| \le C e^{-c_0 T^1 /4} \left(1+|x|\right)
\end{gathered}
\end{equation}
for some $C>0$ depending on the data $\rho_0,H,g$. Moreover, $u$ is a viscosity solution to
\begin{align*}
\left\{
\begin{array}{ll}
-\partial_t  u  -\b\Delta u + H(x,-D  u,\rho) + \l = 0, & \text{in } (0,+\infty)\times \R^d,\\
\ds\sup_{t\in[0,+\infty)}\frac{| u(t,x)|}{1+|x|^2} < \infty
\end{array}
\right.
\end{align*}
where $\rho$ is as in Proposition \ref{prop:unique_X}.
\end{corollary}

\begin{proof}  Let $T^2 > T^1 > 0$ and note first that for all $x \in \R^d$ and $t \le T_1$,
\begin{align*}
&\left| u^{T^1} (t,x) - \lambda (T^1-t) - \big(u^{T^2} (t,x) - \lambda (T^2-t)  \big)\right| = \\
&\left| \tilde u^{T^1} (t,x) + (\lambda^{T^1} - \lambda) (T^1-t) - \tilde u^{T^2} (t,x) - (\lambda^{T^2} - \lambda) (T^2-t)  \big)\right| \le \left| \tilde u^{T^1} (t,x) - \tilde u^{T^2} (t,x)\right| + Ce^{-c_0T^1/4},
\end{align*}
by Proposition \ref{prop:lambda_unique}. Moreover, the Mean Value Theorem, Proposition \ref{prop:tildeulimit0} and Theorem \ref{prop:Du_loc} show that for $t \le T^1/4$,
\begin{align*}
 \left| \tilde u^{T^1} (t,x) - \tilde u^{T^2} (t,x)\right| &\le  \left| \tilde u^{T^1} (t,0) - \tilde u^{T^2} (t,0)\right| + |x| \sup_{|y| \le |x|} \left(\left|D_x \tilde u^{T^1} (t,y) - D_x \tilde u^{T^2} (t,y)\right| \right)\\
&\le C\left[e^{-c_0 T^1 /4} +e^{-c_0(T^1/4 - t)} + e^{-c_0 T^1 /4} |x|(1+|x|)\right].
\end{align*}
By combining the two previous inequalities we get that, for all $x$ and $t \le T^1/8$,
\[
\left| u^{T^1} (t,x) - \lambda (T^1-t) - \big(u^{T^2} (t,x) - \lambda (T^2-t)  \big)\right| \le C e^{-c_0 T^1 /8} (1+|x|^2).
\]
Recall also that Theorem \ref{prop:Du_loc} gives, for all $x$ and $t \le T^1/4$,
\[
| D_x u^{T^1} (t,x) -D_x u^{T^2} (t,x) \big)| \le C e^{-c_0 T^1 /4} (1+|x|).
\]
Therefore, the sequences
$$\left\{u^T  - \lambda (T-t)\right\}_{T\ge 0}, \qquad  \left\{D_x\big(u^T - \lambda (T-t)\big)\right\}_{T\ge 0}$$ 
are Cauchy in $C([0, \tau] \times \Omega)$ for every compact $\Omega \subset \R^d$ and $\tau > 0$ (and $T \ge 8\tau$), hence $u^T  - \lambda (T-t)$ (and its gradient $D_xu^T$) converges locally uniformly to $u$ (to $D_xu$, respectively) as $T \to \infty$.

Finally, since $\left(\rho^T_{s}\right)_{s\in[0,T]}$ converges locally uniformly (in time) to $(\rho_{s})_{s\in[0,+\infty)}$, with $\sup_{s\ge 0}M_{2}(\rho_{s})<+\infty$ and $H$ is assumed to be locally Lipschitz continuous in the measure variable, then $H(x,p,\rho^T_{t})$ converges locally uniformly in $(t,x,p)$) to $H(x,p,\rho_{t})$, and that $u$ is a viscosity solution to the PDE follows by standard stability arguments in the theory of viscosity solutions. The uniform control (in time) on the quadratic growth of $u$ is a consequence of Lemma \ref{lem:bounds_in_T_2n} (2).

\end{proof}

We are now ready to conclude the proof of Theorem \ref{thm:valueconvergence}, which is Theorem \ref{thm:intro2} of the introduction. 

\begin{proof}[Proof of Theorem \ref{thm:valueconvergence}] This almost follows from Proposition \ref{prop:unique_X}, on the convergence of $\rho^T$,  and Corollary \ref{prop:tildeulimit1}, on the convergence of $u^T$. We are just left to show that the limit $(u,\rho)$ satisfies the Fokker--Planck equation in the sense of distributions, but this is a consequence of the convergence of $\rho^T$ in $C_{\rm loc}([0,+\infty), \sP_2(\R^d))$, the estimates on the gradients in \eqref{convuDxu} and the regularity assumptions on $D_pH$. \end{proof}

We conclude this section with a uniqueness result for the limit system satisfied by $(u,\lambda,\rho)$.

\begin{theorem}\label{uniquelim}
For $i=1,2$, let $\l^i \in\R$,  $u^i:[0,+\infty)\times\R^d\to\R$ be $C_{\rm{loc}}^{1,1}$ in space and Lipschitz continuous in time and $\rho^i \in C([0,+\infty); (\sP_2(\R^d),W_2))$ be solutions of
\begin{align*}
\left\{
\begin{array}{ll}
-\partial_t  u  -\b\Delta u + H(x,-D_{x}  u,\rho) + \l = 0, & {\rm{in}\ } (0,+\infty)\times \R^d,\\
\partial_t\rho -\b\Delta \rho + \nabla\cdot (\rho D_pH(x,-D_{x}u,\rho)) = 0, & {\rm{in}\ } (0,+\infty)\times \R^d,\\
\rho(0,\cdot) = \rho_0,\ \
\ds\sup_{t\in[0,+\infty)}\frac{| u(t,x)|}{1+|x|^2} < \infty,
\end{array}
\right.
\end{align*}
Then, $\l^1 = \l^2$, $\rho^1 \equiv \rho^2$ and $u^1 \equiv u^2+c$ for some $c \in \R$.
\end{theorem} 

\begin{proof} Note first that, for any $T>0$, we can make use of  the two triples $(X^{i,T}_s,Y^{i,T}_s,Z^{i,T}_s)_{s\in[0,T]}$, $i=1,2$ that solve the FBSDE system \eqref{eq:FBSDE_intro} with $g(x) = u^1(x,T)$ and $g(x) = u^2(x,T)$ respectively, and $\mathcal L(\xi) = \rho_0$. Now, since for $i=1,2$,
\[
\E \left[u(0,X^{i,T}_0)\right] = \E \left[u(T,X^{i,T}_T)\right] - \lambda^i T + \E \left[\int_0^T  L(X^{i,T}_s,D_pH(X^{i,T}_s,Y^{i,T}_s,\rho^i_s),\rho^i_s)\dd s\right],
\]
we have that
\begin{align*}
\left|\lambda^1 - \lambda^2\right| & \le \frac{\left|\E \left[u(0,X^{1,T}_0)\right]\right| + \left|\E \left[u(T,X^{1,T}_T)\right]\right| + \left|\E \left[u(0,X^{2,T}_0)\right]\right| + \left|\E \left[u(T,X^{2,T}_T)\right]\right| }T \\
&+\frac1T \int_0^T \E\left[ \left|L(X^{1,T}_s,D_pH(X^{1,T}_s,Y^{1,T}_s,\rho^1_s),\rho^1_s) - L(X^{2,T}_s,D_pH(X^{2,T}_s,Y^{2,T}_s,\rho^2_s),\rho^2_s)\right|\right]\dd s.
\end{align*}
The first term of the right hand side of the previous inequality vanishes as $T \to \infty$ by uniform second moment bounds of Proposition \ref{prop:2nd_moment_bound}. Similarly, the integral on $(0,T)$ appearing in the second term remains bounded uniformly in $T$, by means of Theorem \ref{prop:pointwisedecay}(2) and Remark \ref{remL1L2}. Therefore, we conclude that $\lambda^1 = \lambda^2$ by letting $T \to \infty$.

Fix now an arbitrary $t>0$. To show that $\rho^1_t = \rho^2_t$, apply Corollary \ref{newcor}(2) and let $T \to \infty$. Similarly, Theorem \ref{prop:Du_loc} implies that $D_x u^1 \equiv D_x u^2$, which means that $u^1(t, \cdot)$ and $u^2(t, \cdot)$ differ by a constant $c(t)$ that may depend on time. Nevertheless, $u^i$ solve
\[
-\partial_t  u^i  -\b\Delta u^i  = -H(x,-D_{x}  u^1,\rho^1) - \l^1, \qquad {\rm{in}\ } (0,+\infty)\times \R^d
\]
hence, for any $T>0$, $u^1(T, \cdot)-u^2(T, \cdot)=c(T)$ forces $c(t) = u^1(t, \cdot)-u^2(t, \cdot)$ to coincide with $c(T)$ for any $t \le T$ (clearly if $\beta = 0$, and as a consequence of the maximum principle if $\beta > 0$). Therefore, $c(t)$ must be identically constant on $[0,+\infty)$.

\end{proof}


\appendix
\section{Uniform in time semi-concavity and convexity estimates}\label{sec:app}

It is well-known that the value function arising in finite time horizon stochastic or deterministic control problems is in general semi-concave, under suitable semi-concavity assumptions on the data. However, most classical proofs available in the literature (cf. \cite[Theorem 7.4.11]{CanSin}) provide semi-concavity constants that blow up when the time horizon tends to infinity. For our analysis in this paper it is crucial to obtain semi-concavity estimates which are uniform in $T$.

A great number of contributions in the literature address semi-concavity results for value functions. Besides \cite[Theorem 7.4.11]{CanSin}), we refer for instance to \cite{GigGotIshSat}, \cite{BuckCanQui}, \cite[Theorem 1.7]{CarPor20}, \cite[Theorem 5.9]{GomPimVos} or \cite[Theorem 3.11]{CGM}. However, none of these references address these estimates in the precise setting suitable for us, i.e. independently of the time horizon or under the same umbrella regarding first and second order problems. Therefore, for the sake of completeness we provide the desired semi-concavity results, that are uniform with respect to the time horizon $T$ and the noise intensity $\beta$. These results might be known for experts, but in lack of a precise reference, we have decided to include the details of the proof.

\begin{lemma}\label{lem:semi_concave} 
Let $L:(0,+\infty)\times\R^d\times\R^d\to\R$ be continuous, convex and super linear in the velocity variable, and be such that $\R^{2d}\ni(x,v)\mapsto L(t,x,v)$ is semi-concave with a constant $C_L\in\R$, uniformly in $t$ (i.e. $\sup_{t>0}D^2_{(x,v)} L(t,\cdot,\cdot)\le C_L I_{2d}$ in the sense of distributions). Suppose moreover that $g:\R^d\to\R$ is semi-concave with the constant $C_g\in\R$ (i.e. $D^2g\le C_gI_d$ in the sense of distributions). Let $T>0$. We define the value function $u:[0,T]\times\R^d\to\R$ in the classical way as 
\begin{equation*}
u(t,x):=\inf_\a\mathbb{E}\left\{\int_t^T L(s,X_s,\a_s)\dd s + g(X_T)\right\},
\end{equation*}
subject to 
\begin{equation}\label{eq:X}
X_s=x+\int_t^s\a_\t\dd\t+\sqrt{2\b}B^t_s,
\end{equation}
where $(t,x)\in[0,T]\times\R^d$, $\b\ge 0$ and $(B_\t)_{\t\in[0,T]}$ is a given Brownian motion, $B^t_s:=B_s-B_t$, $s\in[t,T]$ and the infimum is taken among all $\R^d$-valued $(\mathcal F_t)$-progressively measurable control processes $\a = (\a_t)_{t \ge 0}$. Then $u(t,\cdot)$ is semi-concave for all $t\in[0,T]$, with a semi-concavity constant depending only on $C_L$ and $C_g$, which is in particular independent of $T$ and $\b$.
\end{lemma}

\begin{proof}

By carefully following the proof of \cite[Theorem 7.4.11]{CanSin} (in the case of $\b=0$), our assumptions imply that that $u(t,\cdot)$ is semi-concave with a constant depending linearly on $(T-t)$. Now, we are going to show that this semi-concavity constant can in fact be chosen independent of $T$.

Let $\delta>0$ be fixed and let $t\in(0,T-\d)$ be fixed. Then the dynamic programming principle yields
\begin{equation*}
u(t,x)=\inf_\a \mathbb{E}\left\{\int_t^{t+\d} L(s,X_s,\a_s)\dd s + u(t+\delta,X_{t+\d})\right\}.
\end{equation*}
Let $x,y\in\R^d$, $\l\in[0,1]$, $\epsilon > 0$, and set $x_\l:=(1-\l)x+\l y$. Let now $(X_s)_{s\in(t,T)}$ and $(\a_s)_{s\in(t,T)}$ be an $\epsilon$-optimal state/control couple, that is $X_t=x_\l$ and
\begin{align*}
u(t,x_\l) \ge - \epsilon+\mathbb{E}\left\{\int_t^{t+\d} L(s,X_s,\a_s)\dd s + u(t+\d,X_{t+\d})\right\}.
\end{align*}
Let us consider the processes 
$$
X^x_s := (1-(s-t)/\d)x + x_\l(s-t)/\d  + \int_t^{s}\a_\t\dd\t + \sqrt{2\b}B^t_s= (1-(s-t)/\d)(x-x_\l) + X_s,
$$
and
$$
X^y_s := (1-(s-t)/\d)y + x_\l (s-t)/\d  + \int_t^{s}\a_\t\dd\t + \sqrt{2\b}B^t_s = (1-(s-t)/\d)(y-x_\l) + X_s.
$$
We notice that $X^x_t = x$, $X^y_t = y$ and $X^x_{t+\d} = X^y_{t+\d} = X_{t+\d}$.  We have furthermore that $(1-\l)X^x_s + \l X^y_s = X_s$, i.e. 
$$
(1-\l)X^x_s + \l X^y_s = x_\l + \int_t^s\a_\t\dd\t + \sqrt{2\b}B^t_s.
$$
Thus
\begin{align*}
&(1-\l)u(t,x) + \l u(t,y)\\
&\le (1-\l)\mathbb{E}\left\{\int_t^{t+\d} L(s,X^x_s,(x_\l-x)/\d +\a_s)\dd s+ u(t+\d,X^x_{t+\d})\right\}\\
&+\l \mathbb{E}\left\{\int_t^{t+\d} L(s,X^y_s,(x_\l-y)/\d +\a_s)\dd s+ u(t+\d,X^y_{t+\d})\right\}\\
&=\mathbb{E}\left\{\int_t^{t+\d} \left[(1-\l)L(s,X^x_s,(x_\l-x)/\d +\a_s) + \l L(s,X^y_s,(x_\l-y)/\d +\a_s)\right]\dd s + u(t+\d,x_{t+\d})\right\}\\
&\le \mathbb{E}\left\{\int_t^{t+\d} L(s, X_s,(1-\l)(x_\l-x)/\d+\l (x_\l-y)/\d+\a_s)\dd s\right\}\\ 
&+ \frac{C}{2}\l(1-\l)\mathbb{E}\left\{\int_{t}^{t+\d}\left[|X^x_s-X^y_s|^2 + \frac{1}{\d^2}|x-y|^2\right]\dd s + u(t+\d,x_{t+\d})\right\}\\
&=u(t,x_\l) +\epsilon + \frac{C}{2}\l(1-\l)\int_t^{t+\d} \left[(1-(s-t)/\d)^2+1/\d^2\right]|x-y|^2\dd s\\
&=u(t,x_\l) +\epsilon + \frac{C}{2}\l(1-\l)|x-y|^2\int_0^{1}\left[\d (1-r)^2+1/\d\right] \dd r = u(t,x_\l) +\epsilon  + \frac{C}{2}(\d/3+1/\d)\l(1-\l)|x-y|^2.
\end{align*}
Letting $\epsilon \to 0$, we conclude that $u(t,\cdot)$ is semi-concave, uniformly with respect to $t\in(0,T-\d)$, with a constant depending on the data and on $(\d/3+1/\d)$.

\medskip

Now, for $t\in[T-\d,T]$, we can use classical arguments to have a semi-concavity estimate with a constant that depends on the length of the interval, i.e. $\d$ in a linear way. Indeed, for $x,y\in\R^d$, as above, we consider $x_\l=(1-\l)x+\l y$ and  $(X_s)_{s\in(t,T)}$ and $(\a_s)_{s\in(t,T)}$ as in \eqref{eq:X} $\epsilon$-optimal for $u(t,x)$. Now, we have
\begin{align*}
&(1-\l)u(t,x) + \l u(t,y)\\
&\le (1-\l)\mathbb{E}\left\{\int_t^{T} L(s,X_s + x-x_\l, \a_s)\dd s+ g(X_{T}+x-x_\l)\right\}\\
&+\l \mathbb{E}\left\{\int_t^{T} L(s,X_s+y-x_\l, \a_s)\dd s+ g(X_{T}+y-x_\l)\right\}\\
&\le \mathbb{E}\left\{\int_t^{T} \left[L(s,X_s, \a_s)+\frac{C}{2}\l(1-\l)|x-y|^2\right]\dd s+ g(X_{T})+\frac{C}{2}\l(1-\l)|x-y|^2\right\}\\
&\le u(t,x_\l) + (1+\d)\frac{C}{2}\l(1-\l)|x-y|^2 +\epsilon,
\end{align*}
which indeed proves our claim.
\medskip

Combining the two previous arguments for $t\in(0,T-\d)$ and $t\in[T-\d,T]$, respectively, we can choose a small (universal) constant $\d$, such that $u(t,\cdot)$ is semi-concave with a constant depending on the data (but independent of $\b$ and $T$).
\end{proof}

It is also well-known that for fully convex control problems, the value function inherits this convexity. We recall this result here (see for instance \cite[Lemma 3.4]{MesMou} or \cite{BanMesMou} on this matter).

\begin{lemma}\label{lem:appdx}
Suppose that we are precisely in the setting of Lemma \ref{lem:semi_concave}. Suppose that the functions $\R^d\times\R^d\ni(x,v)\mapsto L(t,x,v)$ and $\R^d\ni x\mapsto g(x)$ are convex. Then $\R^d\ni x\mapsto u(t,x)$ is convex for all $t\in [0,T]$.

This, together with the implications of Lemma \ref{lem:semi_concave} implies that there exists a constant $C>0$ depending only on $C_L$ and $C_g$ (but independent of $T$ or $\beta$) such that
\begin{align*}
 \esssup_{(t,x)\in[0,T] \times\R^d} |D^2_{xx}u(t,x)|\le C.
\end{align*}
\end{lemma}

\begin{remark}
The previously mentioned two results are used for the solution $u$ to the Hamilton--Jacobi--Bellman equation appearing in \eqref{eq:MFG_intro}, when the solution $(\rho_{t})_{t\in[0,T]}$ to the Kolmogorov--Fokker--Planck equation is fixed and taken as an input in the HJB equation. Here the Lagrangian $L:\R^{d}\times\R^{d}\times\sP_{2}(\R^{d})\to\R$ and Hamiltonian $H:\R^{d}\times\R^{d}\times\sP_{2}(\R^{d})\to\R$ are Legendre duals of each other with respect to the second variable, i.e.
$$
L(x,v,\mu) = \sup_{p\in\R^{d}}\left\{v\cdot p - H(x,p,\mu)\right\}.
$$
\end{remark}

\section{It\^o Lemmas}\label{app:ito}

Based on \cite[Proposition 5.102]{CD1}, we will make use of the following It\^o lemma. 
\begin{lemma}\label{lem:Ito_CD}
Let $F:\R^{d}\times\R^{d}\times\sP_{2}(\R^{d})\to\R$ be a fully {$C^{2}$} function with the property that for any compact set $K\subset \R^{d}\times\R^d\times\sP_{2}(\R^{d})$ we have
$$
\sup_{(x,y,\mu)\in K}\left\{\int_{\R^{d}}|D_{\mu}F(x,y,\mu,\tilde x)|^{2}\dd\mu(\tilde x)+\int_{\R^{d}}|D^{2}_{\mu\tilde x}F(x,y,\mu,\tilde x)|^{2}\dd \mu(\tilde x)\right\}<+\infty. 
$$ 
Let moreover $(\rho_s)_{s\in[t,T]}$ be a flow in $\sP_{2}(\R^{d})$ such that $\rho_{s} = \cL(R_{s})$, where $(R_{s})_{s\in[t,T]}$ solves
\begin{equation}\label{eq:auxiliary}
\dd R_{s} = \alpha_{s}\dd s + \sqrt{2\beta}\dd B_{s},
\end{equation}
with $(\a_s)_{s\in[t,T]}$ a given process.

Let $(X_s,Y_s,Z_s)_{s\in[t,T]}$ be the solution to \eqref{eq:FBSDE}, where $(\rho_{s})_{s\in[t,T]}$ is given as above. Then, we have
\begin{align*}
\dd F(X_{s},Y_{s},\rho_{s}) & = \left[D_{x}F(X_{s},Y_{s},\rho_{s})\cdot D_{p}H(X_{s},Y_{s},\rho_{s}) - D_{y}F(X_{s},Y_{s},\rho_{s})\cdot D_{x}H(X_{s},Y_{s},\rho_{s})\right]\dd s\\
& + \sqrt{2\beta}\left[  D_{x}F(X_{s},Y_{s},\rho_{s}) +  D_{y}F(X_{s},Y_{s},\rho_{s})Z_{t}\right]\cdot \dd B_{s}\\
& +\beta\left\{\Delta_{x} F(X_{s},Y_{s},\rho_{s}) + {{\rm{trace}}\left[D^{2}_{yy}F(X_{s},Y_{s},\rho_{s})Z_{s}Z_s^\top\right]}\right\}\dd s\\
&+ \beta\left\{{2}{\rm{trace}}\left[D^{2}_{xy}F(X_{s},Y_{s},\rho_{s})Z_{s}\right]\right\}\dd s\\
& + \tilde\E \left[D_{\mu}F(X_{s},Y_{s},\rho_{s},\tilde R_{s})\cdot \tilde\alpha_{s}\right]\dd s\\
& + \beta \tilde\E \left\{ {\rm{trace}}\left[D^{2}_{\mu \tilde x}F(X_{s},Y_{s},\rho_{s},\tilde R_{s})\right]\right\}\dd s,
\end{align*}
where the process $(\tilde R_{s},\tilde\a_s)_{s\in[t,T]}$ is an independent copy of of the process $(R_{s},\a_s)_{s\in[t,T]}$ on a copy $(\tilde\Omega,\tilde\cF,\tilde\P)$ of the probability space $(\Omega,\cF,\P)$.
\end{lemma}

Based on this lemma, we have the following results.

\begin{lemma}\label{lem:ito_XD_pH}
Let $H:\R^d\times\R^d\times\sP_2(\R^d)\to\R$ satisfy the assumptions in \eqref{hyp:H_reg}, let $(X_s,Y_s,Z_s)_{s\in[t,T]}$ be a solution to \eqref{eq:FBSDE} and let $(\rho_s)_{s\in[t,T]}$ be a flow in $\sP_{2}(\R^{d})$ such that $\rho_{s} = \cL(R_{s})$, where $R_s$ satisfies \eqref{eq:auxiliary} for some $(\a_s)_{s\in[t,T]}$ given process. Then, we have
\begin{align*}
\dd [X_s&\cdot D_pH(X_s,Y_s,\rho_s)] = \left[|D_{p}H(X_{s},Y_{s},\rho_{s})|^2 +(D_{px}^2H(X_{s},Y_{s},\rho_{s})X_s)\cdot D_{p}H(X_{s},Y_{s},\rho_{s}) \right]\dd s\\
& - \left[ (D_{pp}^2H(X_{s},Y_{s},\rho_{s})X_s)\cdot D_{x}H(X_{s},Y_{s},\rho_{s})\right]\dd s\\
& + \sqrt{2\beta} \left[  D_{p}H(X_{s},Y_{s},\rho_{s}) +D_{px}^2H(X_{s},Y_{s},\rho_{s})X_s +  D^2_{pp}H(X_{s},Y_{s},\rho_{s})Z_{s}X_s \right]\cdot \dd B_{s}\\
& +\beta \left\{ D_p{\rm{trace}}[D^2_{xx}H(X_s,Y_s,\rho_s)]\cdot X_s + 2{\rm{trace}}[D^2_{px}H(X_s,Y_s,\rho_s)] \right\}\dd s\\
& +\beta \left\{ {\rm{trace}}{\left[(D^3_{ppp}H(X_{s},Y_{s},\rho_{s}) X_s)Z_sZ_s^\top\right]}\right\}\dd s\\
& +{2}\beta\, {\rm{trace}}\left\{\left[D^{3}_{ppx}H(X_{s},Y_{s},\rho_{s})X_s+D^2_{pp}H(X_s,Y_s,\rho_s)\right]Z_{s}\right\}\dd s\\
& + \tilde\E \left[(D_{p\mu}H(X_{s},Y_{s},\rho_{s},\tilde R_{s})X_s)\cdot \tilde\alpha_{s}\right]\dd s\\
& + \beta \tilde\E \left\{ {\rm{trace}}\left[D^{3}_{p\mu \tilde x}H(X_{s},Y_{s},\rho_{s},\tilde R_{s})X_s\right]\right\}\dd s.
\end{align*}
\end{lemma}

\begin{proof}
We will simply apply Lemma \ref{lem:Ito_CD} for $F:\R^d\times\R^d\times\sP_2(\R^d)\to\R$ defined as 
$$F(x,y,\rho) = x\cdot D_pH(x,y,\rho).$$
\end{proof} 

\begin{lemma}\label{lem:ito_YD_xH}
Let $H:\R^d\times\R^d\times\sP_2(\R^d)\to\R$ satisfy the assumptions in \eqref{hyp:H_reg}, let $(X_s,Y_s,Z_s)_{s\in[t,T]}$ be a solution to \eqref{eq:FBSDE} and let $(\rho_s)_{s\in[t,T]}$ be a flow in $\sP_{2}(\R^{d})$ such that $\rho_{s} = \cL(R_{s})$, where $R_s$ satisfies \eqref{eq:auxiliary} for some $(\a_s)_{s\in[t,T]}$ given process. Then, we have
\begin{align*}
\dd [Y_s&\cdot D_xH(X_s,Y_s,\rho_s)] = \left[(D^{2}_{xx}H(X_{s},Y_{s},\rho_{s})Y_{s})\cdot D_{p}H(X_{s},Y_{s},\rho_{s}) -|D_{x}H(X_{s},Y_{s},\rho_{s})|^{2} \right]\dd s\\
& - \left[ (D_{xp}^2H(X_{s},Y_{s},\rho_{s})Y_s)\cdot D_{x}H(X_{s},Y_{s},\rho_{s})\right]\dd s\\
& + \sqrt{2\beta} \left\{  D^{2}_{xx}H(X_{s},Y_{s},\rho_{s})Y_{s} +Z_{s}\left[D_{x}H(X_{s},Y_{s},\rho_{s}) +  D^2_{xp}H(X_{s},Y_{s},\rho_{s})Y_{s}\right] \right\}\cdot \dd B_{s}\\
& +\beta \left\{ D_x{\rm{trace}}\left[D^2_{xx}H(X_{s},Y_{s},\rho_{s}\right]\cdot Y_s\right\}\dd s\\
& +\beta \left\{ { {\rm{trace}}\left[\left(D^2_{xpp}H(X_s,Y_s,\rho_s)Y_s + 2D^2_{px}H(X_s,Y_s,\rho_s)\right)Z_sZ_s^\top\right] }\right\}\dd s\\
& +{2}\beta\, {\rm{trace}}\left\{\left[D^{3}_{xxp}H(X_{s},Y_{s},\rho_{s})Y_s+D^2_{xx}H(X_s,Y_s,\rho_s)\right]Z_{s}\right\}\dd s\\
& + \tilde\E \left[D^2_{x\mu}H(X_{s},Y_{s},\rho_{s},\tilde R_{s})Y_s\cdot \tilde\alpha_{s}\right]\dd s\\
& + \beta \tilde\E \left\{ {\rm{trace}}\left[D^{3}_{x\mu \tilde x}H(X_{s},Y_{s},\rho_{s},\tilde R_{s})Y_s\right]\right\}\dd s.
\end{align*}
\end{lemma}

\begin{proof}
We will simply apply Lemma \ref{lem:Ito_CD} for $F:\R^d\times\R^d\times\sP_2(\R^d)\to\R$ defined as 
$$F(x,y,\rho) = y\cdot D_xH(x,y,\rho).$$
\end{proof} 

\begin{remark}
In the statements of the previous lemmas we have abused the notation $D^{3}_{xxp}H$ (also in the case of other similar terms involving third derivatives). In particular in the products involving such terms we denote the action of the tensor $D^{3}_{xxp}H$ in $\R^{d\times d\times d}$ on the vector $Y$ in $\R^{d}$. The resulting object is a standard square matrix in $\R^{d\times d}$. For instance, $\left(D^{3}_{xxp}HY\right)_{ij} = \sum_{k=1}^{d}\partial^{3}_{x_{i}x_{j}p_{k}}HY_{k}$.
\end{remark}

\medskip
\medskip

\noindent {\large\bf Declarations.}

\medskip
\medskip

{\bf Conflict of Interest.} The authors declare no conflict of interest. 

\medskip

{\bf Data availability.} No datasets were generated or analyzed during the current study.

\def\cprime{$'$}

\end{document}